\documentclass[12pt,a4paper]{article}
 
 \usepackage{amssymb, graphics, setspace}
\usepackage{romanbar}
\usepackage{bm}
\usepackage{subfig}
\newcommand{\mathsym}[1]{{}}
\newcommand{\unicode}[1]{{}}

\usepackage[pdftex]{hyperref}
\usepackage{amssymb}
\usepackage[centertags]{amsmath}

\usepackage{amsthm}
\usepackage{graphicx}
\graphicspath{ {./Pictures/} }

\usepackage{xcolor}

\newtheorem{pro}{Proposition}
\newtheorem*{pro*}{Proposition}
\newtheorem{lem}{Lemma}
\newtheorem{rem}{Remark}
\newtheorem*{rem*}{Remark}

\newtheorem{cl}{Claim}

\renewcommand{\proof}{\textsc{proof. }}
\newcommand{\eproof}{\hfill $\Box$}
\newcommand{\bal}{{\bm \alpha}}
\newcommand{\bbe}{{\bm \beta}}

\newcommand{\bphi}{{\bm \varphi}}

\newcommand{\bI}{\mathbf I}

\newcommand{\Ju}{{J_2}}
\newcommand{\Juq}{{J_2^2}}
\newcommand{\II}{{\rm I}}
\newcommand{\cE}{\mathcal E}
\newcommand{\rp}{{\rm p}}
\renewcommand{\rq}{{\rm q}}
\newcommand{\rr}{{\rm r}}
\newcommand{\Emax}{E_{\rm max}}
\newcommand{\Emin}{E_{\rm min}}
\newcommand{\Esad}{E_{\rm sad}}
\newcommand{\PP}{{\mathcal P}}
\newcommand{\Ac}{{\mathcal A}}
\newcommand{\rmuno}{{\rm I}}
\newcommand{\rmdue}{{\rm I\!I}}
\newcommand{\rmtre}{{\rm I\!I\!I}}
\newcommand{\rmquattro}{{\rm I\!V}}
\newcommand{\At}{{\mathtt A}}
\newcommand{\Bt}{{\mathtt B}}
\newcommand{\Ct}{{\mathtt C}}
\newcommand{\Dt}{{\mathtt D}}
\newcommand{\Tt}{{\mathtt T}}
\newcommand{\kt}{{\mathtt k}}
\newcommand{\ct}{{\mathtt c}}

\newcommand{\at}{{\mathtt a}}
\newcommand{\bt}{{\mathtt b}}
 
\makeindex

\title{Three-to-one internal resonances in coupled harmonic oscillators with cubic nonlinearity}

\author{L. Di Gregorio, W. Lacarbonara}

\begin{document}

\allowdisplaybreaks
\maketitle

{\abstract{We investigate a general system of two coupled harmonic oscillators with cubic nonlinearity, a model relevant to various structural engineering applications. As a concrete example, we consider the case of  two oscillators  obtained from the reduction of the wave propagation equations representing a cellular hosting structure with 1-dof resonators in each cell.
Without damping, the system is Hamiltonian, with the origin as an elliptic equilibrium characterized by two distinct linear frequencies. To understand the dynamics, it is crucial to derive explicit analytic formulae for the nonlinear frequencies as functions of the physical parameters involved.
In the small amplitude regime (perturbative case), we provide the first-order nonlinear correction to the linear frequencies. While this analytic expression was already derived for non-resonant cases, it is novel in the context of resonant or nearly resonant scenarios. Specifically, we focus on the 3:1 resonance, the only resonance involved in the first-order correction. Utilizing the Hamiltonian structure, we employ Perturbation Theory methods to transform the system into Birkhoff Normal Form up to order four. This involves converting the system into action-angle variables (symplectically rescaled polar coordinates), where the truncated Hamiltonian at order four depends on the actions and, due to the resonance, on one ``slow" angle.
By constructing suitable nonlinear and not close-to-the-identity coordinate transformations, we identify new sets of symplectic action-angle variables. In these variables, the resulting system is integrable up to higher-order terms, meaning it does not depend on the angles, and the frequencies are obtained from the derivatives of the energy with respect to the actions. This construction is highly dependent on the physical parameters, necessitating a detailed case analysis of the phase portrait, revealing up to six topologically distinct behaviors. In each configuration, we describe the nonlinear normal modes (elliptic/hyperbolic periodic orbits, invariant tori) and their stable and unstable manifolds of the truncated Hamiltonian.
As an application, we examine wave propagation in metamaterial honeycombs with periodically distributed nonlinear resonators, evaluating the nonlinear effects on the bandgap particularly in the presence of resonances.
 }}

 \bigskip
\noindent
{\footnotesize{\textbf{Acknowledgments} Project ECS 0000024 Rome Technopole, CUP  B83C22002820006, National Recovery and Resilience Plan (NRRP) Mission 4 Component 2 Investment 1.5, funded by the European Union - NextGenerationEU.}}

\bigskip

\noindent
{\footnotesize{\textbf{Funder} Project funded under the National Recovery and Resilience Plan (NRRP), Mission 4 Component 2 Investment 1.5 - Call for tender No. 3277 of 30 December 2021 of the Italian Ministry of University and Research funded by the European Union - NextGenerationEU.}}
\tableofcontents

\section{Introduction}

\begin{figure}[h!]
\center
\includegraphics[width=7cm,keepaspectratio]{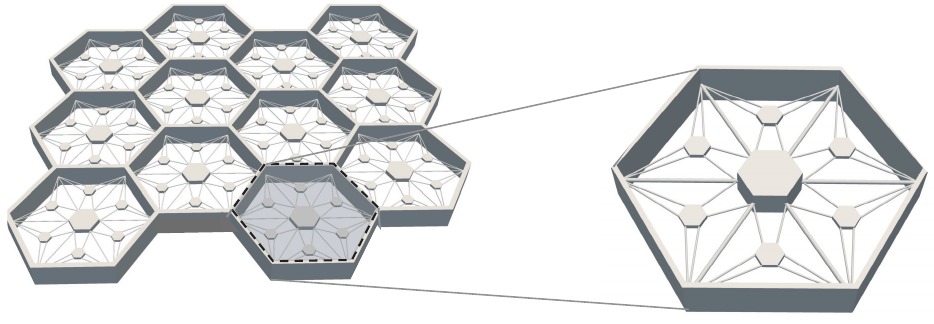}
\caption{Schematic view of the orthotropic plate model with the periodically distributed spider-web resonators, see \cite{SW23jsv} as reference.}
\label{Resonators_model}
\end{figure}

Let us briefly recall the model introduced in 
\cite{SW23jsv}.
Figure \ref{Resonators_model}
 shows  schematic view of the orthotropic plate model with the periodically distributed spider-web resonators. 
 Each multi-frequency resonator should be meant as the multi-mass-spring system resulting from the multi-dof modal reduction of the infinite-dimensional resonator (i.e., the spider webs with a central mass, here represented in the figure, for the sake of graphical clarity, by a single mass-spring system instead of a set of mass-spring systems). The modal reduction is performed via the Galerkin projection method employing a  number of mode shapes of the distributed-parameter resonators. Each resonator is represented by equivalent modal masses and modal springs.

The adopted plate theory  (see \cite{W}) with the elastic constants of the equivalent, homogenized  orthotropic material  describes the motion of the honeycomb with the attached resonators.
By  the Floquet-Bloch Theorem, which states that the solutions of the {corresponding linear} periodic  resonators-plate system are quasi-periodic in space with the fundamental periodicity provided by the lattice period,  
the plate equation of motion can be projected onto the unit cell domain (i.e., the periodically repeated lattice unit).
Then one obtains a system of $2 N$ 
coupled second order 
ODEs, $N$ being
 the number of retained resonators
modes.
For the metamaterial lattice with an array of equally spaced single-dof resonators, i.e., N = 1,  equations 
reduce 
 to the following
 system of second order ODEs
\begin{equation}\label{firenze}
 \left(
\begin{array}{cc}
	\tilde{M}_H(\tilde{k}_1,\tilde{k}_2)&\tilde{M}\\
       \tilde{M}& \tilde{M}
	\end{array}	
\right)
\left(
		\begin{array}{c}
			\ddot{\tilde w}_0 \\
			\ddot{\tilde z}_0
		\end{array}
		\right)
		+
		 \left(
\begin{array}{cc}
	\tilde{K}_H(\tilde{k}_1,\tilde{k}_2)&0\\
       0& \tilde{K}
	\end{array}	
\right)
		\left(
		\begin{array}{c}
			{\tilde w}_0 \\
			{\tilde z}_0
		\end{array}
		\right)
		=
		-\left(
		\begin{array}{c}
			 0 \\
			\tilde N^{(3)} {\tilde z}_0^3
		\end{array}
		\right)\,,
\end{equation}
where
${\tilde w}_0$ and 
			${\tilde z}_0$
denote the nondimensional plate deflection and  resonator relative motion at the origin of the fixed frame; 
	\begin{equation}\label{coefficients3}
		\tilde{M}_H (\tilde{k}_1,\tilde{k}_2) 
		:=
		 \frac{4 \sqrt{3} \sin \left(\frac{\tilde{k}_1}{2}\right) \sin \left(\frac{1}{4}  \left(\tilde{k}_1+\sqrt{3} \tilde{k}_2\right)\right)}{\tilde{k}_1 \left(\tilde{k}_1+\sqrt{3} \tilde{k}_2\right)} 
	\end{equation}
	and 
\begin{equation}\label{coefficients2}
	\begin{split}
		&
\tilde{K}_H(\tilde{k}_1,\tilde{k}_2)= 
\tilde{K}_H(\tilde{k}_1,\tilde{k}_2;
\tilde{D}_{12},\tilde{D}_{66},\tilde{D}_{22})
:=\tilde{M}_H(\tilde{k}_1,\tilde{k}_2) \left[\tilde{k}_1^4  +2 \tilde{k}_1^2 \tilde{k}_2^2 (\tilde{D}_{12} + 2\tilde{D}_{66}) + \tilde{k}_2^4 \tilde{D}_{22} \right]
	\end{split}
\end{equation}
are the nondimensional modal mass and stiffness
as functions of the nondimensional
wave numbers $(\tilde{k}_1,\tilde{k}_2)$,
which stay within the irreducible Brillouin triangle
$\triangle$ (see Figure \ref{Brillouin}):
	\begin{center}
		\begin{figure}[h!]
			\centering
\includegraphics[width=6cm,height=6cm,keepaspectratio]{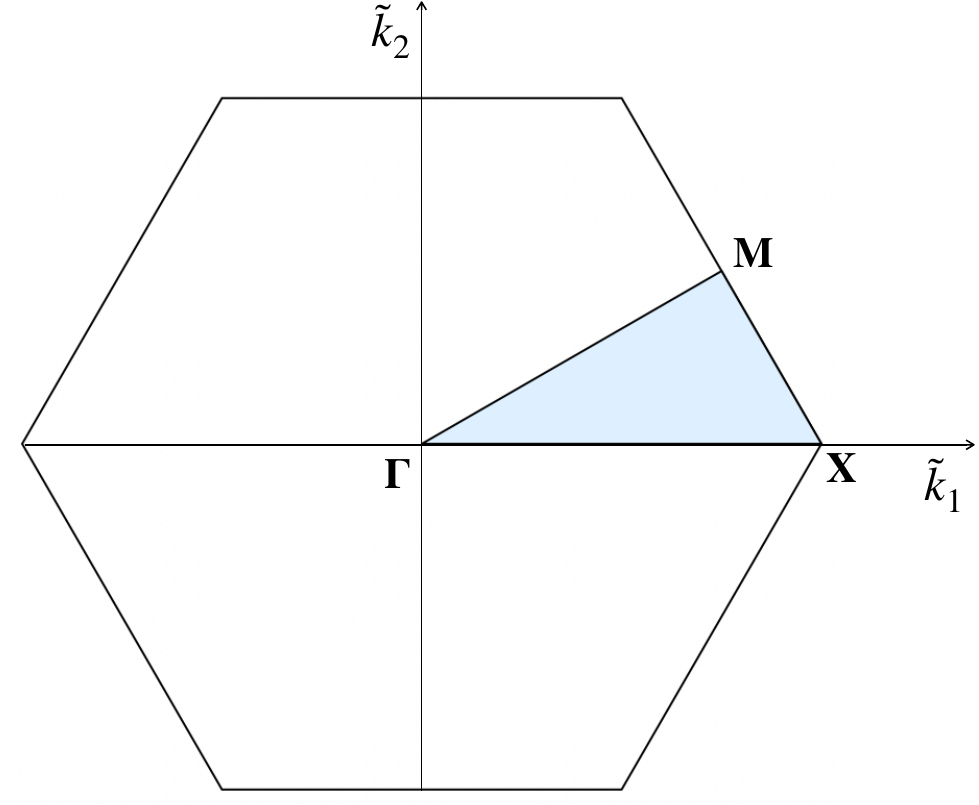}
\caption{The irreducible Brillouin triangle 
$\triangle:=\bf{\Gamma}\overset{{\triangle}}{\bf{X}} \bf{M}$. $\mathbf{\Gamma}=(0,0)$, $\mathbf{X}=(\frac{4}{3}\pi,0)$, $\mathbf{M}=(\pi,\frac{\pi}{\sqrt{3}})$.}
\label{Brillouin}
		\end{figure}
	\end{center}
moreover	
$$
\tilde D_{12}=0.0815599,\qquad \tilde D_{22}=12.48, 
\qquad\tilde D_{66}=0.0000247357\,,
$$
 are the nondimensional plate bending stiffness coefficients; finally
 $\tilde{N}^{(3)}$ is the nondimensional 
 nonlinearity.

\medskip
Actually we  consider the more
 general system of ODEs
\begin{equation}\label{autostrada}
\mathtt M\left(
		\begin{array}{c}
			\ddot v \\
			\ddot y
		\end{array}
		\right)
		+\mathtt K\left(
		\begin{array}{c}
			v \\
			y
		\end{array}
		\right)
		=
		-\left(
		\begin{array}{c}
			 M_3v^3 \\
			N_3 y^3
		\end{array}
		\right)\,,
\end{equation}
where $v(t),y(t)$
are unknown scalar functions,
$M_3,N_3$ are real coefficients,
$\mathtt M$ is a symmetric 
positive definite $2\times 2 $ real matrix
and $\mathtt K$ is a diagonal 
positive definite $2\times 2 $ real matrix.

Note that 
 \eqref{firenze} is a particular
 case of \eqref{autostrada} taking
 $v=\tilde{w}_0$, $z=\tilde{z}_{0}$,
 $M_3=0,$ $N_3=\tilde{N}^{(3)}$
 and
 \begin{equation}\label{frittata}
 \mathtt M=
 \left(
\begin{array}{cc}
	\tilde{M}_H(\tilde{k}_1,\tilde{k}_2)&\tilde{M}\\
       \tilde{M}& \tilde{M}
	\end{array}	
\right)\,,\qquad
 \mathtt K=
 \left(
\begin{array}{cc}
	\tilde{K}_H(\tilde{k}_1,\tilde{k}_2)&0\\
       0& \tilde{K}
	\end{array}	
\right)\,,
 \end{equation}
 with $\tilde{M}_H(\tilde{k}_1,\tilde{k}_2)$
 and 
 $\tilde{K}_H(\tilde{k}_1,\tilde{k}_2)$
 defined in 
 \eqref{coefficients3} and \eqref{coefficients2},
 respectively.

\bigskip

The existing literature on Hamiltonian and dissipative systems covers various topics, including bifurcations, invariant manifolds, and homoclinic and heteroclinic orbits. In \cite{Fontich23}, the authors study a one-parameter family of 2-DOF Hamiltonian systems with an equilibrium point undergoing a Hamiltonian-Hopf bifurcation. They focus on invariant manifolds and the behavior of the splitting of 2D invariant manifolds in the presence of homoclinic orbits. Similarly, \cite{Celletti13} presents a KAM theory for conformally symplectic dissipative systems, demonstrating that solutions with a fixed 
n-dimensional (Diophantine) frequency can be found by
an a-posteriori approach adjusting the parameters.

In \cite{Llave06}, the authors develop numerical algorithms to compute invariant manifolds in quasi-periodically forced systems, focusing on invariant tori and their asymptotic invariant manifolds (whiskers). These algorithms utilize Newton's method and power-matching expansions of parameterizations. \cite{Cabre05} describes a method to establish the existence and regularity of invariant manifolds, simplifying the proof of the stable manifold theorem near hyperbolic points by using the implicit function theorem in Banach spaces.

\cite{H16} proposes a unified approach to nonlinear modal analysis in dissipative oscillatory systems. This approach defines nonlinear normal modes (NNMs) and spectral submanifolds, emphasizing the importance of damping for accurate conclusions about them, and the reduced-order models they produce.
Lastly, \cite{HW95}, \cite{HW96} and \cite{HW93} develop methods to detect orbits asymptotic to slow manifolds in perturbed Hamiltonian systems, revealing complex chaotic behaviors and the creation of homoclinic orbits in resonant Hamiltonian systems through geometric singular perturbation theory and Melnikov-type methods.

\subsection{Main results}

We are interested here in small
amplitude solutions
of \eqref{autostrada}. In the first
approximation the system is linear
with  linear frequencies 
$\omega_-$ and $\omega_+$
and the nonlinearity is a third order
perturbation.
If the linear frequencies 
are non vanishing, distinct
and satisfy the non resonance 
condition 
$3\omega_-\neq \omega_+$
the system can be integrated,
for instance, using the multiple scales method,
up to a smaller  fifth order nonlinear remainder, see \cite{SW23jsv}.
In particular, \cite{SW23jsv} provides explicit expressions for the nonlinear frequencies
of the truncated system (obtained
by disregarding the fifth order perturbation) 
 as functions of the
initial amplitudes. Moreover
the effects on the bandgap
were explored.
\\
In \cite{DL},
we analytically estimated the applicability threshold of the perturbative argument,
specifically the maximal admissible amplitude for which the  above formula
is valid.
It was found that this  applicability  threshold decays to zero in the presence of resonances,
more precisely when the ratio between the optical and acoustic  frequencies is
close to 3; indeed the 3:1 resonance  is the only involved resonance 
in the first order correction. 
\\
The  methodology used is based on techniques from   Hamiltonian Perturbation Theory.  
Since the system is conservative,
 we study it as a Hamiltonian 
 system. The origin is an elliptic
 equilibrium and we put the 
 system in (complete) {\sl Birkhoff Normal Form}
 up to order 4 (3 in the equations
 of motion). 
 The Birkhoff Normal Form
 is a powerful tool in Hamiltonian
 Perturbation Theory
 that, through a suitable symplectic, 
 close-to-the-identity
 nonlinear change of coordinates,
simplifies the Hamiltonian. More precisely, after introducing
 action-angle variables\footnote{Essentially rescaled polar coordinates.},
 in the non resonant case,
 the truncated system at order four
 is \textsl{integrated}, meaning
 its Hamiltonian depends 
 only on the actions,
 which are constant of motion,
 and not on the angles.
 As a consequence the phase space of the truncated Hamiltonian
 is completely foliated by 
 nonlinear normal modes (NNMs), which are 
 two dimensional invariant tori filled with 
 periodic/quasi-periodic orbits 
 depending on whether the frequency ratio is rational/irrational.
 Moreover such tori are (constant) graphs over the angles.
 Finally the nonlinear frequencies
 of the truncated Hamiltonian
  are easily evaluated
  as the derivatives of the Hamiltonian,
  i.e. the energy, with respect
  to the two actions.
  This procedure, being perturbative in nature,
only works in a ball of small radius $\varepsilon$ 
around the origin. More precisely in \cite{DL} we proved that there exists a  constant $c_1$, which
was explicitly estimated as function of the physical parameters, such that the smallness condition
reads 
\begin{equation}\label{def:sigma} 
\varepsilon \leq c_1 \sqrt{|\sigma|}\,,
\qquad \mbox{where}\qquad
\sigma:=\omega_+ - 3 \omega_-\,.
\end{equation}

\begin{figure}[h!]
\begin{minipage}{.45\textwidth}
\center
\includegraphics[width=7.4cm, keepaspectratio]{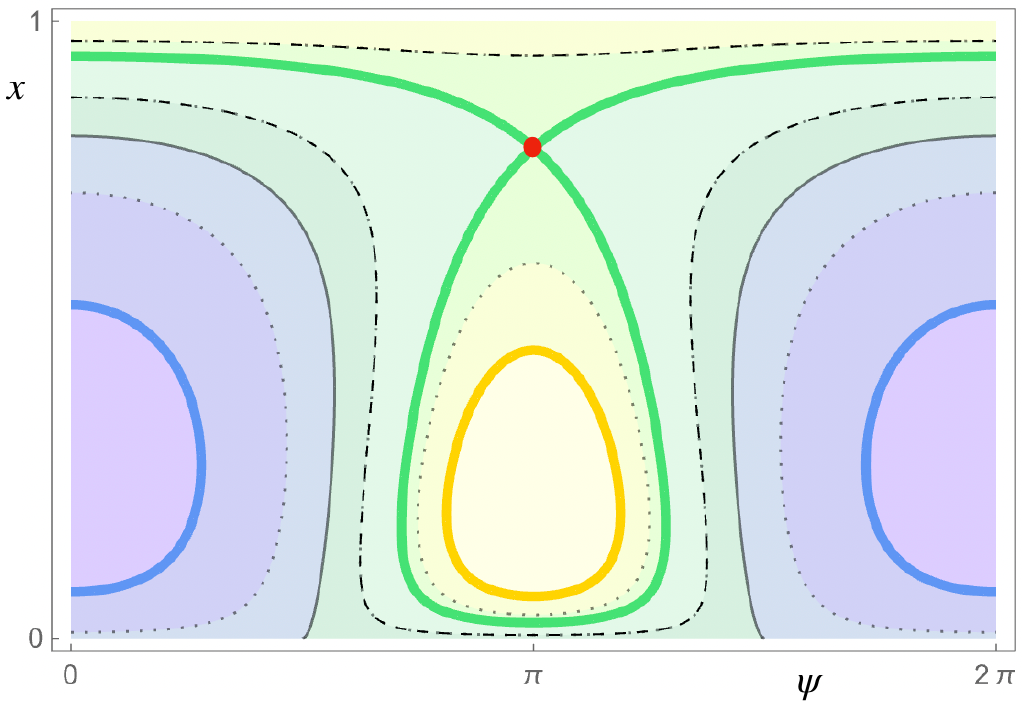}\\
\includegraphics[width=7.5cm, keepaspectratio]{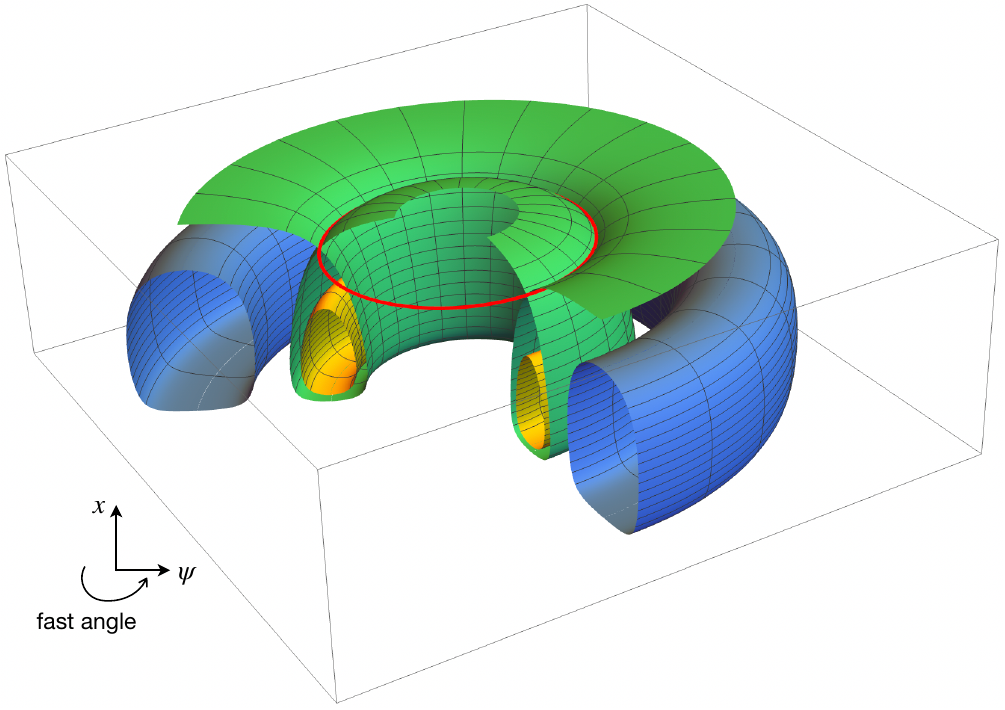}
\end{minipage}
\begin{minipage}{.42\textwidth}
\includegraphics[width=8.3cm,keepaspectratio]{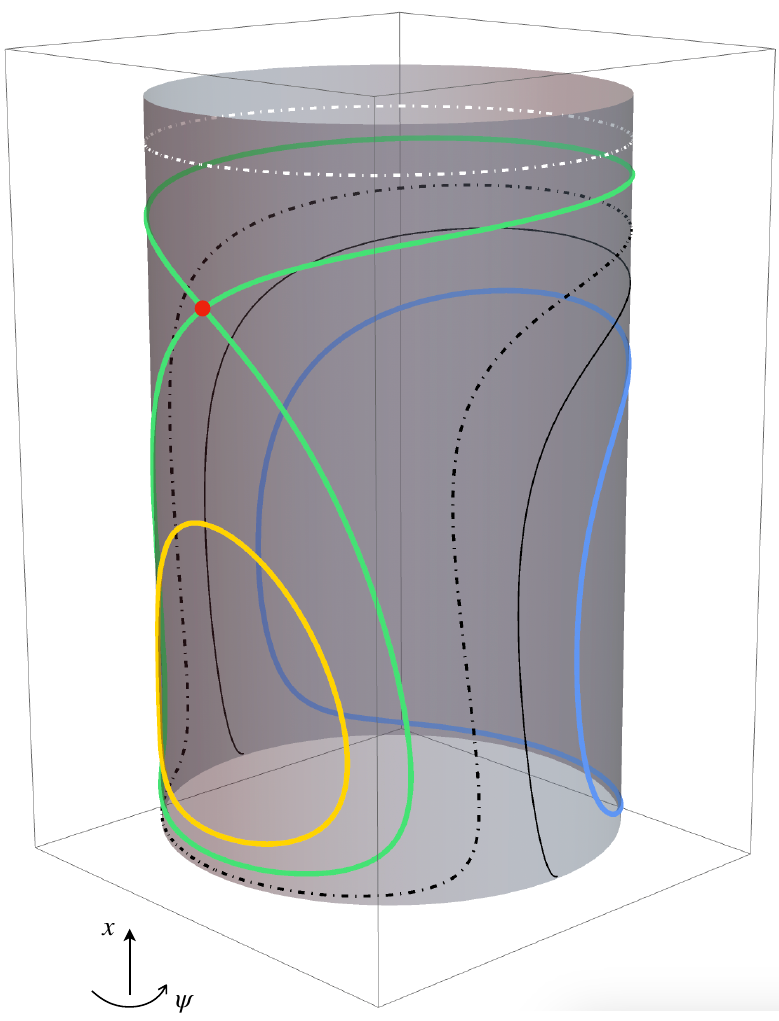}
\end{minipage}
\caption{(Left top) Level curves of the phase space of
the reduced Hamiltonian obtained 
by the fourth order resonant  
$\hat{\mathbb H}_{\rm res}$ (see \eqref{secularRES})
fixing the constant of motion $\Ju=10^{-4}$
(here, e.g., we have chosen the physical parameters as follows: $\tilde{k}_1=4\pi/3, \tilde{k}_2=0, \tilde{M}=0.005862, \tilde{K}=1.73, M_3=0, N_3=-10^4$).
The slow angle is on the horizontal axis and 
its (rescaled) conjugated action is on the vertical axis,
so that the phase space is actually the cylinder shown on the right.
Up to the green and  black curves that 
act as separatrices, the phase
space is divided into four connected components.
Every component is completely foliated  by one dimensional NNMs 
(periodic orbits).
Such NNMs have different topology:
the orbits in the zone above the green curve
or between the green and solid black curves
wrap around the cylinder (dash-dotted curves);
the orbits inside the green or the black curves
do not wrap around the cylinder and are contractible
(dotted, blue and yellow curves).
The red point and the green curve are, respectively,
a zero dimensional NNM (a hyperbolic equilibrium) 
and its one dimensional  (coinciding) 
stable and unstable separatrices.
(Right) The  cylindrical phase portrait
 immersed in the three dimensional space.
\\
(Left bottom) 
A representation
 of the phase space
of the truncated Hamiltonian $\hat{\mathbb H}_{\rm res}$,
once we have fixed the constant of motion $\Ju=10^{-4}$.
 The image is 
obtained by rotating the picture on the top by  the fast angle   from $0$ to $4\pi/3$.
In particular, by rotation,
 the blue and yellow curves become 
two dimensional NNMs (invariant tori) and the red point and the green curve become, respectively,
a one  dimensional NNM (a hyperbolic periodic orbit) 
and its two dimensional (coinciding) 
stable and unstable manifolds. 
}
\label{scacchi}
\end{figure}

In contrast, the main aim of
 the present paper is to investigate
what happens in the complementary regime,
namely 
when the linear frequencies 
are in, or almost in,  3:1 resonance, specifically when
\begin{equation}\label{def:sigma2} 
c_1 \sqrt{|\sigma|}<\varepsilon
\end{equation}
and $\varepsilon$ is small enough.
In this case, only a resonant BNF is available.
This means that, after introducing action-angle variables
and a linear symplectic change of coordinates,
the truncated Hamiltonian at order four, 
$\hat{\mathbb H}_{\rm res}$ (see \eqref{secularRES}), depends on the actions and on  one ``slow'' angle (as its associated  
frequency is small).
The phase  portrait becomes more complicated
and interesting; its topology strongly depends 
on the values of the physical parameters.
The phase space is still  foliated by two dimensional 
NNMs (invariant tori) but 
many of them are no longer graphs over the angles
as in the nonresonant case, exhibiting different topologies.
Moreover, one dimensional NNMs
appear such as: elliptic periodic orbits or even hyperbolic ones
with their two dimensional  
(coinciding) stable and unstable manifolds.
As the parameters vary,
  six possible topologically
  different phase portraits
  appear.
  An example is given in Figure \ref{scacchi}.

\begin{figure}[h!]
\center
\includegraphics[width=7cm, keepaspectratio]{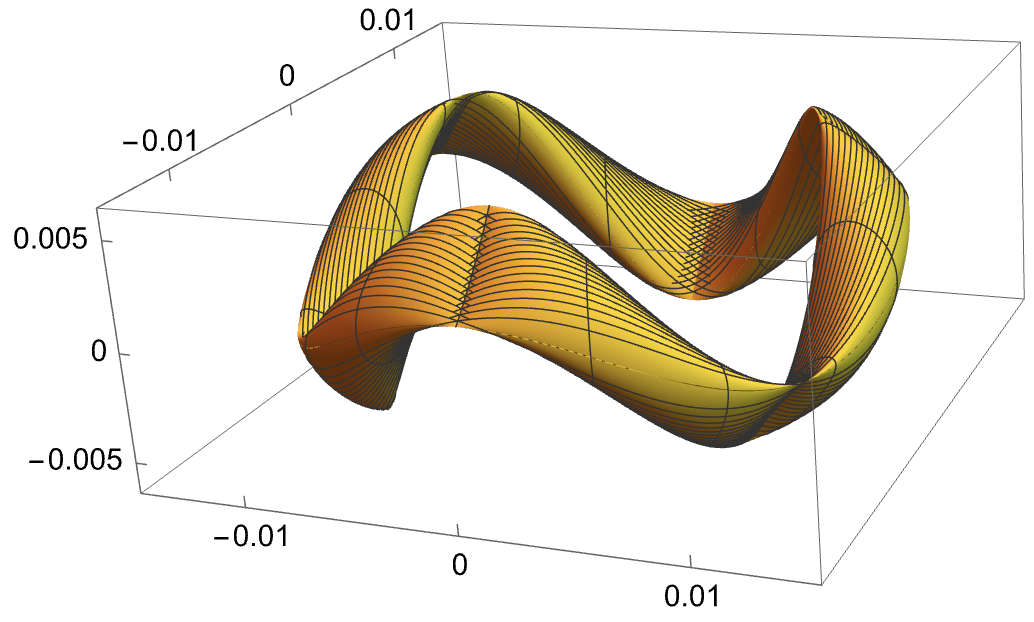} \quad
\includegraphics[width=5cm, keepaspectratio]{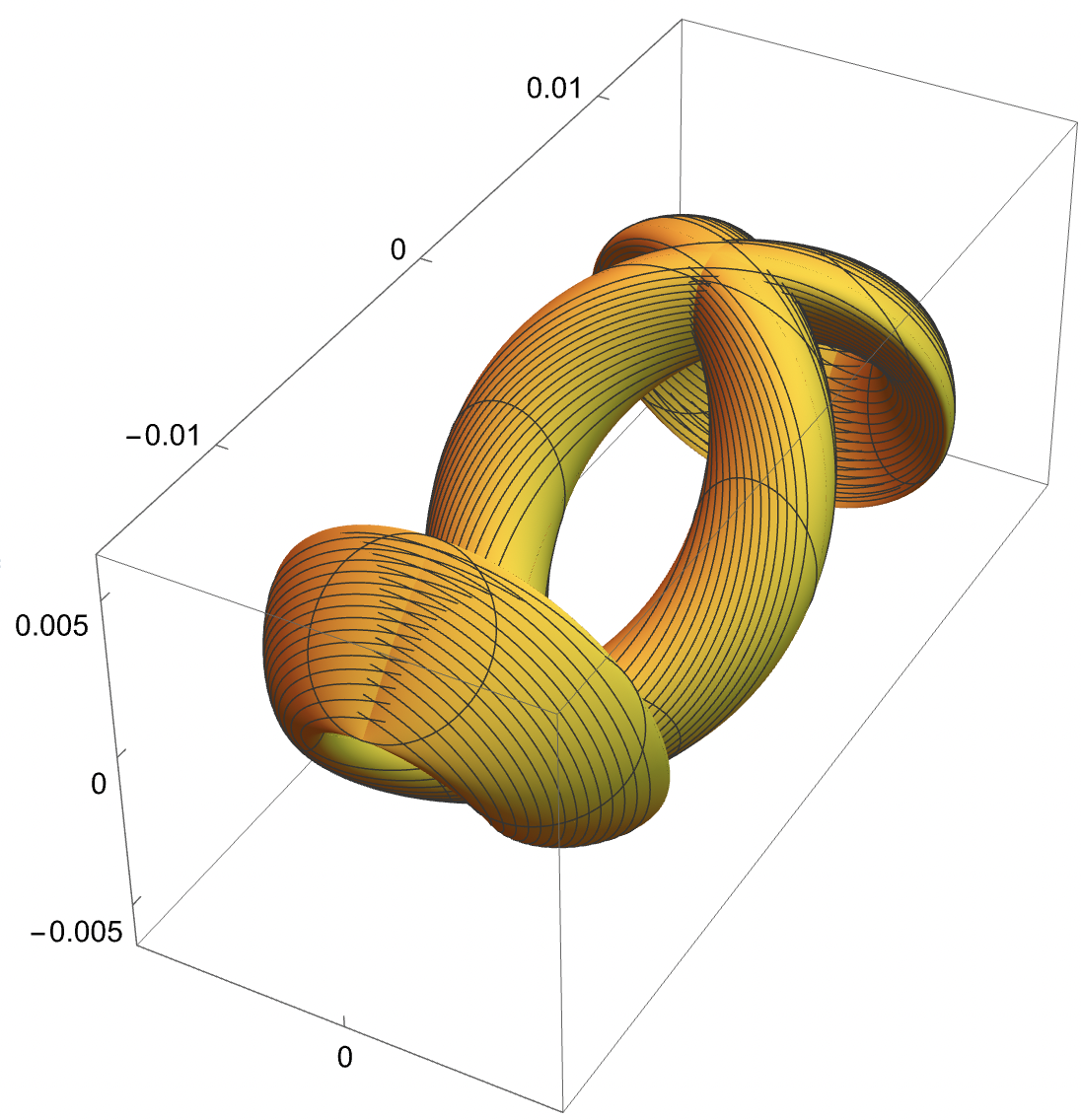} \\
\includegraphics[width=7cm, keepaspectratio]{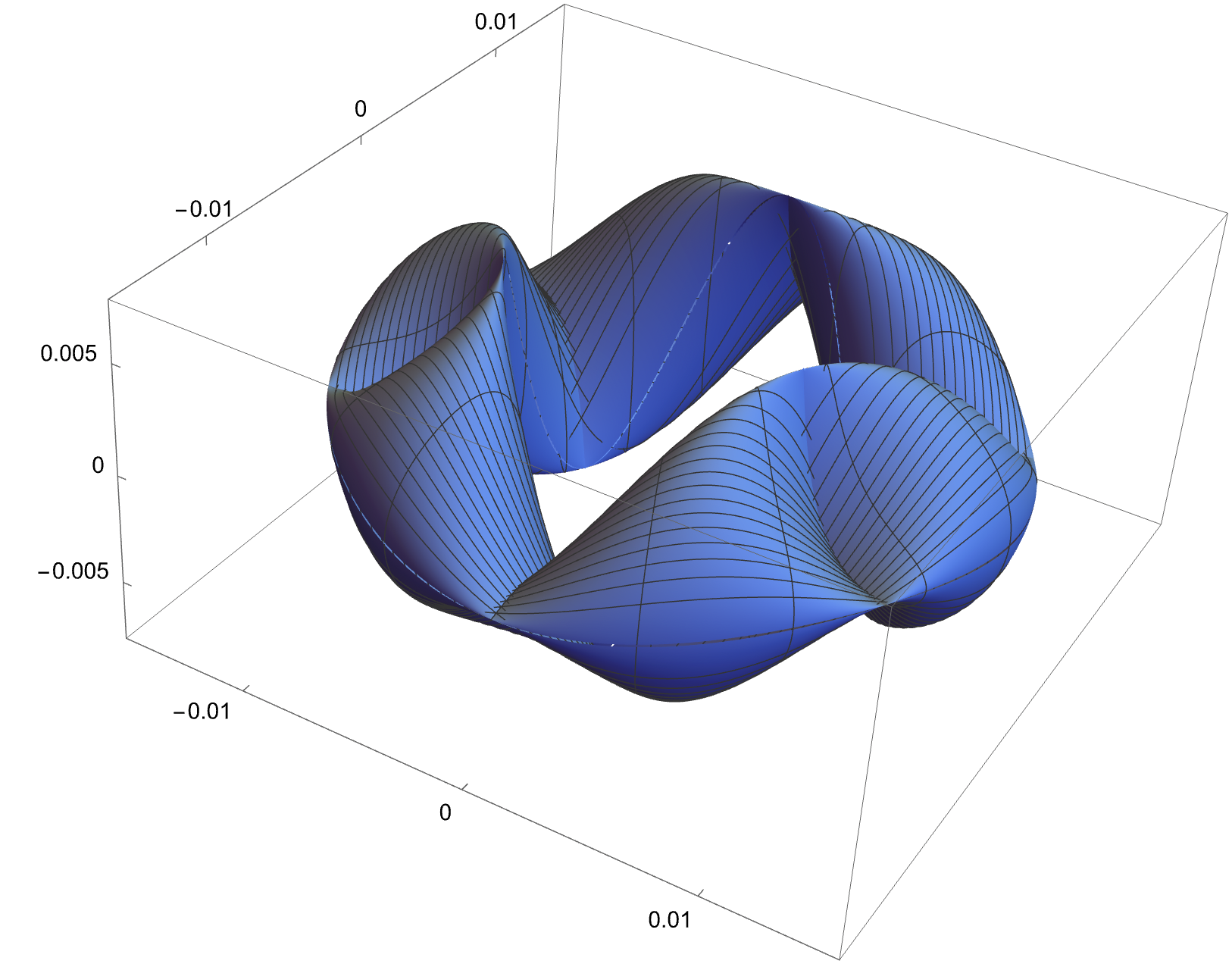}\quad
\includegraphics[width=5cm, keepaspectratio]{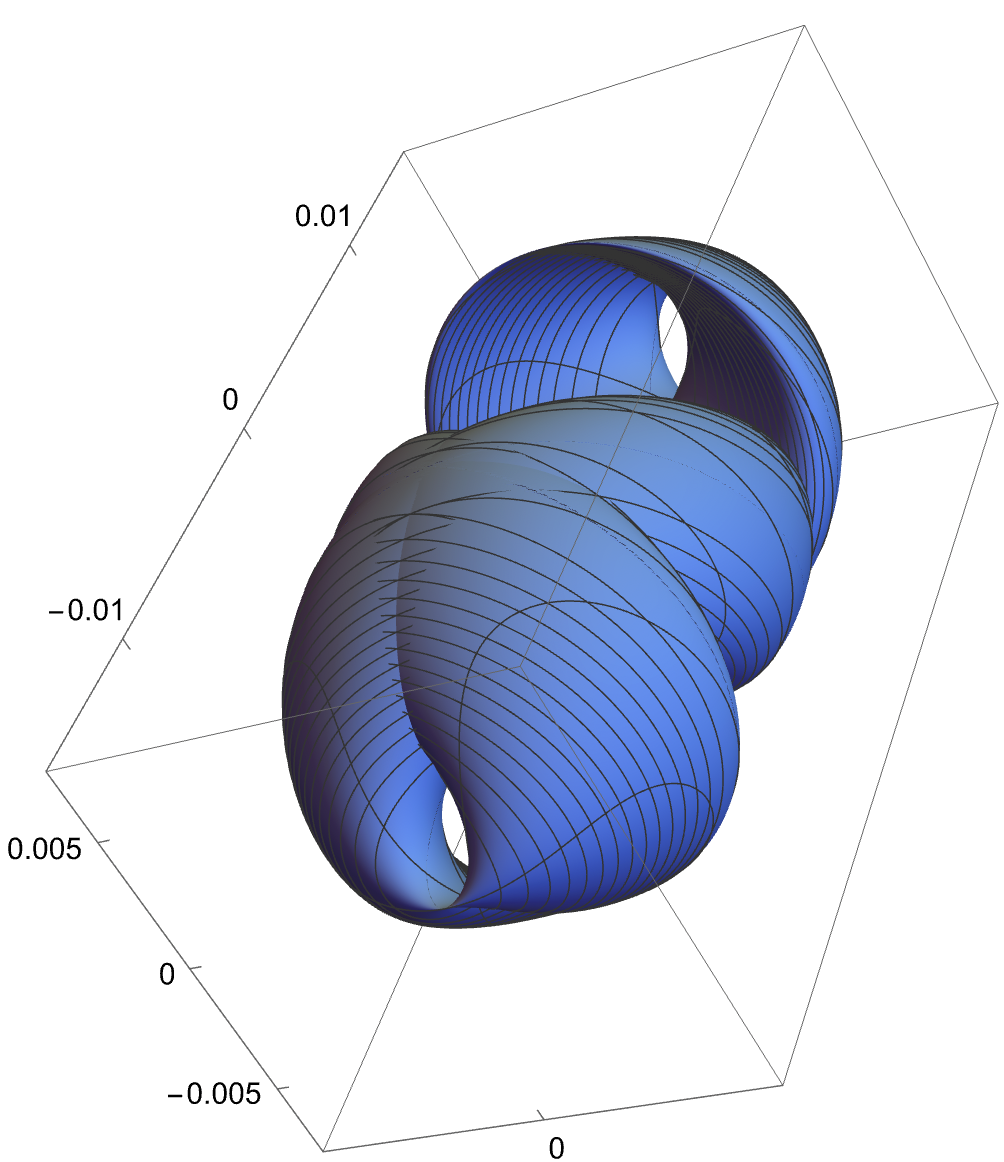}
\caption{The same yellow and blue two dimensional
NNMs  of Figure \ref{scacchi} are plotted here in the modal subspaces $(\dot q_1,q_1,q_2)$
on the left and $(\dot q_2,q_1,q_2)$ on the right (see \eqref{crostini}). Note that, being manifolds, they 
do not have self-intersections in the complete
four dimensional 
modal phase space $(\dot q_1,\dot q_2,q_1,q_2)$.
However one can plot only a 
projection on a three dimensional subspace, where
 self-intersections may occur.
 }
\label{dama}
\end{figure}

\begin{figure}[h!]
\center
\includegraphics[width=4.5cm, keepaspectratio]{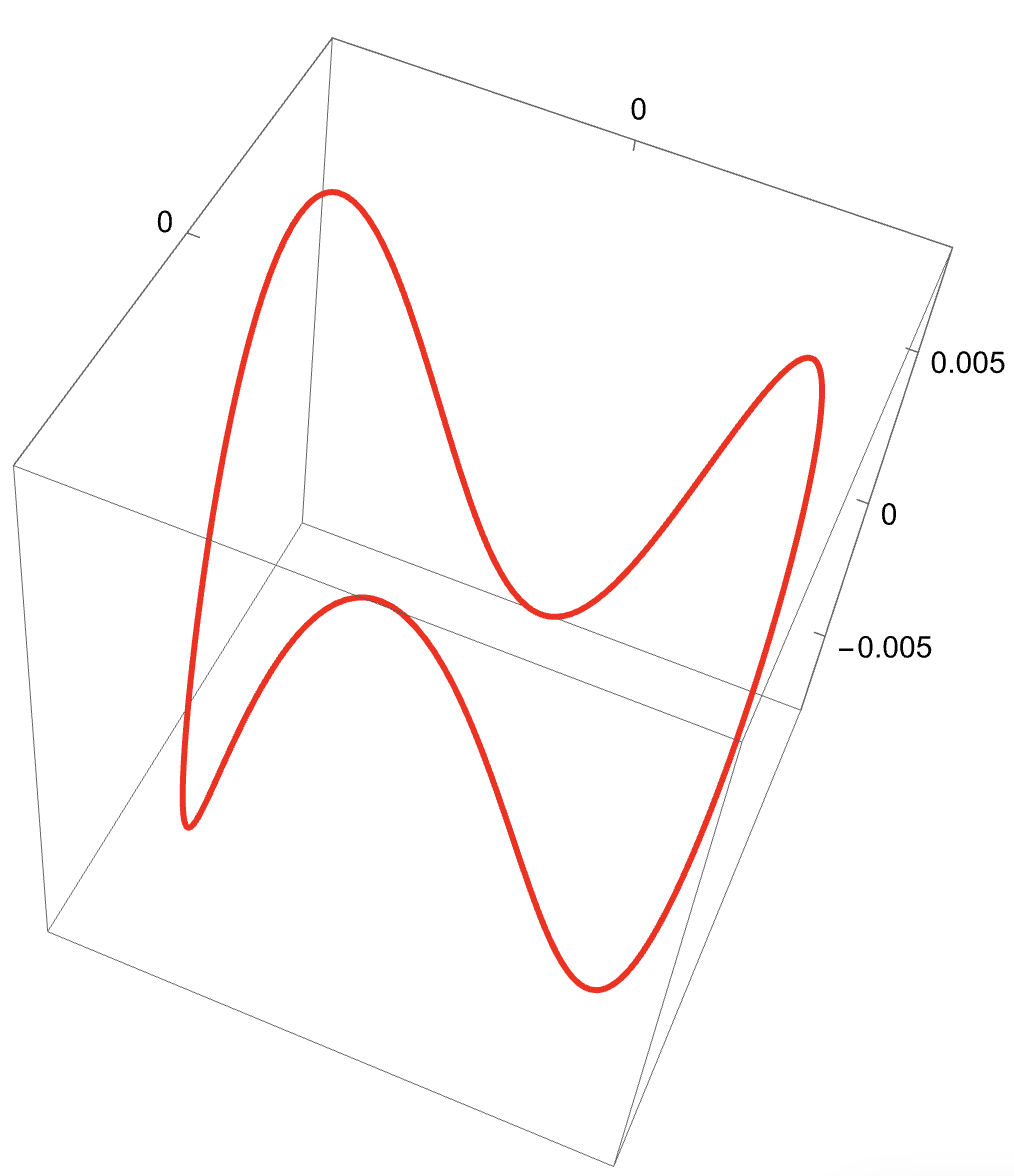} \quad
\includegraphics[width=4.5cm, keepaspectratio]{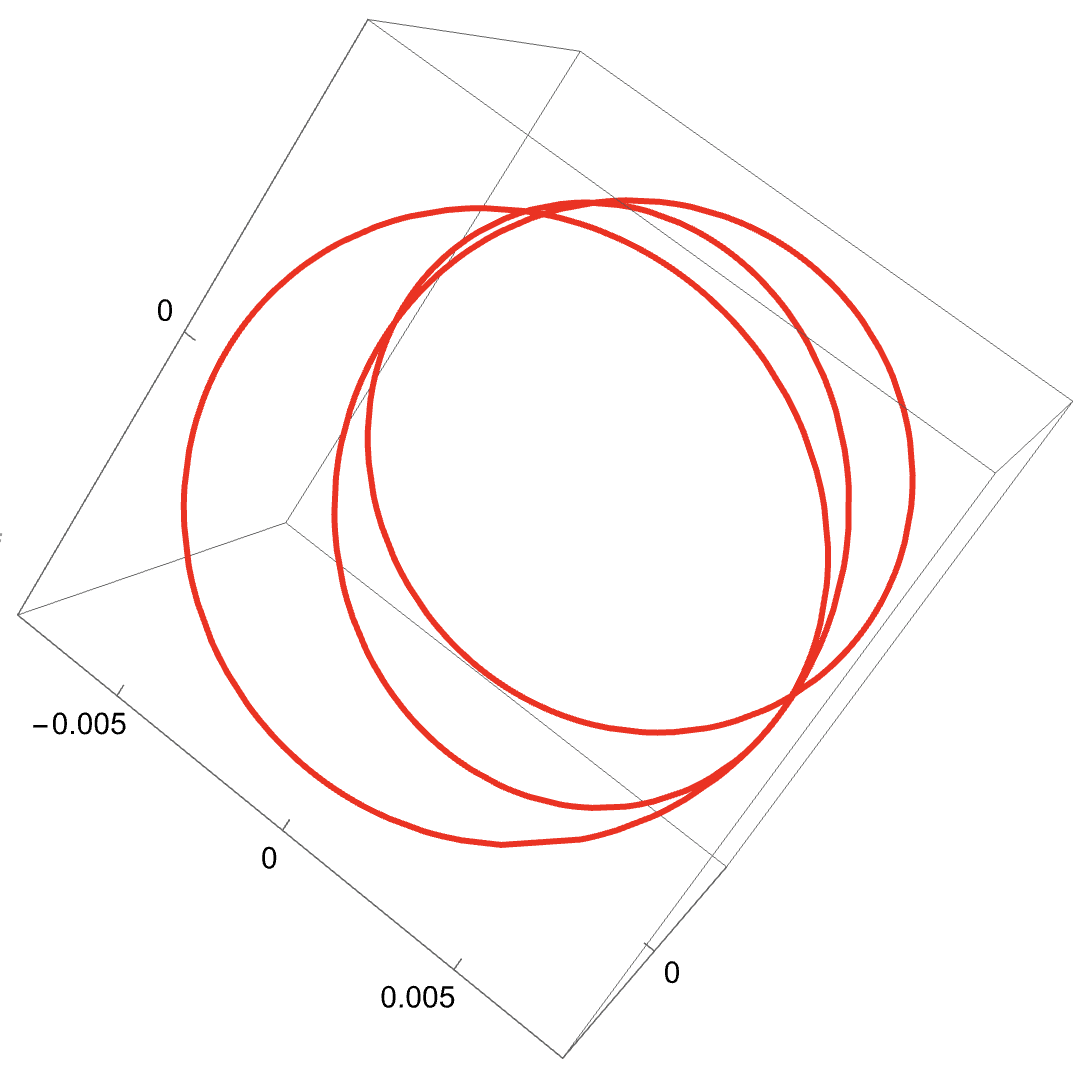} \\
\includegraphics[width=7cm, keepaspectratio]{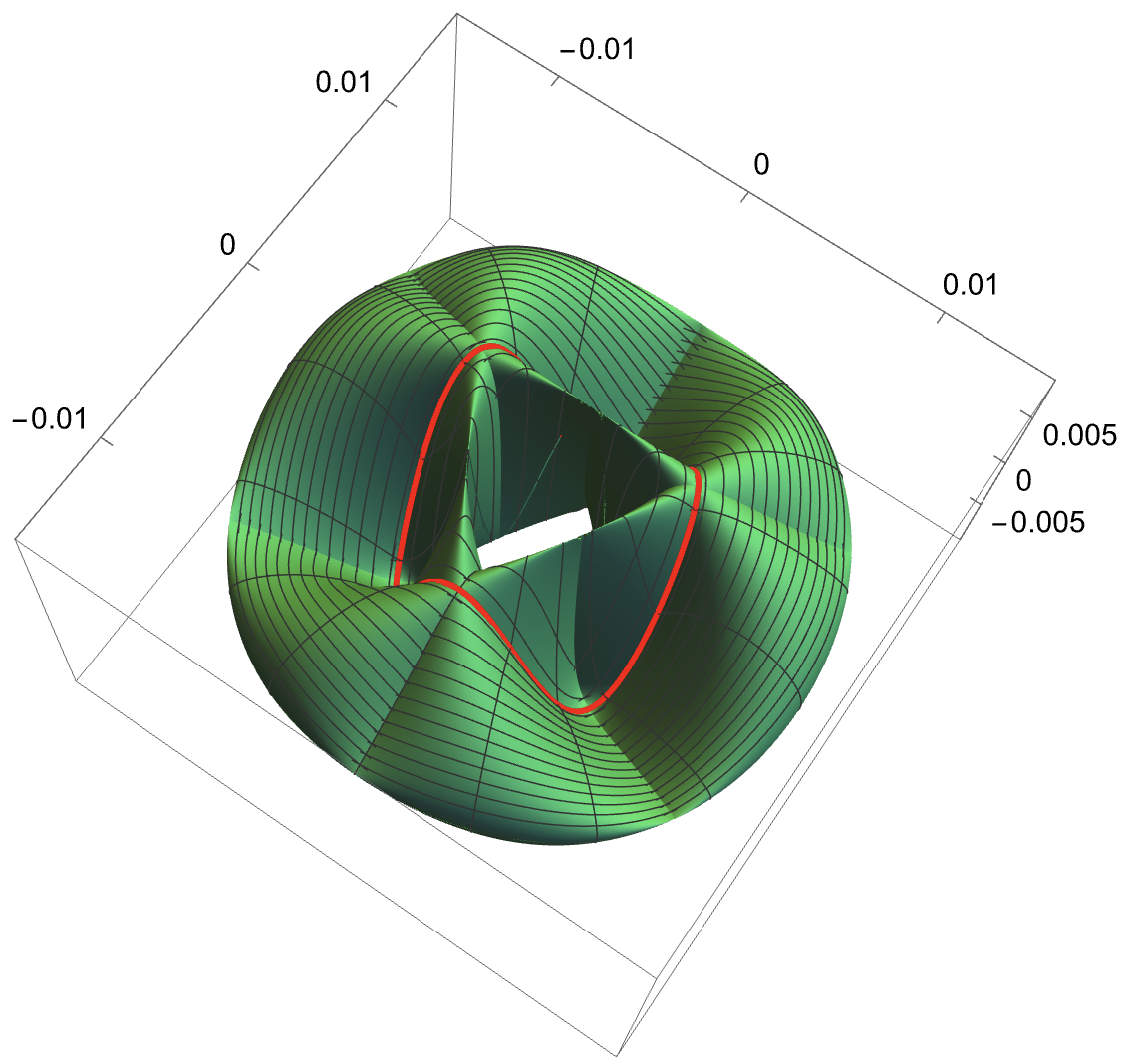}\quad
\includegraphics[width=7cm, keepaspectratio]{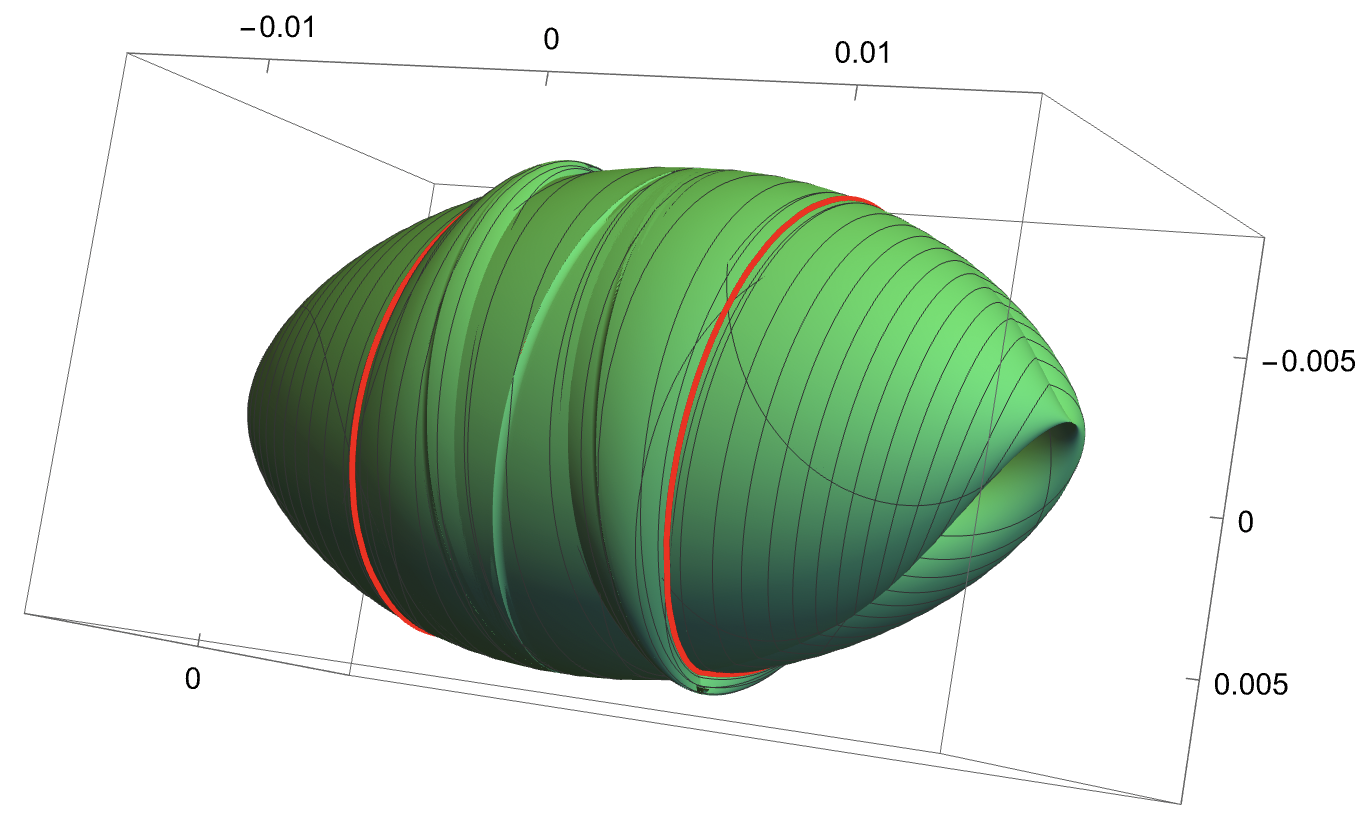}
\caption{The same red one dimensional NNM and 
its green separatrix of Figure \ref{scacchi} are plotted here in the modal subspaces $(\dot q_1,q_1,q_2)$
on the left and $(\dot q_2,q_1,q_2)$ on the right (see \eqref{crostini}).
 Note the triangular symmetry, which is particular
 evident in the green separatrix. It is due to the 3:1
 resonance.}
\label{dama2}
\end{figure}

Let us denote by $\Ju$ the action conjugated to the other angle, the 
``fast'' one, which does not appear in $\hat{\mathbb H}_{\rm res}$.
Then $\Ju$  is a constant of motion for $\hat{\mathbb H}_{\rm res}$.
For every fixed value of $\Ju$, $\hat{\mathbb H}_{\rm res}$
evaluated at $\Ju=const$
in the reduced  bidimensional phase space
 containing only the slow angle
and its conjugated action
is a 1-degree-of-freedom Hamiltonian  system.
In this reduced system, the above two dimensional NNMs
(invariant tori) correspond to  one dimensional NNMs
(periodic orbits), one dimensional NNMs (elliptic/hyperbolic
periodic orbits) correspond to zero dimensional NNMs (elliptic/hyperbolic
fixed points) and, finally, two dimensional  (coinciding) stable and 
unstable manifolds correspond to one dimensional 
(coinciding) stable and 
unstable separatrices, respectively.  
Some examples are shown in Figures \ref{scacchi},
\ref{dama} and \ref{dama2}.

Up to the 
singular\footnote{We call it singular since
it is formed by all the points whose energy
is singular, namely corresponds to some 
critical value of the Hamiltonian. } set 
formed by the union of zero dimensional NNMs (equilibria) and one
dimensional separatrices,
the phase space of 
the reduced Hamiltonian is separated into two or 
four\footnote{According to the different values
of the parameters. In Figure \ref{scacchi}
a case with four regions is shown.}
open connected components having different topologies.
Since the reduced system has one degree of freedom,
on such connected components one can introduce
suitable new action-angle coordinates, integrating the
system. Recollecting, in these new variables, 
$\hat{\mathbb H}_{\rm res}$ depends only on the new
actions and the nonlinear frequencies are simply 
obtained as the derivatives of the Hamiltonian
with respect to the actions.

  However, we note that, at this stage,
   the nonlinear frequencies take the form
   of elliptic integrals, which are not simple to 
   explicitly evaluate since both the integrating
   functions and the domains strongly depend
   on parameters. Nevertheless,
     we calculate them by using suitable
   Moebius transformations.

   Finally, having the explicit formulas available, 
we study the nonlinear bandgap in the  resonant regime.
We found that, while the nonlinearity far from resonances can significantly change the bandgap, in the resonant case, the effect of resonances  results in a less pronounced variation in the bandgap.

 Here we study in details the truncated Hamiltonian
 giving a very precise description of its phase space 
 and  explicitly  integrating the system.
 The case of the complete Hamiltonian is different since the system is genuinely two dimensional
 and, therefore, not integrable\footnote{Since the fast angle
 appears at higher order terms and, therefore, its conjugated
 action $\Ju$ is not more a constant of motion.}.
 However,  using methods of KAM Theory
 one can prove the persistence of  hyperbolic
 periodic orbits with their (local) stable and unstable
 manifolds as well as of the majority of invariant tori.
Indeed, our analysis can bee seen as a 
  necessary 
  preparatory 
  step toward applying  
  KAM techniques  \emph{in the resonant zones}
   (see Remark \ref{KAM}).

Finally, we stress that
our analysis is not limited to the case of the honeycomb metamaterials  but applies directly to a wide range of problems modeled by two harmonic oscillators coupled with cubic nonlinearity as in equation
\eqref{autostrada}.

\begin{figure}[h!]
\begin{minipage}{.45\textwidth}
\center
\includegraphics[width=7.4cm, keepaspectratio]{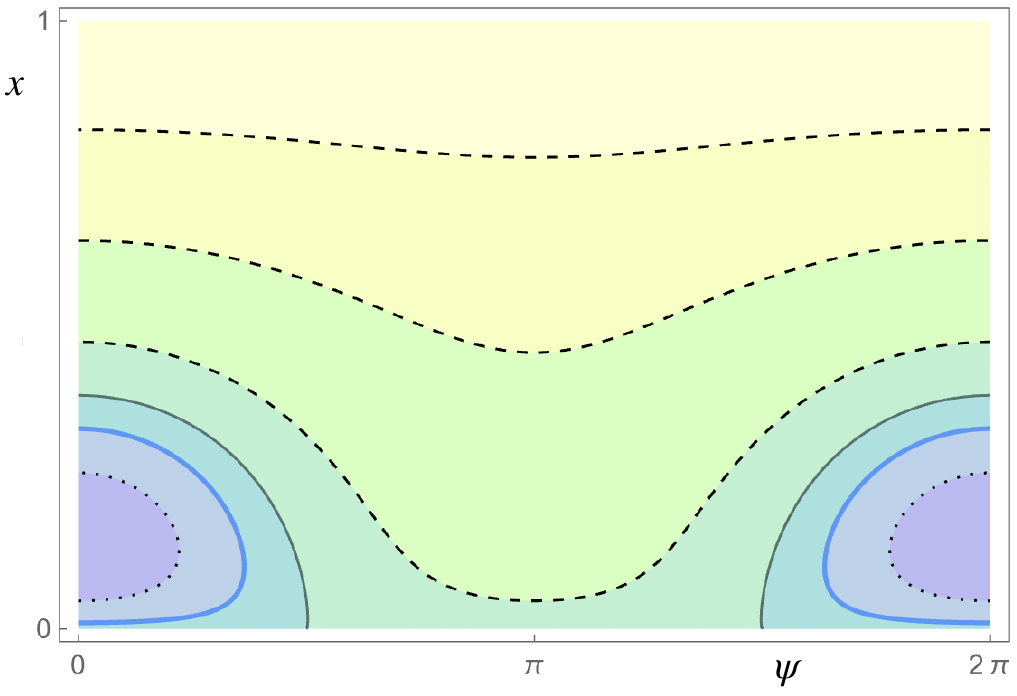}\\
\includegraphics[width=7.5cm, keepaspectratio]{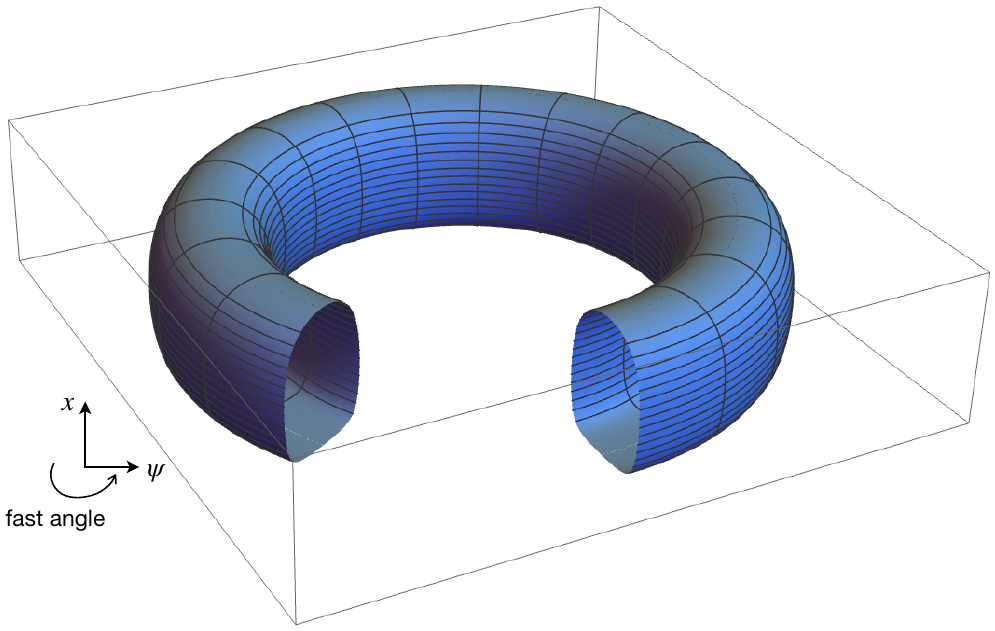}
\end{minipage}
\begin{minipage}{.42\textwidth}
\includegraphics[width=8.2cm,keepaspectratio]{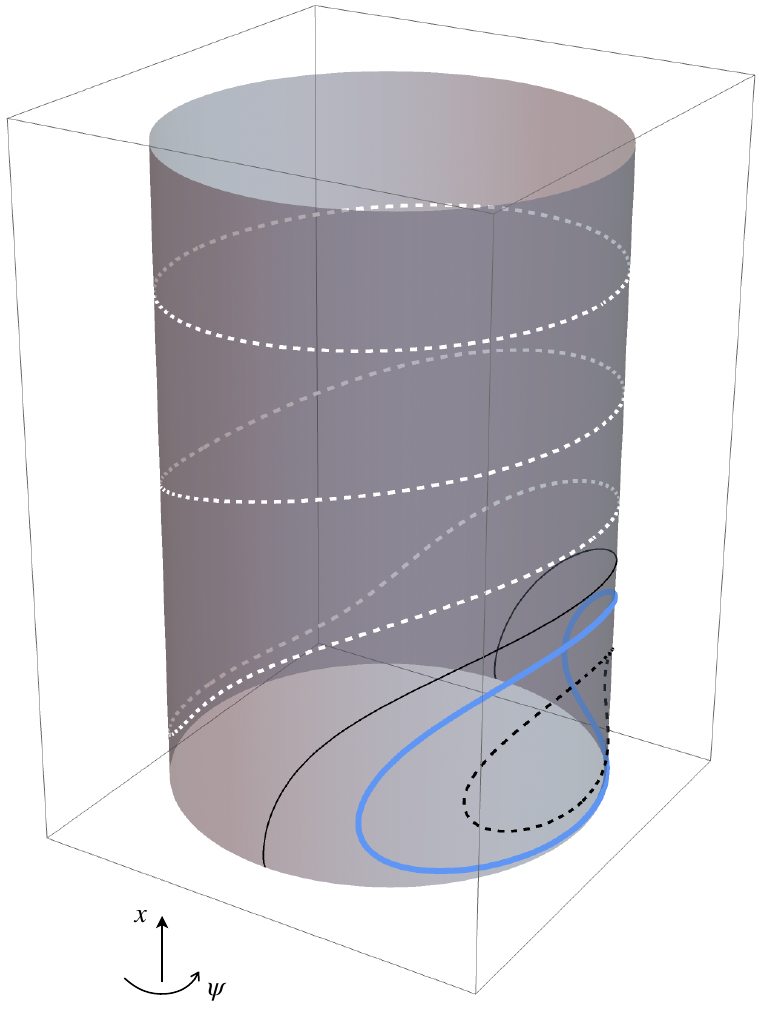}
\end{minipage}
\caption{(Left top) Level curves of 
the reduced Hamiltonian obtained 
by the fourth order resonant  
$\hat{\mathbb H}_{\rm res}$ (see \eqref{secularRES})
fixing the constant of motion $\Ju=10^{-4}$
(here, e.g., we have chosen the physical parameters as follows: $\tilde{k}_1=2.27, \tilde{k}_2=0, \tilde{M}=0.2, \tilde{K}=1.1, M_3=0, N_3=-10^4$).
The slow angle is on the horizontal axis and 
its (rescaled) conjugated action is on the vertical axis,
so that the phase space is actually the cylinder shown on the right.
Except for the   solid black curve that 
acts as separatrix, the phase
space is divided into two connected components.
Every component is completely foliated  by one dimensional NNMs 
(periodic orbits).
Such NNMs have different topology:
the orbits in the zone above the separatrix
wrap around the cylinder (dashed curves);
the orbits inside the separatrix
do not wrap around the cylinder and are contractible
(dotted and blue  curves).
(Right) The  cylindrical phase portrait
 immersed in the three dimensional space.
\\
(Left bottom) 
A representation
 of the phase space
of the truncated Hamiltonian $\hat{\mathbb H}_{\rm res}$,
once we have fixed the constant of motion $\Ju=10^{-4}$.
 The image is 
obtained by rotating the picture on the left by  the fast angle   from $0$ to $4\pi/3$.
In particular, by rotation,
 the blue and yellow curves become 
two dimensional NNMs (invariant tori) and the red point and the green curve become, respectively,
a one  dimensional NNM (a hyperbolic periodic orbit) 
and its two dimensional (coinciding) 
stable and unstable manifolds. 
}
\label{scacchiscacchi}
\end{figure}

\begin{figure}
\center
\includegraphics[width=7cm,keepaspectratio]
{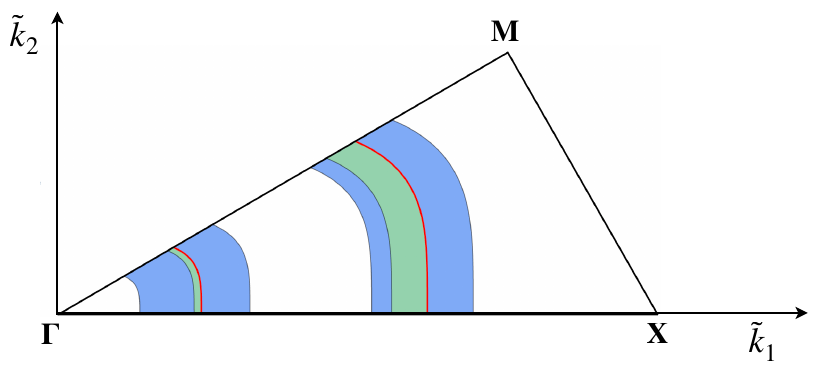}
\qquad \qquad
\includegraphics[width=4cm,keepaspectratio]
{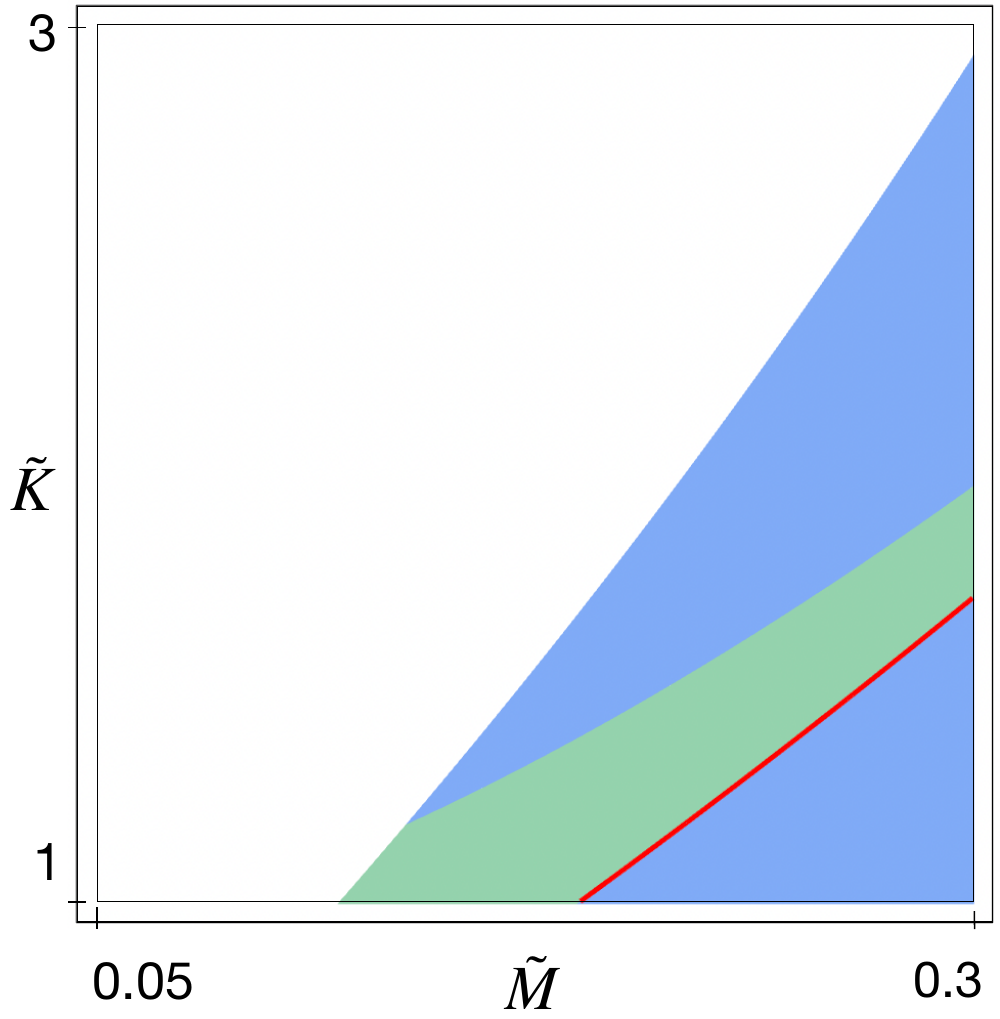}
\caption{(On the left) For $\tilde M=0.2$, $\tilde K=1.1$, $N_3=-10^{4}$
and $\varepsilon\sim 10^{-4}$,
according to the values of the wave numbers
in the Brillouin triangle, we have different BNFs.
In the white regions,
one can construct a  nonresonant BNF.
In the blue and green regions, 
 only a resonant BNF is available. 
 In particular, in the green regions, the phase portrait
 of the Hamiltonian in BNF is as in Figure \ref{scacchi},
 while, in the blue regions, the phase portrait
 of the Hamiltonian in BNF is as in Figure \ref{scacchiscacchi}.
 The red curves represent the pairs $(\tilde k_1,\tilde k_2)$
 for which $\sigma=0$, namely when the exact 3:1 resonance occurs.
 (On the right)
 Here we fix $\tilde k_1=2.58$, $\tilde k_2=0$,  
 and let 
 $(\tilde M,\tilde K)$ vary in the rectangle 
 $[0.05,0.3]\times[1,5]$.
 A more quantitative and precise description 
 will be given in Figure \ref{clarissa}.}
\label{clarissa0}
\end{figure}

\begin{figure}
\center
\includegraphics[width=14cm,keepaspectratio]{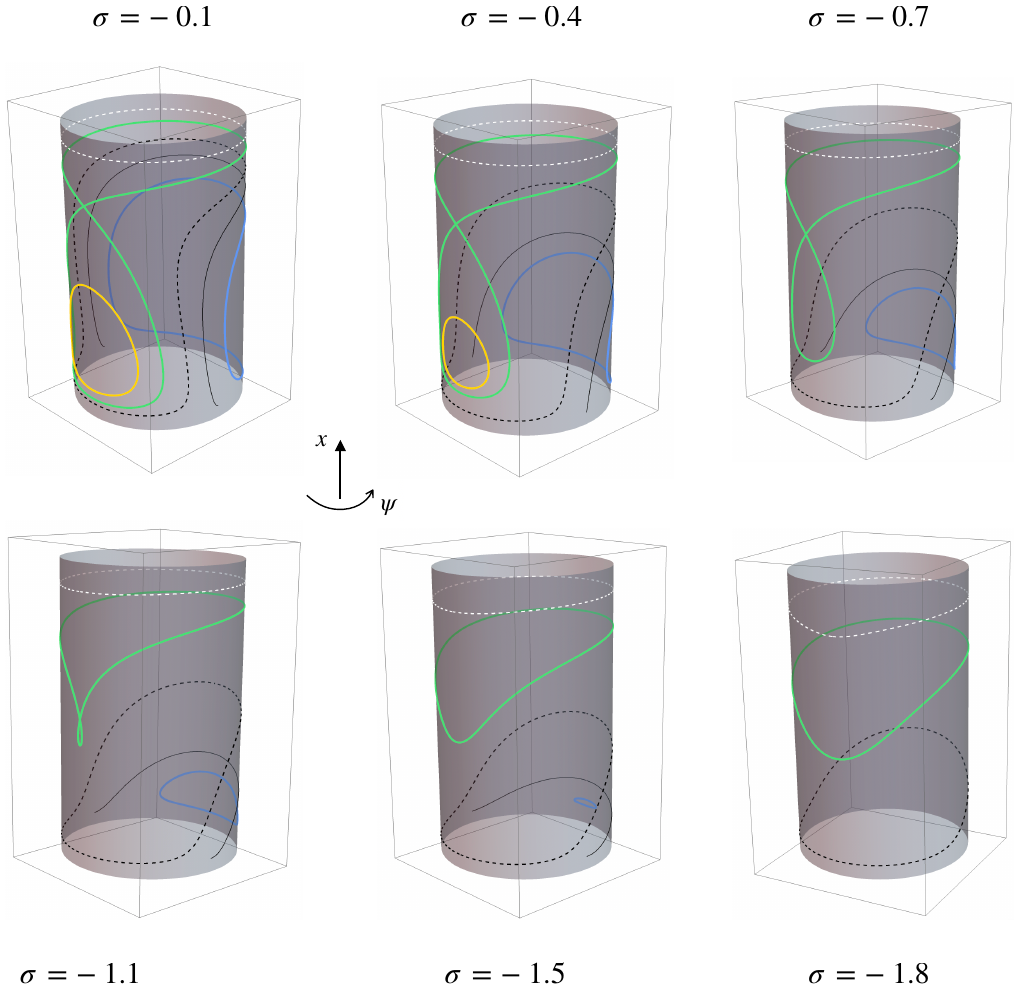}
\caption{Referring to Figure \ref{clarissa0}, six different phase
portraits are shown as $\sigma$ varies from $-0.1$ to $-1.8$.
Note that in Figure \ref{clarissa0}  the red curves on the right
corresponds to the exact 3:1 resonance $\sigma=0$.
Moving on the left and passing through the green, blue,
and, finally, white zone,  the value of $\sigma$ decreases
and the topology of the phase space changes.
For $\sigma=-0.1, -0.4,-0.7,-1.1$ (corresponding to wave numbers
in the green region) the phase space has the same topology
 shown in Figure \ref{scacchi} with a hyperbolic point with its
 coinciding stable and unstable manifolds (green curve) 
 and two periodic orbits not wrapping on the cylinder:
  the yellow curve, that soon contracts and disappears,
   and the blue curve, that reduces.
 For $\sigma=-1.4$ (corresponding to wave numbers
in the blue region), the phase space has the same topology
 shown in Figure \ref{scacchiscacchi}: the hyperbolic point 
 does not exist 
 since the green curve does not self-intersect
 anymore and becomes a simple periodic curve wrapping 
 on the cylinder, while the blue curve becomes smaller and smaller.
  Finally for $\sigma=-1.8$ (corresponding to wave numbers
in the white region), the  blue curve disappears and 
only periodic curves wrapping on the cylinder survive,
showing the typical behavior of integrable systems. }
\label{xxxxx}
\end{figure}

\subsection{Summary of the paper}

\subsubsection*{Section \ref{sec:Ham}: the resonant Birkhoff Normal Form}
\vspace{-0.2cm}

We reinterpret  the problem as a Hamiltonian 
system (see \eqref{ham}).
In Subsection \ref{Sec:BNF},
 we put the system, close to the origin, in resonant BNF. 
 Then, we examine the Hamiltonian truncated
 at fourth order, which is
 equivalent to third order in the 
 equations of motion, 
 as it captures the essential
 characteristics 
  of the overall motion.
 Upon introducing action-angle variables
 it becomes evident that 
 the truncated, or ``effective'',
 Hamiltonian, after a suitable linear change 
 of variables (see \eqref{salmone}),
 also depends on one angle,  known as the  ``slow''
 angle 
(see \eqref{secularRES}), 
as its associated  
frequency is small or even zero on the exact resonance.

After a suitable rescaling, the effective
Hamiltonian,
depending on the slow angle 
$\psi\in[0,2\pi)$
and on the non-dimensional action $x\in(0,1)$,
takes the form $F(\psi,x)=\frac12 a_2x^2+a_1x+b(x)\cos\psi$, where $b(x)=\sqrt{(1-x)^3 x}$
and $a_1,a_2$ depend on the physical parameters and on the other action
(which is a constant of motion);
see Subsection \ref{SubSec:effective Hamiltonian}.

\subsubsection*{Section \ref{sec:phase}: the six possible phase portraits
}
\vspace{-0.2cm}

\begin{figure}
\center
\includegraphics[width=10cm, keepaspectratio]{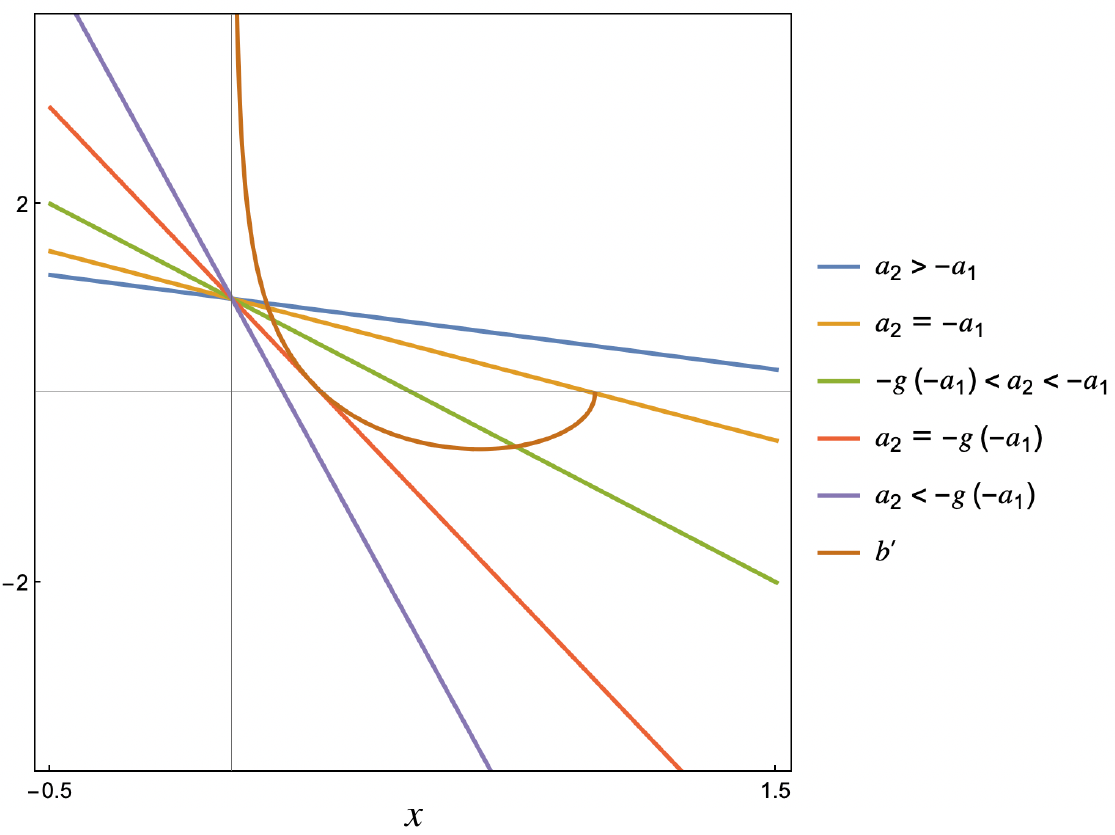}
\caption{Intersections between the function
$b'(x)$ and the straight line 
 $a_2 x + a_1$ for a fixed $a_1$ and different values
 of $a_2$. The green line, passing through $(1,0)$,
 and the purple one, which is tangent to $b'(x)$, 
 separate the half plane $x>0$ in three regions, in which 
 the lines intersect 1,2 or 0 times the curve $b'(x)$.
 More precisely
 for $a_2>-a_1$ there is one intersection,
 for $-g(-a_1)<a_2<-a_1$ there are two intersections and none 
 for $a_2<-g(-a_1)$.}
\label{Intersezione_b'_a'}
	\end{figure}

\begin{figure}
		\center{\includegraphics[width=12cm]{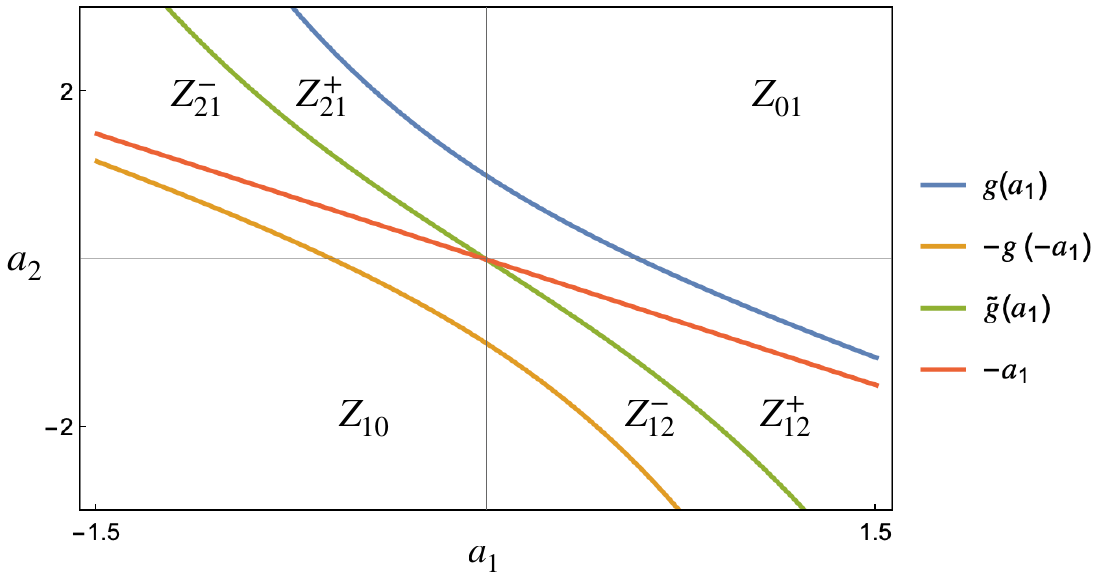}}
		\caption{The six zones
		$Z_{01},Z_{10},Z_{21}^+,Z_{21}^-,
		Z_{12}^+,Z_{12}^-$.}
		\label{Zones_6}
	\end{figure}

The behavior of the system depends
on the number and on the nature
of the critical points of $F$,
which, in turn, depends
on the values of $a_1$ and $a_2$.
The gradient of $F$ can vanish 
only on the lines $\{\psi=0\}$,
when $a_2x+a_1+b'(x)=0$, or 
$\{\psi=\pi\}$, when $a_2x+a_1-b'(x)=0$.
At this point studying the 
solutions of these equations, as
$a_1$ and $a_2$ vary, is crucial
(see Figure \ref{Intersezione_b'_a'}). 
This identifies six zones
in the plane $(a_1,a_2)$, as detailed in
Proposition \ref{pandivia},
Lemma \ref{Lemma6}
and Figure 	\ref{Zones_6}.
Correspondingly we have six possible configurations.
When reached,
the maximum of $F$ is attained  
on the line $\{\psi=0\}$, conversely the minimum
is attained  
on the line $\{\psi=\pi\}$.
E.g. let us briefly describe  the
scenario $a_1+a_2<0$.
By studying $x\to F(0,x)$
we have three  possible cases:
 no critical points,
a maximum and a minimum with negative
energy, a maximum and a minimum with positive
energy. On the other hand
$x\to F(\pi,x)$ has a minimum.
 Note that the maximum of $F(0,x)$
corresponds to a maximum for 
$F(\psi,x)$, the minimum of $F(0,x)$
corresponds to a saddle for 
$F(\psi,x)$ and the minimum of
$F(\pi,x)$ corresponds to a minimum of $F(\psi,x)$.
Analogously, the complementary case 
$a_1+a_2>0$ gives rise to 
three additional configurations.

\subsubsection*{Section \ref{sediablu}: construction of the integrating action variable
}
\vspace{-0.2cm}

Since  the action conjugated to the ``fast''
angle is a 
constant of motion, the
truncated system 
has two independent conserved quantities
(the other one is the energy) and, therefore,
is  integrable (by the Arnold-Liouville
Theorem), in the sense that
one can find a new set of symplectic action-angle variables
in which the new Hamiltonian depends 
{\sl only on the actions}. Although the theoretical construction of the integrating action is classic, finding an {\sl explicit analytical} expression as a function of all the physical parameters involved is rather complicated.

For every value of the energy $E$, 
the new integrating action $\II_1$ is given by
the area enclosed by the level curve
$F(\psi,x)=E$ divided by $2\pi$,
see Section \ref{sediablu}.
Such level curves  are closed and
 can either wrap around the cylinder
$[0,2\pi)\times(0,1)$
or remain confined to its surface
without wrapping around it;
see Figures  
\ref{scacchi} and \ref{scacchiscacchi}.

Since $F$ is even in $\psi$
we can restrict to consider 
$(\psi,x)\in[0,\pi]\times(0,1)$.
In this set the level curves are 
graphs over $x$ and the area enclosed
by them can be computed
by an integral over $x$, whose endpoints
are the $x$-coordinate
of their  intersections with the lines
$\{\psi=0\}$ and $\{\psi=\pi\}$.
It turns out that these 
correspond to
 the roots
$0<x_j(E)<1$, with $j=1,2,3,4$,
of the quartic polynomial
 $\mathbf P(x;E)=
 \big(\frac12a_2x^2+a_1x-E\big)^2-(b(x))^2$, see \eqref{penisola}
and
Figure \ref{4_soluzioni}.
As the energy $E$ varies, 
it is necessary to distinguish 
whether $\mathbf P$ has
$4,2$ or $0$ real 
roots\footnote{Note that, 
excluding
the degenerate case of multiple roots,
the number of real roots
is even.} and
whether  a root corresponds to
an intersection with
$\{\psi=0\}$ or $\{\psi=\pi\}$.
Explicit formulae for the roots
are given in Subsection \ref{Sec:radici},
see Figure \ref{xPiuMenoSoluzioni}.

Once we have defined the integrating
action $\II_1$ as a function of $E$
(and of the ``dumb'' action, let us say, $\II_2$),
the resulting integrated Hamiltonian will be its inverse
$E=E(\II_1,\II_2)$.
The nonlinear frequencies
are given by the derivatives of the energy
with respect to the actions,
see \eqref{california}, expressed  through integrals,
see Proposition \ref{WWW}.
Such integrals are evaluated by suitable Moebius transformations  in terms of 
 elliptic functions, see Subsections
 \ref{Sec:elliptic} and \ref{somorta}.

\subsubsection*{Section \ref{sec:bandgap}: evaluation of the 
nonlinear bandgap for the honeycomb metamaterial}
\vspace{-0.2cm}

Finally, 
 having the explicit formulas for the nonlinear frequencies available, 
we discuss the nonlinear bandgap 
for the honeycomb metamaterial, 
especially
 in the  resonant regime.
We found that, while nonlinear effects far from resonances can significantly alter the bandgap, in the resonant case the nonlinear frequencies, especially the acoustic one, closely align with the linear frequencies, resulting in a less pronounced variation in the bandgap.

\section{The Hamiltonian structure and resonant BNF}
\label{sec:Ham}

In this section, after
introducing optical and acoustic 
modes,
we identify the system
in \eqref{autostrada}  as  Hamiltonian,
see \eqref{ham} below, and we evaluate
the coefficients of the Hamiltonian, see
\eqref{Gab}.
Set
$$
\mathbf \Lambda:=
 \left(
\begin{array}{cc}
	\omega_-^2&0\\
       0& \omega_+^2
	\end{array}	
\right)\,, 
$$
where
$\omega_-^2<\omega_+^2$ are the positive
eigenvalues of $\mathtt M^{-1}\mathtt K$ and
 $0<\omega_-<\omega_+$.
Since $\mathtt M$ is symmetric and 
$\mathtt K$ is diagonal,
there exists a $2\times 2$ matrix $\mathbf  \Phi$  such that
\begin{equation}\label{simultaneo}
\mathbf \Phi^T \mathtt M\mathbf \Phi=\mathbf  I\,,
\qquad
\mathbf \Phi^T \mathtt K \mathbf \Phi=\mathbf \Lambda\,,
\qquad
\mathbf  \Phi
 =\left(
\begin{array}{cc}
	\phi_1^-& \phi_1^+\\
       \phi_2^-& \phi_2^+
	\end{array}	
\right)\,,
\end{equation}
where $\mathbf  I$ is the identity matrix. 
\noindent
Consider the change of variables
\begin{equation}\label{crostini}
\left(
		\begin{array}{c}
			v \\
			y
		\end{array}
		\right)
		=\mathbf \Phi \mathbf q\,,\qquad
		\mathbf q:=\left(
		\begin{array}{c}
			q_1 \\
		q_2
		\end{array}
		\right)\,.
\end{equation}
	By Lemma \ref{simultaneo}
	the system in \eqref{autostrada}
	is transformed into 
	\begin{equation}\label{fettuccine}
	\ddot {\mathbf q}+\mathbf \Lambda \mathbf q=\mathbf c(\mathbf q)\,,\qquad
	\mathbf c(\mathbf  q)=\left(
		\begin{array}{c}
			c_1 \\
		c_2
		\end{array}
		\right):=-\mathbf \Phi^T\left(
		\begin{array}{c}
	M_3 (\phi_1^- q_1+\phi_1^+ q_2)^3 \\
	N_3 (\phi_2^- q_1+\phi_2^+ q_2)^3 \\
	\end{array}	
	\right)\,.
\end{equation}
In particular
$$
			\begin{array}{ll}
					c_1=-\phi_1^- M_3(\phi_1^- q_1+\phi_1^+ q_2)^3-\phi_2^- N_3(\phi_2^- q_1+\phi_2^+ q_2)^3\\
					c_2=
					-\phi_1^+M_3(\phi_1^- q_1+\phi_1^+ q_2)^3-\phi_2^+N_3(\phi_2^- q_1+\phi_2^+ q_2)^3\,.
			\end{array}		
	$$
Introducing the momenta 
$\dot {\mathbf q}=\mathbf p=\left(
		\begin{array}{c}
			p_1 \\
		p_2
		\end{array}
		\right)$,	 the system in \eqref{fettuccine}
	is Hamiltonian with Hamiltonian
	\begin{equation}\label{ham}
		H(\mathbf p,\mathbf q)=\frac 12 (p_1^2+p_2^2) + \frac12 \omega_-^2 q_{1}^2+ \frac12 \omega_+^2 q_{2}^2+f(q)\,,
	\end{equation}
where 
	\begin{equation}
				f(\mathbf q)
			:= \frac 14 M_3 (\phi_1^- q_1+\phi_1^+ q_2) ^4+\frac 14 N_3 (\phi_2^- q_1+\phi_2^+ q_2)^4\,.
				\label{potential}
	\end{equation}
Indeed it is immediate to see that the Hamilton's
equations
$\dot {\mathbf p}=-\partial_{\mathbf q} H$, $\dot {\mathbf q}=\partial_{\mathbf p} H=\mathbf p$
are equivalent to the system 
in \eqref{fettuccine}.
Since $f(\mathbf q)$ is a homogeneous polynomial
of degree $4$ we write
	\begin{equation}\label{fij}
	f(\mathbf q)=\sum_{i+j=4}f_{i,j}q_1^i q_2^j\,,
	\qquad
	\mbox{with}
	\quad			f_{i,j}:=\frac{6}{i!j!}
				\Big(
				(\phi_1^-)^i(\phi_1^+)^jM_3+(\phi_2^-)^i(\phi_2^+)^jN_3
				\Big)\,.
	\end{equation}
Introducing coordinates 
$\mathbf Q=(Q_1,Q_2)$, $\mathbf P=(P_1,P_2)$ through 
\begin{equation}\label{islanda}
			p_1=\sqrt{\omega_-}P_1\quad p_2=\sqrt{\omega_+}P_2\quad q_1=\frac1{\sqrt{\omega_-}}Q_1\quad q_2=\frac1{\sqrt{\omega_+}}Q_2
		\end{equation}
we have  that the Hamiltonian in the new variables
reads
	\begin{equation}\label{hamPQ}
		\mathtt H(\mathbf P,\mathbf Q):=\omega_- \frac{P_1^2+Q_1^2}{2}+\omega_+ \frac{P_2^2+Q_2^2}{2}+ f\left(\frac{Q_1}{\sqrt{\omega_-}},\frac{Q_2}{\sqrt{\omega_+}}\right)\,.
	\end{equation}
In complex coordinates, 
$\mathrm i=\sqrt{-1}\in \mathbb C $, $\mathbf z=(z_1,z_2)\in \mathbb C^2$
		\begin{equation}\label{islanda2}
			z_j=\frac{Q_j+\mathrm i P_j}{\sqrt 2}\quad \bar z_j=\frac{Q_j-\mathrm i P_j}{\sqrt 2}\qquad j=1,2
	\end{equation}
		the Hamiltonian reads
\begin{equation}\label{hamiltonian}
		\mathtt H(\mathbf z,\bar{\mathbf z})
			=
				\mathtt N(\mathbf z,\bar{\mathbf  z})
				+
				{\mathtt G}(\mathbf z,\bar{\mathbf  z})
\end{equation}
where
\begin{equation}\label{GG}
\mathtt N(\mathbf z,\bar{\mathbf  z}):=
				\omega_- z_1\bar z_1+
				\omega_+ z_2\bar z_2\,,
				\qquad
{\mathtt G}(\mathbf z,\bar{\mathbf  z})
			:=
			f\left(\frac{z_1+\bar z_1}{\sqrt{2\omega_-}},\frac{z_2+\bar z_2}{\sqrt{2\omega_+}}\right)\,.
\end{equation}
Note that in complex coordinates
the  Hamilton's equations of motion are
	\begin{equation}\label{HE}
		 \dot z_j=-\mathrm{i} \partial_{\bar z_j}\mathtt H\,, \ \ \dot{\bar{z}}_j= \mathrm{i} \partial_{z_j}\mathtt H\,.
	\end{equation}
In the following
we use the multi-index notation
\begin{equation}\label{fuga}
		P(\mathbf z,\bar{\mathbf  z})=\sum_{(\bal,\bbe)\in\mathbb N^2\times\mathbb N^2} P_{\bal,\bbe}
		\mathbf  z^\bal
		\bar{\mathbf z}^\bbe
		\end{equation}
for suitable coefficients $P_{\bal,\bbe}\in\mathbb C$
with $\mathbf z^\bal=z_1^{\alpha_1}z_2^{\alpha_2}$ (analogously
for $\bar{\mathbf z}^\bbe$). 
In these notation, recalling
\eqref{fij} and \eqref{GG},
 we rewrite ${\mathtt G}$ 
 as\footnote{Where, for integer vectors
 $\bal=(\alpha_1,\alpha_2),\bbe=(\beta_1,\beta_2)$ we set
 $|\bal+\bbe|:=\alpha_1+\alpha_2+\beta_1+\beta_2$.
}
	\begin{equation}\label{G}
			{\mathtt G}(\mathbf z,\bar{\mathbf  z})
						=
				\sum_{i+j=4}\frac{f_{i,j}}{4(\sqrt{\omega_-})^{i}
				(\sqrt{\omega_+})^{j}}
				(z_1+\bar z_1)^i
				(z_2+\bar z_2)^j
			=
				\sum_{|\bal+\bbe|=4 }
				{\mathtt G}_{\bal,\bbe} \mathbf z^\bal \bar{\mathbf z}^\bbe
	\end{equation}
where	
\begin{equation}\label{Gab}
				{\mathtt G}_{\bal,\bbe}:=
				\frac{f_{\alpha_1+\beta_1,\alpha_2+ \beta_2}}
				{4(\sqrt{\omega_-})^{\alpha_1+\beta_1}
				(\sqrt{\omega_+})^{\alpha_2+\beta_2}}
				\frac{(\alpha_1+\beta_1)!}{\alpha_1!\beta_1!}
				\frac{(\alpha_2+\beta_2)!}{\alpha_2!\beta_2!}\,.
	\end{equation}
Note that ${\mathtt G}_{\bal,\bbe}=
{\mathtt G}_{\bbe,\bal}\in\mathbb R$.

\subsection{Resonant BNF}\label{Sec:BNF}

The aim of the BNF is to construct a symplectic change
of variables that ``simplifies'' the Hamiltonian
$\mathtt H$ in \eqref{hamiltonian}.
First note that a Hamiltonian $H$ depending 
only on 	$|z_1|^2$ and $|z_2|^2$
writes 
$H=\sum_{\bal} H_{\bal,\bal}
|\mathbf  z|^{2\bal}$
and is integrable; in particular $|z_1|^2$ and $|z_2|^2$
 are constants of motion.
In light of the above considerations
 we guess if it is possible to find,
in a sufficiently small neighborhood of the 
origin
\begin{equation}\label{birbe}
\|\mathbf z\|\leq \epsilon\,,
\end{equation}
 a close-to-the-identity symplectic transformation that ``integrates''  	
	$\mathtt H$ up to terms of degree 6 in $(\mathbf z,\bar{\mathbf  z})$,
	which are smaller.	
	This amounts to transform 
	 $\mathtt H$ into
	${\mathtt N}+\bar{\mathtt H}_4+O(\|\mathbf z\|^6)$, with
	\begin{equation}\label{mozart}
		\bar{\mathtt H}_4:=
		\sum_{ |\bal|=2} {\mathtt G}_{\bal,\bal} |\mathbf  z|^{2\bal}
		=
		{\mathtt G}_{(2,0),(2,0)} 
		|z_1|^4
			+
		{\mathtt G}_{(1,1),(1,1)}
		|z_1|^2 |z_2|^2
		+
		{\mathtt G}_{(0,2),(0,2)} 
		|z_2|^4\,,
	\end{equation}
where, recalling \eqref{Gab},
	\begin{eqnarray}
				{\mathtt G}_{(2,0),(2,0)} 
			&=& 
				\frac{3 f_{4,0}}{2 \omega_-^2}
				=\frac{3}{8\omega_-^2}
				\Big(
				(\phi_1^-)^4M_3+(\phi_2^-)^4N_3
				\Big)\,,
				\nonumber
				\\
				{\mathtt G}_{(1,1),(1,1)} &=&
				\frac{f_{2,2}}{\omega_- \omega_+}
				=\frac3{2\omega_- \omega_+}\Big(
				(\phi_1^-)^2(\phi_1^+)^2M_3+(\phi_2^-)^2(\phi_2^+)^2 N_3
				\Big)\,,
				\nonumber
				\\
				{\mathtt G}_{(0,2),(0,2)} 
			&=&
				\frac{3 f_{0,4}}{2 \omega_+^2}
				=\frac{3}{8\omega_+^2}
				\Big(
				(\phi_1^+)^4M_3+(\phi_2^+)^4N_3
				\Big)\,.
	\label{calippo}			
	\end{eqnarray}
As well known, this is possible if the nonresonance  condition
$\omega_+ k_1+\omega_- k_2\neq 0$ is satisfied for every couple
of integers $k_1,k_2$ with $|k_1|+|k_2|=4$ {\sl and}
$\epsilon$ is small enough.
It is simple to show (see, e.g. Proposition 1 in \cite{DL})
that 
$$
\min_{|k_1|+|k_2|=4}|\omega_+ k_1+\omega_- k_2|\geq 
\min\{\omega_-, \omega_+-\omega_-, |3 \omega_--\omega_+|\}\,.
$$
While, by hypothesis, $\omega_-, \omega_+-\omega_->0$,
 $\sigma=\omega_+ - 3\omega_-$ (introduced in \eqref{def:sigma})
could be zero or small.
It turns out that there exists a  constant $C_1$
(see \cite{DL} for a proof and the evaluation of $C_1$)  such that, if
\begin{equation}\label{peperoni1}
 \epsilon \leq C_1 \sqrt{|\sigma|}\,,
\end{equation} 
 then it is possible to construct
a symplectic transformation putting $\mathtt H$
in (complete) BNF up to order 4, namely
${\mathtt N}+\bar{\mathtt H}_4+O(\|\mathbf z\|^6)$.
 Otherwise, if $|\sigma|$ is too small with respect to\footnote{In
 particular we can assume that 
 $|\sigma|\leq\min\{\omega_-, \omega_+-\omega_-\}$.}
 $\epsilon$, namely if 
$\epsilon > C_1 \sqrt{|\sigma|}$, but 
 $\epsilon$ still satisfies  
 a suitable (weaker\footnote{With $C_2>C_1 \sqrt{|\sigma|}$.}) smallness condition $\epsilon\leq C_2$,
  only a {\sl resonant} BNF is available. 
This means that,
in the case
\begin{equation}\label{peperoni2}
 C_1 \sqrt{|\sigma|}\leq \epsilon \leq C_2\,,
\end{equation}
through   a symplectic transformation,
  the Hamiltonian takes the form
 ${\mathtt N}+\bar{\mathtt H}_{4,{\rm res}}+O(\|\mathbf z\|^6)$,
		where
	\begin{equation}\label{H4RES}	
		\bar{\mathtt H}_{4,{\rm res}}
		:=
			\bar{\mathtt H}_4+
	\mathtt G_{(0,1),(3,0)}z_2\bar z_1^3+
	\mathtt G_{(3,0),(0,1)}z_1^3\bar z_2
	\stackrel{\eqref{Gab}}=
	\bar{\mathtt H}_4+
	\frac{f_{3,1}}
				{4(\sqrt{\omega_-})^{3}
				\sqrt{\omega_+}}
	(z_2 \bar z_1^3+z_1^3\bar z_2)
	\,.
	\end{equation}

\begin{rem}\label{restopaper1}
The construction of the above symplectic transformation in the resonant
case was given in \cite{DL}, where the remainder 
$O(\|\mathbf z\|^6)$ was explicitly estimated.
This means that we found a concrete constant $c_*$ depending
on the parameters such that $O(\|\mathbf z\|^6)\leq c_* \epsilon^6$.
\end{rem}
\begin{rem}\label{belgiomagari}
 $\epsilon$ introduced in \eqref{birbe}
is simply related to $\varepsilon$ introduced in
\eqref{def:sigma} by the change of variables \eqref{crostini},
\eqref{islanda},  \eqref{islanda2}. 
This means that there exist two constants
$\underline c<\bar c$ such that $\underline c\leq \epsilon/\varepsilon\leq
\bar c$. 
Then \eqref{peperoni1} and \eqref{peperoni2}
justify \eqref{def:sigma} and \eqref{def:sigma2},
respectively. 
\end{rem}


We now introduce action-angle variables\footnote{$\mathbb T:=\mathbb R/{2\pi \mathbb Z}$, $\mathbb T^2:=\mathbb R^2/{2\pi \mathbb Z^2}$.}
$(\mathbf I,\bphi)=(I_1,I_2,\varphi_1,\varphi_2)\in \mathbb R^2\times \mathbb T^2$
through the transformation
	\begin{equation}\label{amazon}
		 z_j=\sqrt {I_j} e^{-\mathrm i \varphi_j}\,,\qquad
			I_j>0\,,
			\qquad j=1,2.
	\end{equation}
	
\begin{rem}
 Note that the above map is singular at 
	$\mathrm  z_1$ or $\mathrm  z_2=0$ 
and is defined for $I_1,I_2>0$.
\end{rem}
\noindent
In the symplectic variables in \eqref{amazon}
 the truncated Hamiltonians
 ${\mathtt N}+{\mathtt H}_4$ and
${\mathtt N}+\bar{\mathtt H}_{4,{\rm res}}$
take the final forms
\begin{eqnarray}
\hat{\mathcal H}_{\rm res}(\bI,\bphi)
					&:=&
					\hat{\mathcal H}(\bI,\bphi)
					+
						\frac{f_{3,1}}
						{2(\sqrt{\omega_-})^{3}
						\sqrt{\omega_+}}
						\sqrt{I_1^3I_2}
						\cos (\varphi_2-3\varphi_1)\,,
							\label{polloRES}
					\,,
					\\
\hat{\mathcal H}(\bI,\bphi)
					&:=&\omega_-I_1+\omega_+ I_2
					+\mathcal H_4(\bI)
					\label{HamAARES}
							\\
					\label{pollo}
						\mathcal H_{4}(\bI)
					&:=&
						{\mathtt G}_{(2,0),(2,0)} 
						I_1^2
						+
						{\mathtt G}_{(1,1),(1,1)}
						I_1 I_2
							+
						{\mathtt G}_{(0,2),(0,2)} 
						I_2^2\,.
\end{eqnarray}


The frequencies of the integrable
 nonresonant truncated Hamiltonian
$ \hat{\mathcal H}$ in \eqref{HamAARES}
			 are
			 the derivatives
of the energy with respect to the actions, namely, by \eqref{pollo},
\begin{eqnarray}\label{omega_nonlineare}
			\omega_-^{\rm nl}
&:=&\partial_{I_1} \hat{\mathcal H}=	
\omega_-+
			2{\mathtt G}_{(2,0),(2,0)}I_1
			+
	{\mathtt G}_{(1,1),(1,1)}I_2
				\,,
				\nonumber
					\\
			\omega_+^{\rm nl}
			&:=&\partial_{I_2} \hat{\mathcal H}
			=
			\omega_++2{\mathtt G}_{(0,2),(0,2)}I_2+
{\mathtt G}_{(1,1),(1,1)}I_1
				\,,
				\nonumber
	\end{eqnarray}	
In particular, when 
$M_3=0$, by \eqref{calippo} we 
have
\begin{eqnarray}\label{formulaSL1}
			\omega_-^{\rm nl}
&=&		
\omega_-+N_3\left(
	\frac{3}{8\omega_-}
				(\phi_2^-)^4 a_-^2
				+
				\frac3{4\omega_-}
			(\phi_2^-)^2(\phi_2^+)^2
				a_+^2
				\right)
				\,,
				\nonumber
					\\
			\omega_+^{\rm nl}
			&=&
			\omega_+ +N_3\left(
	\frac{3}{8\omega_+}
				(\phi_2^+)^4
				a_+^2
				+
				\frac3{4 \omega_+}(\phi_2^-)^2(\phi_2^+)^2 a_-^2
				\right)\,,
	\end{eqnarray}
where $a_-, a_+>0$
are the initial amplitudes.
Note that in the original variables $q_1$ and $q_2$, 
one has
\begin{equation}\label{cerniera}
q_1(0)=a_-\,,\quad 
q_2(0)=a_+\,,\quad
\dot p_1(0)=\dot p_2(0)=0\,,
\end{equation}
that correspond, by \eqref{islanda}, \eqref{islanda2} and
\eqref{amazon}, in initial 
action-angle variables: 
\begin{equation}\label{ampiezze}
I_1(0)=\frac12\omega_- a_-^2\,,\qquad
I_2(0)=\frac12\omega_+ a_+^2\,,\qquad
\varphi_1(0)=\varphi_2(0)=0\,.
\end{equation}
Formula \eqref{formulaSL1}
was already known (see \cite{SW23jsv} or \cite{DL}),
but it
does not hold close to resonances.
To obtain the analogous of formula  
\eqref{formulaSL1} in the resonant case is much
more complicated since one has to 
integrate the Hamiltonian
$\hat{\mathcal H}_{\rm res}$ in \eqref{polloRES}.
This is exactly what we are going to do in the following 
sections. The analogous of \eqref{formulaSL1} in the resonant case
are the formula 
\eqref{californiaQUATER}-\eqref{sublime2} 	below.


\begin{rem}[Reversibility]\label{reverso}
 Since the Hamiltonian $\hat{\mathcal H}_{\rm res}(\bI,\bphi)$
in \eqref{polloRES} is even in $\bphi$ the system is reversible, namely
 if $\big(\bI(t),\bphi(t)\big)$ is a solution
 the same holds true for $\big(\bI(-t),-\bphi(-t)\big)$.
 In particular if $\bphi(0)=0$ the solution is
 even in the actions and odd in the angles, namely
 $\bI(t)=\bI(-t)$ and  $\bphi(t)=-\bphi(-t)$.
\end{rem}

\subsection{The slow angle and the effective Hamiltonian}
\label{SubSec:effective Hamiltonian}

 It is convenient to introduce the 
adimensional effective Hamiltonian
$F$ depending solely on one angle
$\psi_1$,
namely the ``slow angle''.
Let us consider the canonical 
transformation 
		\begin{equation*}
				\Phi_*:\mathbb R^2\times \mathbb T^2\to \mathbb R^2\times \mathbb T^2\qquad \Phi_*(J,\psi)=(I,\bphi):=( \mathcal M^T J,  \mathcal M^{-1}\psi)\qquad \mathcal M=
		\left(
			\begin{array}{cc}
			-3&1\\
				1&0
			\end{array}
		\right)
		\end{equation*}
		so that
	\begin{equation}\label{salmone}
	\left\{
			\begin{array}{l}
				I_1=J_2-3 J_1\\
				I_2=J_1\,,
			\end{array}
			\right.
			\qquad 
			\left\{
			\begin{array}{l}
			\varphi_1=\psi_2\\
			\varphi_2=\psi_1+3\psi_2
				\,.
			\end{array}
			\right.
	\end{equation}
Note that $ \mathcal M$ has integer entries 
and $\mathrm{det} \mathcal M=-1$
 so that the inverse $ \mathcal M^{-1}$ has also integer entries.
 This implies that $\psi= \mathcal M\varphi$ and its inverse
 $\varphi=\mathcal M^{-1}\psi$
are well defined on the 
torus
 $\mathbb T^2$. 
 Note also that, by \eqref{amazon},
 we have
	 \begin{equation}\label{amazon2}
			J_2>3J_1\,, \qquad J_1>0\,. 
	\end{equation}
Let us write $\hat{\mathcal H}_{{\rm res}}$
in \eqref{polloRES} in the $(J,\psi)$-variables
\begin{equation}\label{secularRES}
	\hat{\mathbb H}_{\rm res}(J,\psi_1):=
	\hat{\mathcal H}_{{\rm res}}\big(\Phi_*(J,\psi_1)\big):=
			\omega_- J_2+\sigma J_1
			+\mathbb H_{4,{\rm res}}(J,\psi_1)\,,
	\end{equation} 
			where
	\begin{equation}\label{connessione}
			\mathbb H_{4,{\rm res}}(J,\psi_1)
			:=\mathcal H_{4}(J_2-3J_1,J_1)
			+
			\frac{f_{3,1}}
				{2(\sqrt{\omega_-})^{3}
				\sqrt{\omega_+}}
			\sqrt{(J_2-3J_1)^3 J_1}
			\cos (\psi_1)\,.
	\end{equation}
	Note that
	$\hat{\mathbb H}_{\rm res}$
	is reversible in the sense of Remark \ref{reverso}.
	Moreover it
	depends only on the ``slow angle'' $\psi_1$,
	that evolves by a small frequency 
	$\sigma+O(|J|)\sim0$ (recall \eqref{def:sigma2}),
	but
does not depend on the ``fast angle'' $\psi_2$,
that, on the contrary,  evolves by a
  frequency $\omega_->0$, which is definitively 
  different from zero.
So the partial derivative w.r.t. $\psi_2$
of $\hat{\mathbb H}_{\rm res}$
vanishes and, by the Hamilton's equations,
$\dot J_2=0$, so that $J_2$ is a constant of motion, namely
	$$
		J_2(t)=J_2(0)=:\Ju\,.
	$$
Moreover the fast angle $\psi_2$ simply evolves as
 $\psi_2(t)=\omega_- t+\psi_2(0)$.
 It remains to study the evolution of the $(J_1,\psi_1)$
 variables.

Being $J_2$ a constant of motion the dynamic of the ``resonant truncated Hamiltonian''
$\hat{\mathbb H}_{\rm res}$
in \eqref{secularRES}
  is simply generated by the 
one-degree-of-freedom
``effective Hamiltonian''
	 \begin{equation}\label{1DHAM}
			  \mathbb H_{\Ju}(J_1,\psi_1):=
				\sigma J_1
				+\mathbb H_{4,{\rm res}}(J_1,\Ju,\psi_1)
				\,,
	\end{equation}
	with $\mathbb H_{4,{\rm res}}$ defined in \eqref{connessione}.
At this point it is convenient to introduce the ``adimensional Hamiltonian''\footnote{Also
$\chi$ is adimensional.}
		\begin{equation}\label{pluto}
			\hat H(J_1,\psi_1)=
						\hat H_\Ju(J_1,\psi_1):=\frac{1}{\chi \Juq}\mathbb H_{\Ju}(J_1,\psi_1)\,,
						\qquad
						\mbox{with}\ \ 
							\chi:=\frac{f_{3,1}}
							{2\sqrt 3(\sqrt{\omega_-})^{3}
							\sqrt{\omega_+}}\neq 0\,,
		\end{equation}
		and rewrite $\hat H_\Ju$ as a function of the
``adimensional action'' 
\begin{equation}\label{albero}
x:=3J_1/\Ju\qquad \mbox{with}
\quad0<x<1\,, 
\end{equation}
by \eqref{amazon2}.
We have the following
\begin{lem}\label{lem:pollo}
It results that
\begin{equation}\label{loacker}
\hat H(J_1,\psi_1)=
\hat H_\Ju(J_1,\psi_1)
	=F(\psi_1,3J_1/\Ju;\Ju)+ a_0\,,
\end{equation}
where 
	\begin{eqnarray}\label{sublime}
					F(\psi,x)=F(\psi,x;\Ju)
					&:=&a(x;\Ju)+b(x)\cos \psi\,,
					\nonumber
					\\
					a(x)=a(x;J_2)
				&:=&
					\frac12 a_2 x^2+ a_1 x\,,					\nonumber
					\\
					b(x)
				&:=&
					\sqrt{(1-x)^3 x}>0\,,\qquad 0< x< 1\,,
					\nonumber
						\\
					a_0
				&:=&
					\frac{{\mathtt G}_{(2,0),(2,0)} }{\chi}\,,
					\nonumber
						\\
					a_1&=&a_1(\Ju)
				:=
					-2\frac{{\mathtt G}_{(2,0),(2,0)} }{\chi}+
					\frac{{\mathtt G}_{(1,1),(1,1)} }{3\chi}
					+\frac{\omega_+-3\omega_- }{3 \Ju\chi}\,,
					\nonumber
					\\
					a_2
				&:=&
					2\frac{{\mathtt G}_{(2,0),(2,0)} }{\chi}
						-
							2\frac{{\mathtt G}_{(1,1),(1,1)} }{3\chi}
									+
							2\frac{{\mathtt G}_{(0,2),(0,2)} }{9\chi}\,.
		\end{eqnarray}
\end{lem}	
\proof 
Multiplying the right hand side of \eqref{loacker}
by $\chi\Juq$
 we have, by \eqref{sublime}, 
 \begin{eqnarray*}
&& 	\chi\Juq \,F(\psi_1,3J_1/\Ju)+
\chi\Juq \, a_0
=
\chi\Juq \,a(3J_1/\Ju)
+
\chi\Juq \, a_0+ \chi\Juq \,b(3J_1/\Ju)\cos \psi_1
\\
&&\qquad
=\frac92 \chi a_2 J_1^2+3 \chi a_1 J_1\Ju
+
\chi\Juq \, a_0
+\chi\sqrt{3(\Ju-3J_1)^3J_1}
\cos\psi_1
\\
&&\qquad
={\mathtt G}_{(2,0),(2,0)} (\Ju-3J_1)^2
+{\mathtt G}_{(1,1),(1,1)} J_1 (\Ju-3J_1)
+{\mathtt G}_{(0,2),(0,2)} \Juq
\\
&&\qquad \quad
+\sigma J_1
+\chi\sqrt{3(\Ju-3J_1)^3J_1}
\cos\psi_1
\\
&&\qquad
\stackrel{\eqref{pollo}}=
\mathcal H_{4}(J_2-3J_1,J_1)
+\sigma J_1
+\chi\sqrt{3(\Ju-3J_1)^3J_1}
\cos\psi_1
\\
&&\quad \ \stackrel{\eqref{connessione},\eqref{pluto}}=
\mathbb H_{4,{\rm res}}(J,\psi_1)
+\sigma J_1
\\
&&\qquad\stackrel{\eqref{1DHAM}}=
\mathbb H_{\Ju}(J_1,\psi_1)
\stackrel{\eqref{pluto}}=
\chi\Juq 
\hat H_\Ju(J_1,\psi_1)\,,
\end{eqnarray*}
proving \eqref{loacker}.
\eproof

\begin{rem}
 Note that $F(\psi,x;\Ju)$ depends on $\Ju$ only through 
 $a(x;\Ju)$, which depends on $\Ju$ only through 
 $a_1(\Ju)$. Moreover at the exact resonance $\omega_+=3\omega_-$
 the dependence on $\Ju$ disappears.
\end{rem}

\medskip	
\noindent
Since $\psi_2$ does not appear in $\hat H$,
the conjugated action $\Ju$ is a constant of motion	
and $\hat H$ is actually a one degree of freedom 
Hamiltonian system depending on $\Ju$ as a parameter.
		From now on we  consider the one degree
of freedom Hamiltonian
	$
		\hat H(J_1,\psi_1)=\hat H_\Ju(J_1,\psi_1)
	$
on the phase space 
	$
		(0,3\Ju)\times \mathbb T\ni (J_1,\psi_1)
	$
with 
	$
		\mathbb T:=\mathbb R/2\pi \mathbb Z
	$.

\section{The phase portrait}\label{sec:phase}

\begin{figure}[h]
		\center
		\includegraphics[width=6cm,keepaspectratio]{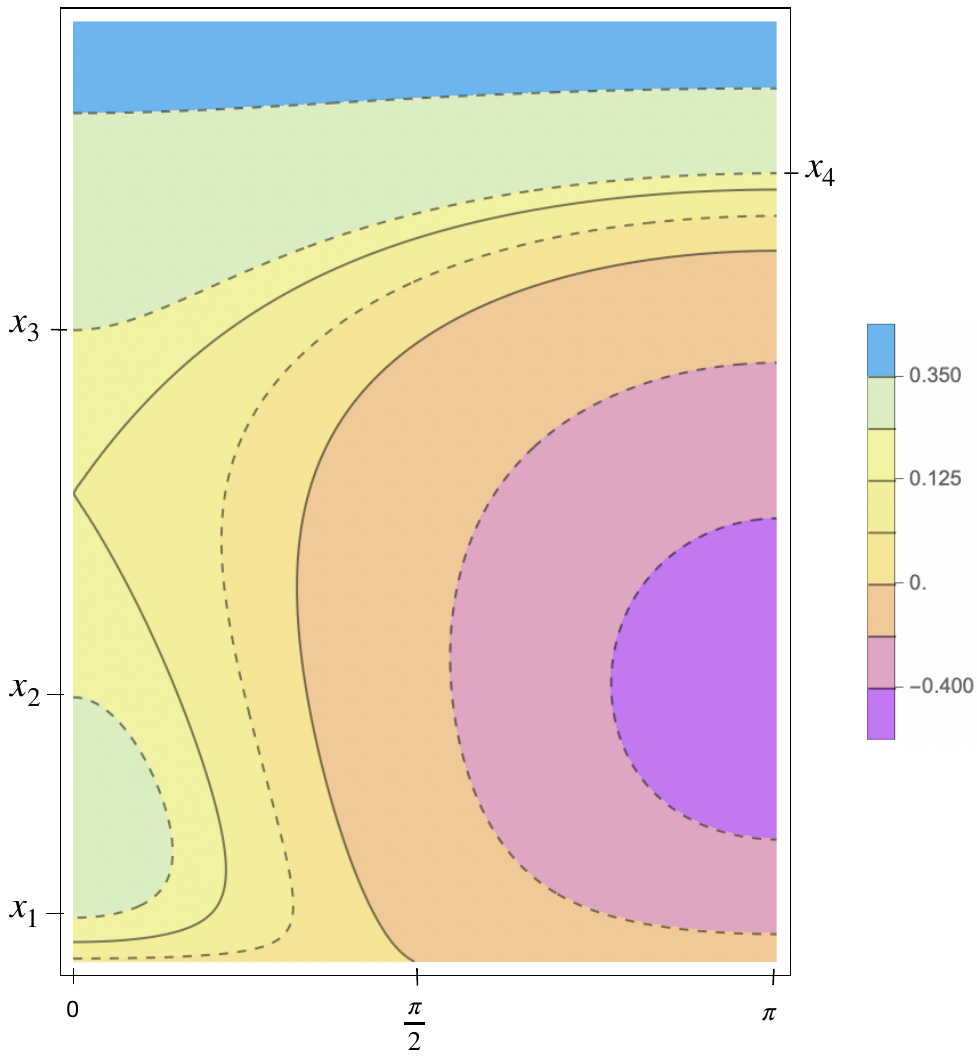}
		\caption{The zone 
		$Z_{21}^+\cap\{0\leq \psi\leq \pi\}$ for $a_1=1$ and $a_2=3$ filled by the level curves of $x-\frac 32 x^2+b(x)\cos \psi=E$, for different values
			of the energy $E$.
	The intersections $x_1,x_2,x_3$, respectively $x_4$, of the two 
	level curves of energy $E=0.16$ with
	the line $\psi=0$, respectively  $\psi=\pi$,
	are shown.  }
	\label{4_soluzioni}
	\end{figure}

In this section we study the
phase portrait of 
the adimensional Hamiltonian
$\hat H$ in \eqref{pluto}
describing level curves, critical points and extrema.
An important remark, that simplifies the treatment, 		
is the fact that, thanks to 
 \eqref{loacker}, $\hat H$
has, up to the rescalings $x=3J_1/J_2$ and $F=\hat H+a_0$,
  the same 
  level curves, critical points and extrema as 
 the auxiliary function
$F$ in \eqref{sublime}.
  Such objects
 are studied in Subsections
 \ref{sec:cc}, \ref{sec:ex}
 and \ref{sec:lc}, respectively.
As usual, the new action coordinates,
 that integrate the system, 
 are defined as the areas enclosed by the level curves.
  In order to evaluate them
 its important
 to  determine the intersections
 between the level curves and the lines
 $\{\psi=0\}$ and  $\{\psi=\pi\}$
 since they appear as endpoints of 
 the involved integrals.
 It turns out that such intersections
 correspond to the real roots
 of the quartic polynomial 
  $\mathbf P(x)$, see \eqref{penisola}
and
Figure \ref{4_soluzioni}.
As the energy $E$ varies, 
it is necessary to distinguish 
whether $\mathbf P$ has
$4,2$ or $0$ real 
roots\footnote{Note that, 
excluding
the degenerate case of multiple roots,
the number of real roots
is even.} and
whether  a root corresponds to
an intersection with
$\{\psi=0\}$ or $\{\psi=\pi\}$.
Explicit formulae for the roots
are given in Subsection \ref{Sec:radici},
see Figure \ref{xPiuMenoSoluzioni}.
In Subsections
\ref{sec:disney} and 
\ref{6zones}
as the parameters vary,
 six topologically different
scenarios appear.

 \subsection{Critical points, elliptic and hyperbolic zones}
 \label{sec:cc}

We now describe how the critical point of $\hat H$
 depends on the values of the parameters
$a_2$ and $a_1$ in \eqref{sublime}.
First we note that 
$(J_1,\psi_1)$ is a critical point of $\hat H$
if and only if $(\psi,x):=(\psi_1,3J_1/J_2)$
is a critical point of the auxiliary function
$F(\psi,x)=a(x)+b(x)\cos \psi$
 defined on $\mathbb T
 \times(0,1) $ (recall \eqref{sublime}).
 Moreover the nature of a critical point
 (maximum, minimum or saddle) is the same
 for $\hat H$ and $F$.
 Then in the following we will study critical points of $F$
 as the parameters $a_2$ and $a_1$ vary.

It is immediate to see that, since
	$$
		\partial_x F(\psi,x)=a'(x)+b'(x)\cos \psi\,,
		\qquad
		\partial_\psi F(\psi,x)=-b(x)\sin \psi
	$$
and $b(x)>0$,
 the critical points of $F$ have the form
 $(x,0)$ with $a'(x)+b'(x)=0$ or $(x,\pi)$
 with $a'(x)-b'(x)=0$.
 Namely
 	\begin{eqnarray}
		\nabla F(0,x)=0 \ &\iff&\
			 -a_2 x-a_1=b'(x)\,,
			 \label{gingerina1}
				 \\
		 \nabla F(\pi,x)=0 \ &\iff&\
			 a_2 x+a_1=b'(x)\,
			 \label{gingerina2}
	\end{eqnarray}
	where 
	 \begin{equation}\label{petrolio}
 				b'(x)=\frac{(1-4x)\sqrt{1-x}}{2\sqrt x}.
	 \end{equation}
The number of solutions of equations
\eqref{gingerina1},\eqref{gingerina2} depends
on the parameters $a_1,a_2$.

\noindent

 Set
 \begin{equation}\label{gigi}
g(a_1):=\frac{1}{27}
\textstyle
\big(\sqrt{9+4a_1^2}-2 a_1\big)
\big(9-4a_1^2-4a_1\sqrt{9+4a_1^2}\big)
\end{equation}
and\footnote{Note that $g(a_1)>-a_1$.}
\begin{eqnarray}
Z_{10}
&:=&\{(a_1,a_2)\ \ :\ \ a_2<-g(-a_1)\}\,,
\nonumber
\\
Z_{12}
&:=&\{(a_1,a_2)\ \ :\ \ -g(-a_1)<a_2<-a_1\}\,,
\nonumber
\\
Z_{21}
&:=&\{(a_1,a_2)\ \ :\ \ -a_1<a_2<g(a_1)\}\,,
\nonumber
\\
Z_{01}
&:=&\{(a_1,a_2)\ \ :\ \ a_2>g(a_1)\,.
\label{kiwi}
\end{eqnarray}

\begin{figure}
\center
\includegraphics[width=12cm]{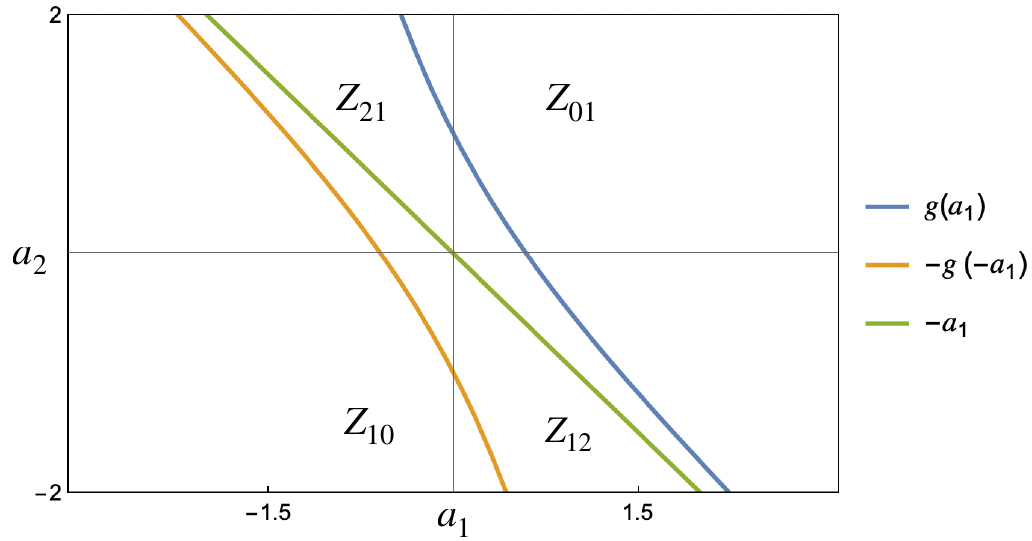}
		\caption{The four zones
		$Z_{10},Z_{12},Z_{21},Z_{01}$ in the space $(a_1,a_2)$. See formula \eqref{kiwi}.}
	\end{figure}

In particular the following result holds
\begin{pro}\label{pandivia}
If $(a_1,a_2)\in Z_{ij}$ then $F(0,x)$ has $i$ critical points
 and $F(\pi,x)$
$j$ critical points.
More precisely:
\begin{itemize}
\item If $(a_1,a_2)\in Z_{10}$ then $F(0,x)$ has a positive maximum
at some  $x_1^{(0)}$ and $F(\pi,x)$ is strictly decreasing;

\item If $(a_1,a_2)\in Z_{01}$ then $F(\pi,x)$ has a 
negative minimum
at some  $x_1^{(\pi)}$
 and $F(0,x)$ is strictly increasing;

\item If $(a_1,a_2)\in Z_{12}$ then $F(0,x)$ has a
positive maximum
at some  $x_1^{(0)}$, while $F(\pi,x)$ has
a negative minimum  at
some $x_1^{(\pi)}$
and a maximum
at some $x_2^{(\pi)}$, with
$x_1^{(\pi)}<x_2^{(\pi)}$;

\item If $(a_1,a_2)\in Z_{21}$ then $F(0,x)$ has
 a positive maximum at
some $x_1^{(0)}$ and  a minimum
 at some $x_2^{(0)}$, with
$x_1^{(0)}<x_2^{(0)}$,
while $F(\pi,x)$ has 
 a negative minimum
at some  $x_1^{(\pi)}$.
\end{itemize}
As a corollary,
if $(a_1,a_2)\in Z_{ij}$ then $F$ has $i$ critical points
of the form $(0,x)$ and 
$j$ critical points
of the form $(\pi,x)$.
More precisely:
\begin{itemize}
\item If $(a_1,a_2)\in Z_{10}$ then $F$ has a positive maximum
at   $(0,x_1^{(0)})$;
\item If $(a_1,a_2)\in Z_{01}$ then $F$ has a negative  minimum
at   $(\pi,x_1^{(\pi)})$;
\item If $(a_1,a_2)\in Z_{12}$ then $F$ has a positive  maximum
at   $(0,x_1^{(0)})$,
a negative  minimum  at
 $(\pi,x_1^{(\pi)})$
and a saddle
at  $(\pi,x_2^{(\pi)})$;
\item If $(a_1,a_2)\in Z_{21}$ then $F$ has
 a positive  maximum at
 $(0,x_1^{(0)})$,  a saddle at  $(0,x_2^{(0)})$ and
 a negative minimum
at  $(\pi,x_1^{(\pi)})$.
\end{itemize}
\end{pro}
\proof
See Appendix.
\eproof

\begin{figure}
\center 
{\includegraphics[width=7cm,keepaspectratio]{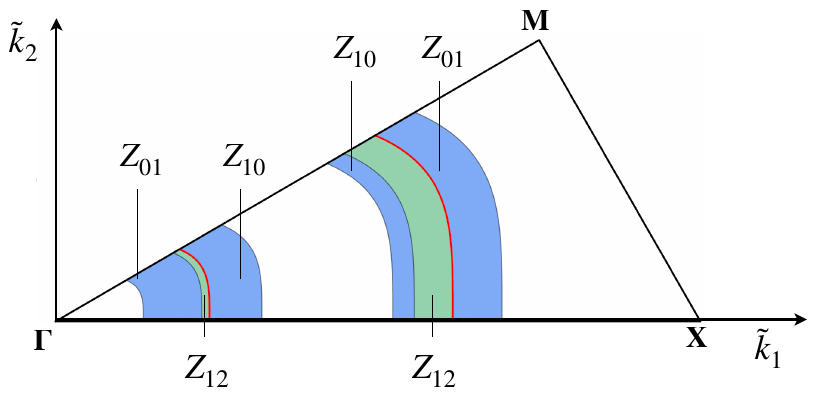}
\qquad
\includegraphics[width=5cm,keepaspectratio]{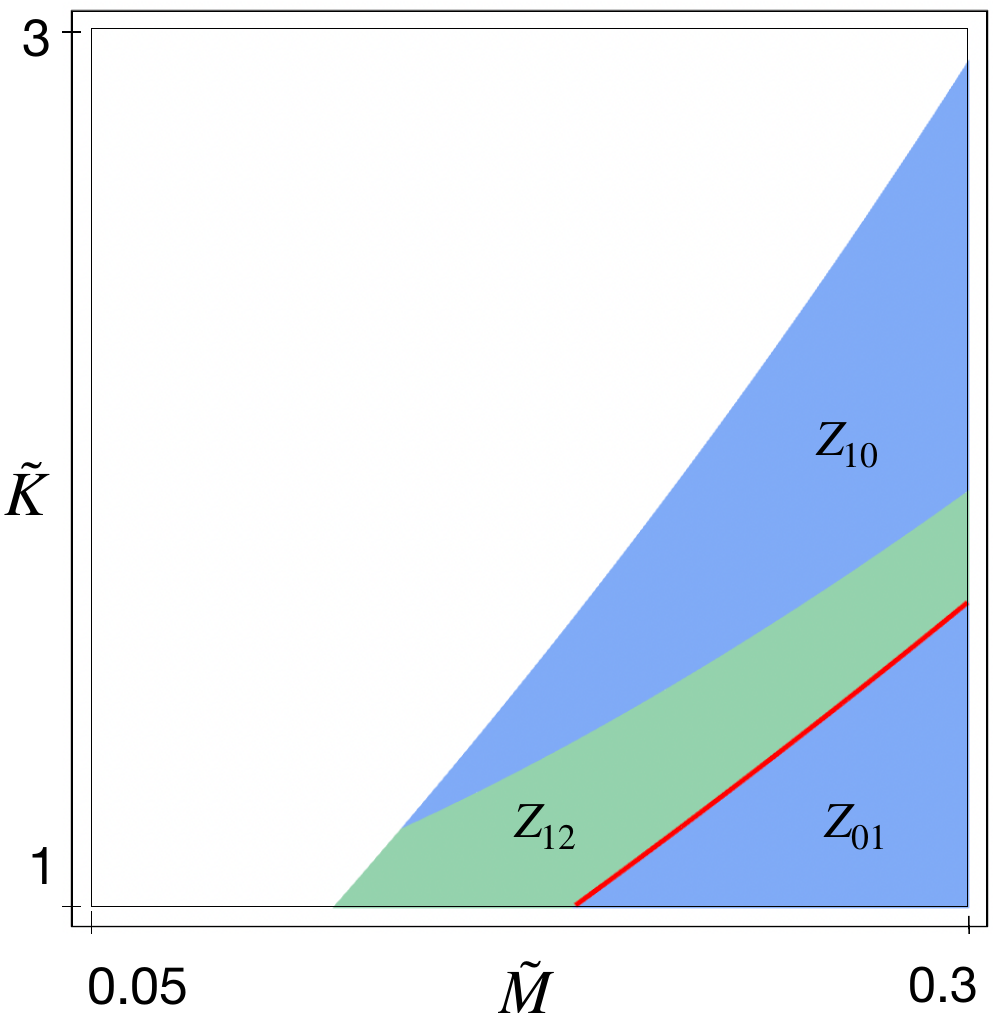}}
\caption{(On the left) For $\tilde M=0.2$, $\tilde K=1.1$,
$N_3=-10^{4}$,
one can take $C_1$ and $C_2$ in \eqref{peperoni2} as 
$C_1=5\times 10^{-4}$ and $C_2=2.5\times 10^{-3}$.
Choose $\epsilon=6.5\times 10^{-4}$. Then the smallness
condition  
\eqref{peperoni1} for the nonresonant BNF is
satisfied only   when $|\sigma|\geq 1.7$.
On the other hand for $|\sigma|< 1.7$ only a resonant BNF is available. The first condition
is satisfied only by wave numbers $(\tilde k_1,\tilde k_2)$
in the white regions in the  Brillouin triangle. On the contrary,
in the green regions we have a resonant BNF 
 of  type $Z_{12}$
  while, in the blue regions,
we have a resonant  BNF 
 of  type $Z_{01}$ or $Z_{10}$.
Actually, close to the red curves, representing the exact 3:1 
resonance, there is also a very tiny strip (of width $10^{-6}$, not shown in the figure)
corresponding to a BNF of type $Z_{21}$.
(On the right) Here we fix $\tilde k_1=2.58$, $\tilde k_2=0$,  
 and let 
 $(\tilde M,\tilde K)$ vary in the rectangle 
 $[0.05,0.3]\times[1,5]$.}
\label{clarissa}
\end{figure}

\medskip
We call $Z_{21}$, $Z_{12}$ \emph{hyperbolic zones},
since they contain hyperbolic equilibria,
 and $Z_{01}$, $Z_{10}$ \emph{elliptic zones},
since they contain only elliptic equilibria.
For any fixed pair $(\tilde M,\tilde K)$, 
it is possible to identify which wave numbers
$(\tilde k_1,\tilde k_2)$ in the Brillouin triangle
give rise to resonant normal forms with  different 
phase portraits. In particular if the corresponding
values of $a_1$ and $a_2$ belong to 
$Z_{21}$, $Z_{12}$, then the phase portrait
contains one hyperbolic and two elliptic equilibria, while, 
for $a_1$ and $a_2$ belonging to 
$Z_{01}$, $Z_{10}$
only elliptic equilibria appear (see Figure \ref{clarissa}).
For brevity we denote by \emph{BNF of type $Z_{ij}$}
the corresponding Birkhoff Normal Form.

\begin{figure}[h!]
			\centering		\includegraphics[width=8cm,height=8cm,keepaspectratio]{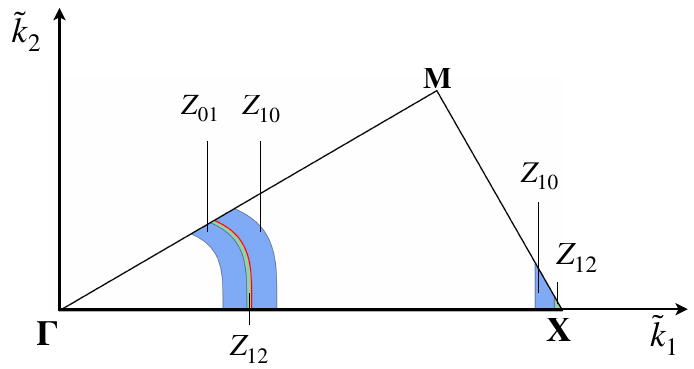}
			\caption{ For $\tilde M=0.146$, $\tilde K=5.73$,
$N_3=-10^{4}$,
one can take $C_1$ and $C_2$ in \eqref{peperoni2} as 
$C_1=9\times 10^{-4}$ and $C_2=4.5\times 10^{-3}$.
Choose $\epsilon=1.2\times 10^{-3}$. Then the smallness
condition  
\eqref{peperoni1} for the nonresonant BNF is
satisfied only   when $|\sigma|\geq 1.78$.
On the other hand for $|\sigma|< 1.78$ only a resonant BNF is available. The first condition
is satisfied only by wave numbers $(\tilde k_1,\tilde k_2)$
in the white regions in the  Brillouin triangle. On the contrary,
in the green regions we have a resonant BNF 
 of  type $Z_{12}$
  while, in the blue regions,
we have a resonant  BNF 
 of  type $Z_{01}$ or $Z_{10}$.
Actually, close to the red curves, representing the exact 3:1 
resonance, there is also a very tiny strip (of width $10^{-6}$, not shown in the figure)
corresponding to a BNF of type $Z_{21}$.
Note that the point $\bf X$ is on an exact 3:1 resonance.
}			\label{Triangolo_Brillouin_zone_iperbolica_ellittica_con_Zone_X}
		\end{figure}

\begin{rem}\label{maniavanti}
 In the following for simplicity we restrict to the case in which 
 $(a_1,a_2)\in  Z_{ij}$ for some $0\leq i,j\leq 2$. 
 This means that
 we avoid the {\sl degenerate} cases $a_2+a_1=0$,
 when $x=1$ is a solution of
 \eqref{gingerina1}-\eqref{gingerina2},
 $a_2-g(a_1)=0$ and $a_2+g(-a_1)=0$, 
 when two solutions coincide. We will briefly discuss
 such degenerate cases in Subsection \ref{degenere}.  
\end{rem}

%
%

\begin{figure}[h!]
	\center
\includegraphics[width=6cm,keepaspectratio]{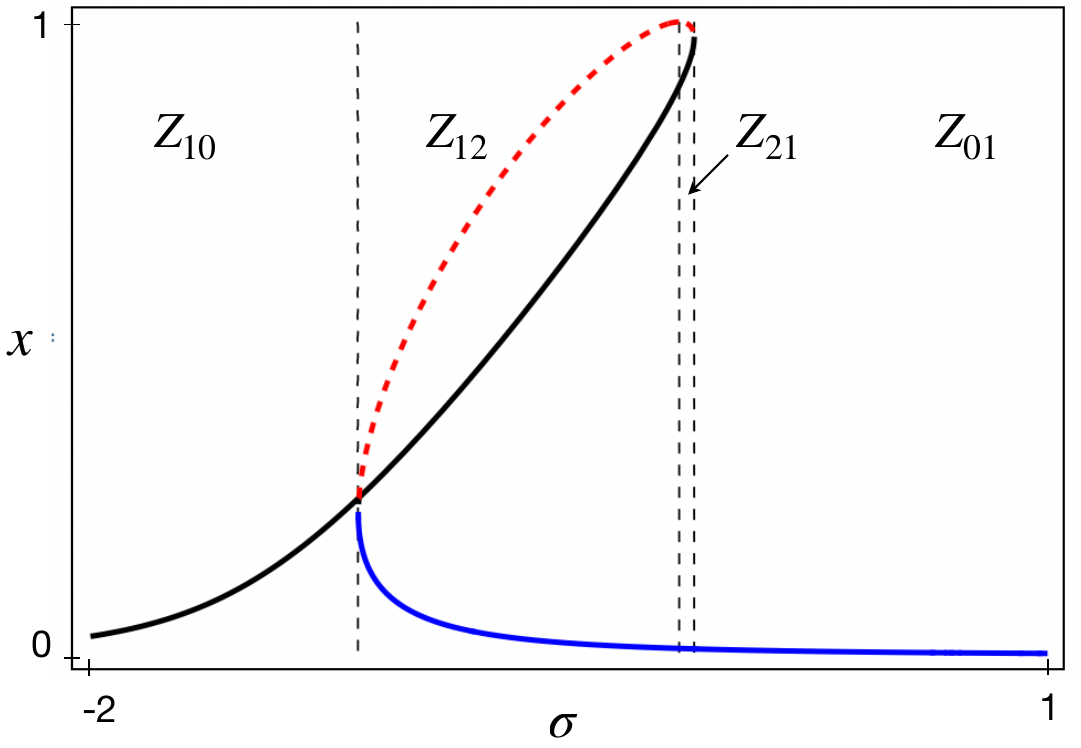}
		\caption{
		Frequency response curve:  critical values of 
		$x$ as function of $\sigma$ and their linear stability.
		The black/blue lines correspond to
		marginally stable equilibria (maxima/minima).
		The red dotted curve corresponds to unstable equilibria (saddles).
		Here $\Ju=10^{-4}, k_1=4 \pi/3, k_2=0, N_3=-10^4, M_3=0, \tilde M=0.15$ and $\tilde K$
		is chosen as a suitable function of $\sigma$;
		more precisely $\tilde K(\sigma)$ is the inverse
		function of 
		 $\tilde K\to(\omega_+-3\omega_-)$.
		The four different zones correspond, from left to right,
		to $Z_{10},Z_{12},Z_{21},Z_{01}$.}
	\end{figure}
%

\subsection{Extrema}\label{sec:ex}

We now discuss the extrema of $F$ 
in \eqref{sublime}
and their dependence on the parameters $a_1,a_2$.
Following the notation of Proposition \ref{pandivia}
we set
\begin{eqnarray}\label{maxmin}
\Emax&:=&F(0,x_1^{(0)})\,,\qquad
{\rm if}\ \ (a_1,a_2)\in Z_{10},\ Z_{12},\ Z_{21}\,,
\nonumber
\\
\Emin&:=&F(\pi,x_1^{(\pi)})
\,,\qquad
{\rm if}\ \ (a_1,a_2)\in Z_{01},\ Z_{12},\ Z_{21}\,,
\nonumber
\\
\Esad&:=&F(0,x_2^{(0)})\,,\qquad
\, {\rm if}\ \ (a_1,a_2)\in Z_{21}\,,
\nonumber
\\
\Esad&:=&F(\pi,x_2^{(\pi)})\,,\qquad
{\rm if}\ \ (a_1,a_2)\in Z_{12}\,.
\end{eqnarray}
Then define
\begin{equation}\label{extrema}
E_+:=\sup_{\mathbb T\times (0,1)} F\,,
\qquad
E_-:=\inf_{\mathbb T\times (0,1)} F\,.
\end{equation}
Since $F(\pi,x)<F(\psi,x)<F(0,x)$ for every $0<x<1$,
$0<\psi<2\pi$, $\psi\neq \pi$,
we have that
$$
E_+:=\sup_{(0,1)} F(0,x)\,,
\qquad
E_-:=\inf_{(0,1)} F(\pi,x)\,.
$$
Note that
$$
E_-<0<E_+\,,
$$
since $F(0,0)=F(0,\pi)= 0$ and 
$F(\pi,x)<0<F(0,x)$ for $x>0$ small enough since
 $\lim_{x\to 0^+}\partial_x F(0,x)=+\infty$
 and $\lim_{x\to 0^+}\partial_x F(\pi,x)=-\infty$.
Note that in the cases $Z_{10}, Z_{12}$
we have $E_+=\Emax$, since the function
$x\to F(0,x)$ has only one critical point (a maximum);
analogously  in the cases $Z_{01}, Z_{21}$
we have $E_-=\Emin$, since the function
$x\to F(\pi,x)$ has only one critical point (a minimum).
Moreover $E_+=a(1)$ in the case $Z_{01}$; indeed
 the function
$x\to F(0,x)=a(x)+b(x)$ has no critical points
then $E_+=\max\{a(0)+b(0),a(1)+b(1)\}=\max\{a(0),a(1)\}$,
moreover $a(x)+b(x)$ is increasing close to
zero since $\lim_{x\to 0^+}\big( a'(x)+b'(x)\big)=+\infty$.
Analogously $E_-=a(1)$ in the case $Z_{10}$.
Finally 
 in the case $Z_{21}$ we have
$E_+=\max\{a(1), \Emax\}$, since the function
$x\to F(0,x)$ has a maximum at 
$x_1^{(0)}$ and  a saddle at  $x_2^{(0)}$
with $x_1^{(0)}< x_2^{(0)}$.
Analogously
in the case $Z_{12}$ we have
$E_-=\min\{a(1), \Emin\}$.

\subsection{Level curves}\label{sec:lc}

Since $F$ is even  with respect to $\psi$
we can reduce to consider the ``half phase space''
$[0,\pi]\times (0,1)$. Take an energy $E_-<E<E_+$
with $E\neq \Emax,\Emin,\Esad$, and consider the level set
$\{F=E\}$. If $(\psi_0,x_0)\in\{F=E\}$, namely
$F(\psi_0,x_0)=E$, since $(\psi_0,x_0)$ is not a critical point
(being $E\neq \Emax,\Emin,\Esad$ and recalling
Proposition \ref{pandivia} and \eqref{maxmin}),
we can locally\footnote{Namely in a 
sufficiently small neighborood of $(\psi_0,x_0)$.} express
$\{F=E\}$ as a curve by the implicit function theorem.  
In particular, in the half phase space
$[0,\pi]\times (0,1)$, we can always express $\psi$
as a function of $x$, indeed the equation
$
		F(\psi,x)=a(x)+b(x)\cos \psi=E
		$
		has the unique solution
	\begin{equation}\label{psi(x)}
		\psi(x)=\psi(x;E;\Ju)=\arccos \left( \frac{E-a(x)}{b(x)} \right)\,.
	\end{equation}
	Since the domain of definition of the $\arccos$
	is $[-1,1]$, the domain of $\psi(x)$ is
	$$
	D:=\{x\in(0,1) \ | \ -b(x)\leq E-a(x)\leq b(x)\}\,.
	$$ 
We now discuss the structure of $D$.
Consider first the case in which 0 is an accumulation point for $D$; then it must be $E=0$.
Indeed, taking the limit for $x\to 0^+$, $x\in D $
in the inequality $-b(x)\leq E-a(x)\leq b(x)$,
we get $E=0$. 
Moreover, when $E=0$,
\begin{equation}\label{biforco}
\lim_{x\to 0^+} \psi(x;0)=
\lim_{x\to 0^+}\arccos\big(-a(x)/b(x)\big)=\arccos (0)=\pi/2\,.
\end{equation}
	\begin{cl}
		1 cannot be 
		an accumulation point for $D$,
		 since we are assuming that
		$a_2+a_1\neq 0$ (recall Remark \ref{maniavanti}).
	\end{cl}
\proof Indeed
		assume, by contradiction, that 
		1 is an accumulation point for $D$.
		Then
		 taking the limit for $x\to 1^-$, $x\in D$,
		in the inequality $-b(x)\leq E-a(x)\leq b(x)$
		we get $E=a(1)$. Substituting
		$E=a(1)$ in the above inequality and dividing by
		$1-x$ we get
	$$
			-\sqrt{x(1-x)}\leq \frac{a(1)-a(x)}{1-x}\leq \sqrt{x(1-x)}\,,
			\qquad \forall x\in[x_0,1)\,.
	$$
Taking again the limit for $x\to 1^-$ we get
$0=a'(1)=a_2+a_1$, which contradicts the assumption
$a_2+a_1\neq 0$.
\eproof

\medskip
\noindent
As a consequence, assuming $E\neq 0$,
 we have that $D$ is a compact set contained in 
 $(0,1)$;
moreover it is not difficult to see that it is formed by
a finite number of  closed intervals (possibly
  isolated points),
whose endpoints satisfy one of the equations
 \begin{equation}\label{isola}
 a(x)-E=\mp b(x)\,.
\end{equation}
This amounts to find the roots of the 
 quartic polynomial
 \begin{equation}\label{penisola}
 \mathbf P(x)=
 (a(x)-E)^2-(b(x))^2=\left(\frac12 a_2 x^2+a_1 x-E\right)^2-(1-x)^3x=0\,,
 \end{equation}
with $0<x<1$.

\begin{lem}\label{ottobre}
If $E$ is not a critical energy for\footnote{Namely
the energy of a critical point of $F$.} $F$, 
the  roots of the  quartic polynomial $\mathbf P(x)$ in \eqref{penisola} 
with $0<x<1$ are simple. 
\end{lem}
\noindent
\proof
By contradiction, if $0<x_0<1$ is a multiple root of $\mathbf P$,
then $\mathbf P(x_0)=\mathbf P'(x_0)=0$.
Write
$$
\mathbf P(x)=
 \big(a(x)-E-b(x)\big)\big(a(x)-E+b(x)\big)\,.
$$
Assume that $a(x_0)-E-b(x_0)=0$, the case $a(x_0)-E+b(x_0)=0$
being analogous. By $\mathbf P'(x_0)=0$ it follows that
\begin{equation}\label{litfiba}
	\mathbf P'(x_0)=
	 \big(a'(x_0)-b'(x_0)\big)\big(a(x_0)-E+b(x_0)\big)=0\,.
\end{equation}
Since $a(x_0)-E-b(x_0)=0$ and $b(x_0)>0$, by \eqref{litfiba}
we get $a'(x_0)-b'(x_0)=0$. This means that 
$(x_0,\pi)$ is a critical point of $F$, which is a contradiction
since $E$ is a not critical energy.
\eproof

\medskip
From now on we will assume that $E$ is not a critical energy of $F$.
We denote the roots of
$\mathbf P(x;E)$ with $0<x<1$
 by $x_i=x_i(E)$ with 
$i\in\{1,2,3,4\}$. We label the roots in increasing order, namely
$x_i<x_{i+1}$.

\subsection{The quartic equation}
\label{sec:quartic}

In studying the solutions of \eqref{isola} (equivalently of
\eqref{penisola}) on $0<x<1$, it is convenient
to consider the real variable $t\in\mathbb R$
and make the substitution 
	$$
		x=\frac{t^2}{1+t^2}\,.
	$$
Since
	$$
		a(x)-E=\frac{1}{(1+t^2)^2}\left[\frac12 a_2 t^4+a_1 t^2(1+t^2)
		-E(1+t^2)^2\right]
		\,,\qquad
		b(x)=\frac{|t|}{(1+t^2)^2}
	$$
and
	$$
		\frac12 a_2 t^4+a_1 t^2(1+t^2)
		-E(1+t^2)^2=
		\left(\frac12 a_2+a_1-E\right)t^4
		+(a_1-2E)t^2-E
		\,,
	$$
the two equations in \eqref{isola} are equivalent to
	\begin{equation}\label{isola2}
		\left(\frac12 a_2+a_1-E\right)t^4
		+(a_1-2E)t^2-E=\mp |t|\,.
	\end{equation}

	\begin{lem}\label{nuovazelanda}
			Let $t_0$ be a root of the polynomial
					 \begin{equation}\label{groenlandia}
						P(t):=\left(\frac12 a_2+a_1-E\right)t^4
							+(a_1-2E)t^2-t-E
					\end{equation}
			and set
					\begin{equation}\label{gricia}
							x_0:=\frac{t_0^2}{1+t_0^2}\,.
					\end{equation}
			If $t_0<0$, resp. $t_0>0$, then
					 $x_0$ solves $F(0,x_0)=E$, resp. $F(\pi,x_0)=E$.
					 Conversely  if $0<x_0<1$ solve the equation
					 in \eqref{isola} with the $\mp$ sign,
					 then $t_0:=\mp x_0/(1-x_0^2)$
					 solves \eqref{isola2} with the $\mp$ sign
					 and, therefore, is a root of $P(t)$. 
	\end{lem}

\noindent
	\proof
		If $P(t_0)=0$ for some $t_0>0$, then $t_0$ satisfies 
		\eqref{isola2} and, therefore \eqref{isola}, with the plus sign.
		As a consequence $F(\pi, x_0)=E$.
		The proof in the case $t_0<0$ is analogous.
	\eproof
\medskip

When $E=\frac12 a_2+a_1$ the polynomial
$P(t)$  reduces to 
$
(a_2+a_1)t^2-t+a_1+a_2/2
$, whose two roots are easily evaluated. Then we can reduce to the case 
$E\neq\frac12 a_2+a_1$ and consider the equivalent monic polynomial
${\rm P}(t):=P(t)/(\frac12 a_2+a_1-E)$, namely
		\begin{equation}\label{polinomio}
				{\rm P}(t)=t^4+\rp t^2+\rq t+\rr,\quad
				\rp:=\frac{a_1-2E}{\frac12 a_2+a_1-E},
				\ \
				\rq:=\frac{-1}{\frac12 a_2+a_1-E},
				\ \ 
				\rr:=\frac{-E}{\frac12 a_2+a_1-E}.
		\end{equation}
The above quartic polynomial is called 
``depressed'' since it is monic and  its third order coefficient vanishes.
Obviously $P(t)$ and ${\rm P}(t)$ have the same roots.
An immediate corollary of Lemma \ref{nuovazelanda}
is the following

\begin{lem}\label{nuovazelanda2}
Fix $E\neq  \frac12 a_2+a_1$.
Let $t_0$ be a root of ${\rm P}(t)$ in \eqref{polinomio}. 
If $t_0<0$, resp. $t_0>0$,
then $x_0$ in \eqref{gricia}
 solves $F(0,x_0)=E$, resp. $F(\pi,x_0)=E$.
 In particular $x_0$ is a root of 
 $\mathbf P(x)$ in \eqref{penisola}.
\end{lem}

\begin{rem}\label{mimi}
 If ${\rm P}(t)$ has four real distinct roots and $E\neq 0$
 (so that $t=0$ is not a root), then
 the number of positive/negative roots
 depends on the sign of $\rr$ defined in \eqref{polinomio}. Indeed, since $\lim_{t\to\pm\infty}
 {\rm P}(t)=+\infty$, the number of
 positive/negative roots is even if $\rr>0$ and odd 
 otherwise.
\end{rem}

\subsection{Finding the roots of the quartic equation}\label{Sec:radici}

Following \cite{CP23}
we find the roots of the quartic polynomial ${\rm P}(t)$
in \eqref{polinomio}.
First set\footnote{Compare formulas (20) and (10) in
 \cite{CP23}.}
$$
p_*:=-\frac{\rp^2+12 \rr}{3}\,,\quad
q_*:=-\frac{2\rp^3-72 \rp \rr+ 27\rq^2}{27}\,,\quad
\Delta:=-4 p_*^3-27 q_*^2\,.
$$
Let us define the positive\footnote{Compare Theorem 8 in
 \cite{CP23}.} number $s_*>0$ as
 \begin{equation}
s_*:=\left\{
			\begin{array}{ll}
	\sqrt[3]{-\frac{q_*}{2}+\sqrt{-\frac{\Delta}{108}}}	
	+\sqrt[3]{-\frac{q_*}{2}-\sqrt{-\frac{\Delta}{108}}}	
	-\frac{2\rp}{3}
	& {\rm if}\ \ \Delta\leq 0\,,	\\
	2\sqrt{-\frac{p_*}{3}}
	\cos\left(
	\frac13 \arccos\left(
	-\frac{q_*}{2}\sqrt{\left(-\frac{3}{p_*}\right)^3}
	\right)
	\right)
	-\frac{2\rp}{3}
		& {\rm if}\ \ \Delta>0		\,.
			\end{array}
			\right.
\end{equation}
Then the roots of ${\rm P}(t)$ are given 
by\footnote{Compare formula (9) in
 \cite{CP23}.}
\begin{equation}\label{radici}
t^\pm_\varsigma:=\frac{-\varsigma \sqrt{s_*}\pm\sqrt{\delta_\varsigma  }}{2}\,,
\qquad
\delta_\varsigma:=\varsigma 2
\rq (s_*)^{-1/2}-2\rp -s_* \,, 
\qquad
\varsigma=\pm \,.
\end{equation}
The number of {\it real} roots of ${\rm P}(t)$
is:

4 if $\delta_\pm>0$,

2 if $\delta_+ \delta_-< 0$,

0 if $\delta_\pm<0$.

\begin{figure}[h!]
 \center
			\includegraphics[width=10cm]
			{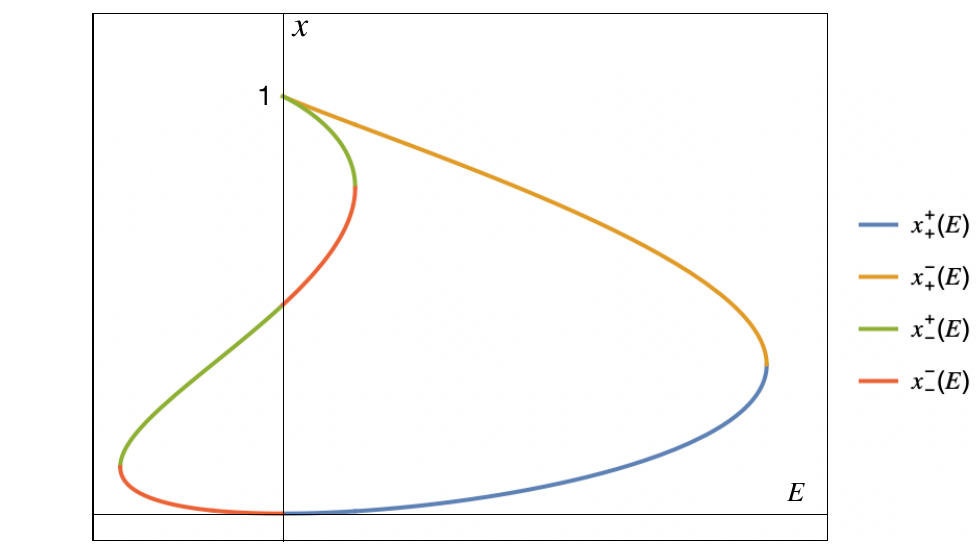}
			\caption{The mutual position
			of the roots
			$x^\pm_\varsigma$, $\varsigma\in\{-1,+1\}$, as in \eqref{radicix}, used in \eqref{ravenna} to define the roots of the quartic polynomial
			for $a_1=1$, $a_2=-2$
			and different value of the energy $E$.}
			\label{xPiuMenoSoluzioni}
	\end{figure}

\noindent
Let us now define 
\begin{equation}\label{radicix}
x^\pm_\varsigma:=
\frac{(t^\pm_\varsigma)^2}{1+
(t^\pm_\varsigma)^2}
=
\frac{|t^\pm_\varsigma|^2}{1+
|t^\pm_\varsigma|^2}
 \,.
\end{equation}
Note that $x^\pm_\varsigma$ 
is an increasing function of 
$|t^\pm_\varsigma|$.
By Lemma \ref{nuovazelanda2}
$x^\pm_\varsigma$ are the roots of 
$\mathbf P(x)$ in \eqref{penisola}.
We now want to order the real roots
$x^\pm_\varsigma$ in increasing order
$x_1<x_2<\ldots$.
We have different cases (see Figure 
\ref{xPiuMenoSoluzioni}):
 \begin{eqnarray}
x_1
&:=&
\left\{
			\begin{array}{ll}
	\min\{x_+^+,x_-^-\}
	& {\rm if}\ \ \delta_\pm>0\,,	\\
	\min\{x_+^+,x_+^-\}
		& {\rm if}\ \ 	\delta_-\leq 0
		\leq \delta_+
			\,,
		\\
		\min\{x_-^+,x_-^-\}
		& {\rm if}\ \ \delta_+\leq 0
		\leq \delta_-		\,,
			\end{array}
			\right.
			\nonumber
\\
x_2
&:=&
\left\{
			\begin{array}{ll}
	\max\{x_+^+,x_-^-\}
	& {\rm if}\ \ \delta_\pm>0\,,	\\
	\max\{x_+^+,x_+^-\}
		& {\rm if}\ \ 	\delta_-\leq 0
		\leq \delta_+
			\,,
		\\
		\max\{x_-^+,x_-^-\}
		& {\rm if}\ \ \delta_+\leq 0
		\leq \delta_-		\,,
			\end{array}
			\right.
			\nonumber
\\
x_3
&:=& \ \quad\min\{x_-^+,x_+^-\}\quad
{\rm if} \quad \delta_\pm>0\,,
\nonumber
\\
x_4
&:=& \ \quad\max\{x_-^+,x_+^-\}\quad
{\rm if} \quad \delta_\pm>0\,.
\label{ravenna}
\end{eqnarray}

 \subsection{The separatrices at the saddle points}\label{sec:disney}

Recall the definition of the zones 
$Z_{ij}$ given in \eqref{kiwi}.
We now consider the curve with zero energy 
bifurcating from the point $(\pi/2,0)$ (recall \eqref{biforco})
in the ``half phase space''
$[0,\pi]\times (0,1)$.
In the case $Z_{10}$
such curve ``turns left''
and touches the line $\{\psi=0\}$ at some point
$>x_1^{(0)}$.
Analogously 
in the case $Z_{01}$
such curve ``turns right''
and touches the line $\{\psi=\pi\}$ at some point
$>x_1^{(\pi)}$.
\\
The situation in the cases $Z_{21}$ and $Z_{12}$ is more 
involved; more precisely
 it depends on the sign of $\Esad$. 
 In particular for $(i,j)\in\left\{(2,1),(1,2)\right\}$ we set
 	 \begin{equation}\label{parigi}
		Z_{ij}^\pm:=\{(a_1,a_2)\in Z_{ij}\ :\ 
		\pm\Esad>0
		\}\,,\qquad
		Z_{ij}^0:=\{(a_1,a_2)\in Z_{ij}\ :\ 
		\Esad=0
		\}\,,
	\end{equation}
so that
$$
Z_{ij}=Z_{ij}^+\cup Z_{ij}^-\cup Z_{ij}^0\,.
$$

The next result characterises      the sets in \eqref{parigi}

\begin{lem}\label{Lemma6}
 Setting
 \begin{equation}\label{gtilde}
\tilde g (a_1):=-\frac2{27} a_1(4 a_1^2+27)
\end{equation}
we have
\begin{eqnarray}
&&		
		Z_{21}^+
=Z_{21}\cap 
\{		a_2>\tilde g(a_1)
		\}\,,\quad
		\quad	Z_{12}^+
=Z_{12}\cap 
\{		a_2>\tilde g(a_1)
		\}\,,
		\nonumber
		\\
		&&
		Z_{21}^-
=Z_{21}\cap 
\{		a_2<\tilde g(a_1)
		\}\,,
		\qquad
		Z_{12}^-
=Z_{12}\cap 
\{		a_2<\tilde g(a_1)
		\}\,,
		\nonumber
		\bigskip
		\\
&&Z_{21}^0
=Z_{21}\cap 
\{		a_2=\tilde g(a_1)
		\}=
		\{a_2=\tilde g(a_1)		,\ \ a_1<0\}
		\,,
		\nonumber
		\\
&&		Z_{12}^0
=Z_{12}\cap 
\{		a_2=\tilde g(a_1)
		\}
=\{a_2=\tilde g(a_1),\ \ a_1>0
		\}		
		\,.	
		\label{zeppole}		
\end{eqnarray}
\end{lem}

\noindent
Note that, since $\tilde g$ is odd and $\tilde g(a_1)\leq -a_1$
for $a_1\geq 0$, by the definition of
$Z_{21}$ and $Z_{12}$ it follows that
$Z_{21}^-\subset \{a_1<0\}$ and 
$Z_{12}^+\subset \{a_1>0\}$.

\noindent
\proof
We discuss only the case 
$Z_{21}$, the study of $Z_{12}$ being analogous.
As we said above, the picture of the phase space in the case
$Z_{21}$ strongly depends on the sign of the energy
of the saddle point $\Esad=F(0,x_2^{(0)})$,
where $x_2^{(0)}$ is the minimum of the function
$x\to F(0,x)$. In particular 
we claim that $\Esad\lesseqqgtr 0$ if and only if
$a_2\lesseqqgtr \tilde g(a_1)$.
In particular we note that 
 $x=x_2^{(0)}$ satisfies the system  	\begin{equation*}
		\left
		\{
		\begin{array}{ll}
		F(0,x)=\frac 12 a_2 x^2+a_1 x+b(x)=0\,,\\
		\partial_x F(0,x)=a_2 x +a_1 +b'(x)=0\,.
		\end{array}
		\right.
	\end{equation*}
By algebraic manipulation we get
	\begin{equation*}
		\left
		\{
		\begin{array}{ll}
		\frac 12 a_2 x^2+ x  b'(x)-b(x)=0\,,\\
		a_1 x +2b(x)-x  b'(x)=0\,,
		\end{array}
		\right.
	\end{equation*}
by which we finally have 
 	\begin{equation*}
		a_2=2\frac{b(x)-x b'(x)}{x^2}\,, \qquad \quad a_1=\frac{x b'(x) -2 b(x)}{x}\,.
	\end{equation*}
By using \eqref{petrolio}
 	\begin{equation}\label{oliva}
		a_1=-\frac 32 \sqrt{\frac{1}{x}-1}\,,
		\qquad
		a_2=
		\sqrt{\frac 1x -1}\left( \frac 1x +2\right)\,.
	\end{equation}
 Note that $a_1<0$. By inverting 
 the first expression in \eqref{oliva}
 	we get
		$	\frac 1x=\frac 49 a_1^2 +1$;
 substituting in the second expression we obtain
 that $a_2=\tilde g(a_1)$ defined in \eqref{gtilde}.
 Therefore in $Z_{21}^0$, namely when
 $a_2=\tilde g(a_1)$, 
 the value of the function
$x\to F(0,x)=\frac12 a_2 x^2+a_1 x+b(x)$
at its minimum  $x_2^{(0)}$   is exactly 0.
On the other hand in
$Z_{21}^+$, namely when
 $a_2>\tilde g(a_1)$, 
one has
$F(0,x_2^{(0)})>0$.
Finally in $Z_{21}^-$ it is 
$F(0,x_2^{(0)})<0$.
	 \eproof

  \subsection{Different topologies of the level curves}\label{6zones}

 Let us consider the energy level sets in the  phase spase
 $$
 \PP:=\mathbb T \times (0,1)\,,
 $$
  which is a cylinder.
 The points where the level curves $\{F=E\}$
 touch  the lines $\psi=0$ or $\psi=\pi$ are the solutions
 of the equation $F(0,x)=a(x)+b(x)=E$ and 
 $F(\pi,x)=a(x)-b(x)=E$, respectively; equivalently
 they are the roots of the quartic
 polynomial in \eqref{isola}.

 \medskip

We note that in the cases $Z_{10}$, $Z_{01}$
the set $\{F=E\}$ has only one connected component.
The same holds in the case  $Z_{21}$ except for
 $\Esad<E<\min\{a(1), \Emax\}$
when $\{F=E\}$ possesses two connected components.
Analogously in the case $Z_{12}$ the level set
$\{F=E\}$ possesses two connected component for
 $\max\{a(1), \Emin\}<E<\Esad$ and only one otherwise.

\begin{rem}\label{badedas}
 Up to the energy level corresponding to
$E=0$ and to the critical energies\footnote{Namely the energy
of critical points of $F$.
In the case
$a_1+a_2=0$, that we are actually excluding (recall Remark \ref{maniavanti}),
 there is also a curve which 
touches the line $x=1$.}, the level sets are curves 
 of three types: 
\\
{\rm (i)} a homotopically trivial, namely contractible, 
curve making a loop around the maximum
$(0,x_1^{(0)})$
intersecting twice the line $\psi=0$;
\\ 
{\rm (ii)} a curve wrapping on the cylinder; in particular it intersects once the line $\psi=0$ and once the line $\psi=\pi$;
\\ 
{\rm (iii)} a  homotopically trivial
curve making a loop around the minimum
$(\pi,x_1^{(\pi)})$
intersecting twice the line $\psi=\pi$.
\end{rem}

In the following we will always 
 label the 
roots of  the quartic
 polynomial in \eqref{isola}
 so that $x_i(E)<x_{i+1}(E)$.
 Recall Proposition \ref{pandivia}.

\medskip
{\bf Case $Z_{10}$.}
 \begin{figure}[h!]
 \center
			\includegraphics[width=10cm]
			{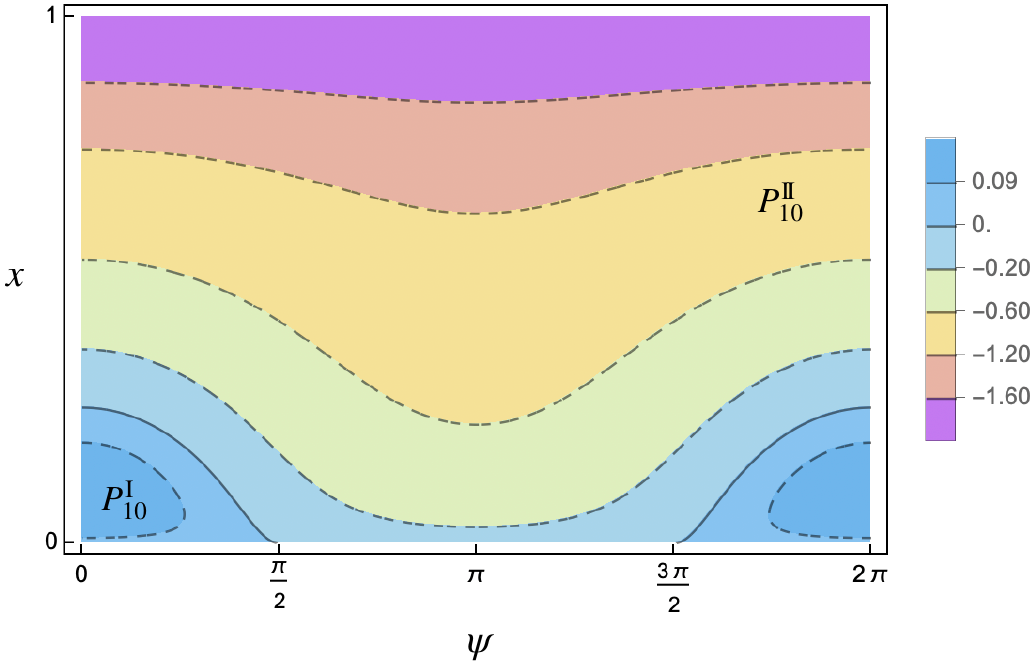}
			\caption{Phase portrait of $Z_{10}$: $(\psi,x)\in [0,2\pi] \times (0,1)$. The zone $Z_{10}$
			for $a_1=-1$ and $a_2=-2$ filled by the level curves of $-x-x^2+b(x)\cos \psi=E$, for different values
			of the energy $E$.}
			\label{Z10_final.pdf}
	\end{figure}
The zero level separatrix 
actually separates the phase space $\PP$
into two open connected components 
$\PP_{10}^{\rmuno}$ and $\PP_{10}^{\rmdue}$
supporting two different kind of motions\footnote{Where
$\{F>0\}:=\{(\psi,x)\in\PP \ :\ F(\psi,x)>0\}$.}
\begin{equation}\label{P10}
\PP_{10}^{\rmuno}:=\{F>0\}\,,\qquad
\PP_{10}^{\rmdue}:=\{F<0\}\,,
\end{equation}
with $\PP=\PP_{10}^{\rmuno}\cup\PP_{10}^{\rmdue}\cup\{F=0\}$.
Indeed in  $\PP_{10}^{\rmuno}$ the level curves have the form in 
case (ii) above, while in $\PP_{10}^{\rmdue}$ they have the form in 
case (i). In the present case the quartic polynomial
in \eqref{isola} possesses, for $E\neq 0$ and not critical,
only two real roots
 $x_1(E)<x_2(E)$.
Note that
$x_1(E)=x_1(E;\Ju)$ and $x_2(E)=x_2(E;\Ju)$.
If $E>0$ 
the $E$--level curve starts at
$(0,x_1(E))$ and come back on the line $\psi=0$
at $(0,x_2(E))$, otherwise, for $E<0$, it joints the line 
$\psi=\pi$ at $(\pi,x_1(E))$ and the line 
$\psi=0$ at $(\pi,x_2(E))$.
Recalling \eqref{psi(x)}, the level curve  $\{F=E\}$
can be expressed as a graph over $x_1(E)<x
<x_2(E)$
by the function $\psi(x;\Ju)$.

\medskip

{\bf Case $Z_{01}$.}	
	\begin{figure}[h]
	\center
			\includegraphics[width=10cm,keepaspectratio]{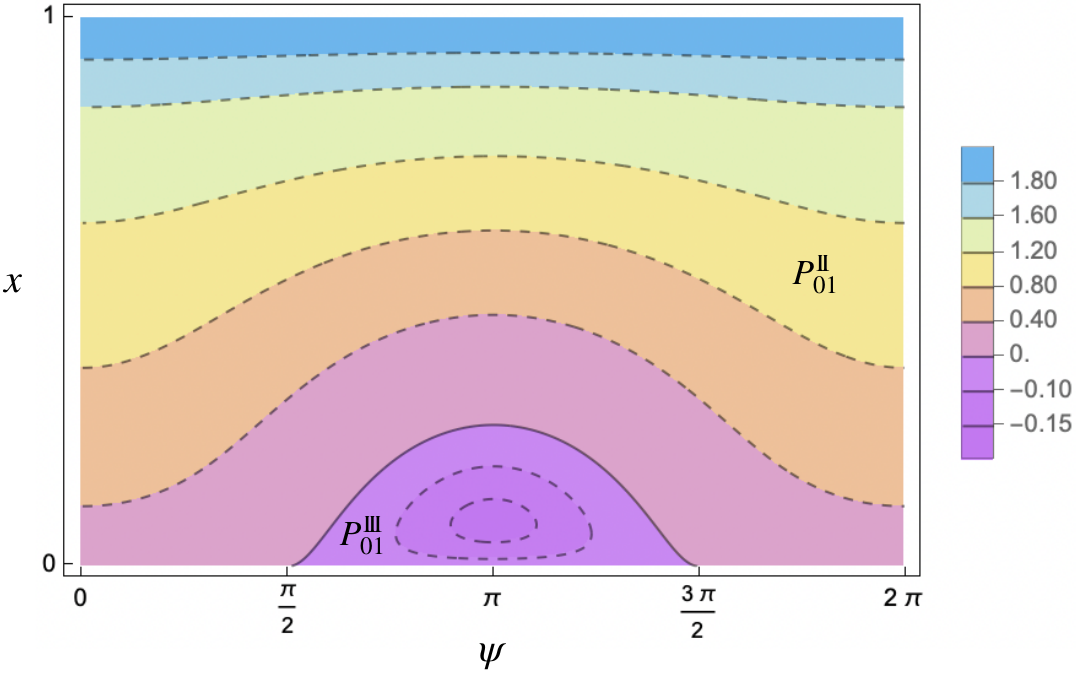}
			\caption{Phase portrait of $Z_{01}$: $(\psi,x)\in [0,2\pi] \times (0,1)$. The zone $Z_{01}$
			for $a_1=1$ and $a_2=2$ filled by the level curves of $x+x^2+b(x)\cos \psi=E$, for different values
			of the energy $E$.}
	\end{figure}
	We set
	\begin{equation}\label{P01}
	\PP_{01}^{\rmdue}:=\{F>0\}
\,,\qquad
\PP_{01}^{\rmtre}:=\{F<0\}
\,,
\end{equation}
with $\PP=\PP_{01}^{\rmtre}\cup\PP_{01}^{\rmdue}\cup\{F=0\}$.
Again
the zero level separatrix 
actually separates the two different kind of motions:
 in  $\PP_{01}^{\rmtre}$ the level curves have the form in 
case (iii) above, while in $\PP_{01}^{\rmdue}$ they have the form in 
case (ii).

\medskip
{\bf Case $Z_{21}^+$.}
\begin{figure}[h]
\center
			\includegraphics[width=10cm,keepaspectratio]
			{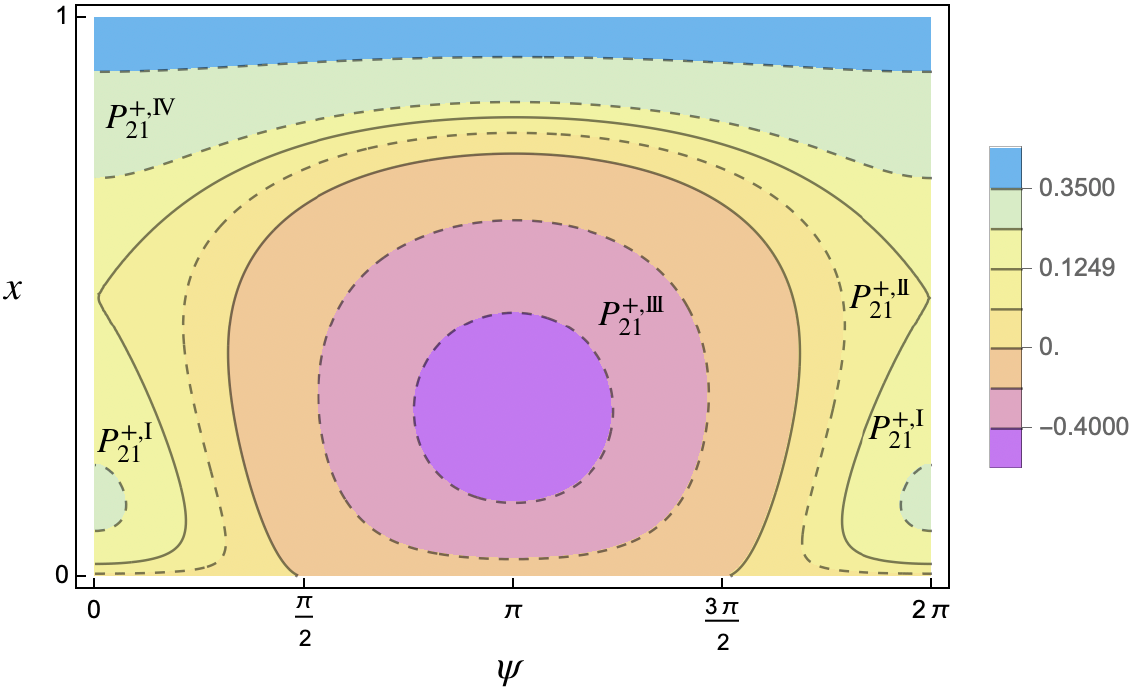}
			\caption{Phase portrait of $Z_{21}^+$: $(\psi,x)\in [0,2\pi] \times (0,1)$. The zone $Z_{21}^+$
			for $a_1=-1$ and $a_2=3$ filled by the level curves of $-x+\frac 32 x^2+b(x)\cos \psi=E$, for different values
			of the energy $E$.}
	\end{figure} 
The zero level separatrix  and the two separatrices 
emanating from the saddle point
$(0, x_2^{(0)})$ with energy\footnote{Recall
Proposition \ref{pandivia} and \eqref{maxmin}.
}
 $F(0,x_2^{(0)})=\Esad>0$ (recall \eqref{parigi})
 separate
 the phase space $\PP$
into 4 open connected components:
 \begin{eqnarray}\label{P21+}
&&\!\PP_{21}^{+,\rmuno}:=\{F>\Esad\ \mbox{containing}\  
(0,x_1^{(0)})\}\,,
\qquad 
\PP_{21}^{+,\rmdue}:=\{0<F<\Esad\}\,,
\\
&&\!\PP_{21}^{+,\rmtre}:=\{F<0\}\,, 
\quad
\PP_{21}^{+,\rmquattro}:=\{F>\Esad\ \mbox{not containing}\  
(x_1^{(0)},0)\}
\nonumber
\end{eqnarray}
with $\PP=\PP_{21}^{+,\rmuno}\cup\PP_{21}^{+,\rmdue}
\cup\PP_{21}^{+,\rmtre}\cup\PP_{21}^{+,\rmquattro}\cup\{F=0\}\cup\{F=\Esad\}$.
The level curves  in $\PP_{21}^{+,\rmdue}$ and $\PP_{21}^{+,\rmquattro}$ have the form
case (ii) above, the ones in $\PP_{21}^{+,\rmuno}$  have the form in 
case (i), finally the ones in $\PP_{21}^{+,\rmtre}$ are as in (iii).
 In particular the level curves  in $\PP_{21}^{+,\rmuno}$
 pass through the points
 $(0,x_1(E))$ and $(0,x_2(E))$; the ones in 
 $\PP_{21}^{+,\rmdue}$ through 
 $(0,x_1(E))$ and $(\pi,x_2(E))$;
 the ones in 
 $\PP_{21}^{+,\rmtre}$ through 
 $(\pi,x_1(E))$ and $(\pi,x_2(E))$;
 the ones in 
 $\PP_{21}^{+,\rmquattro}$ through 
 $(0,x_3(E))$ and $(\pi,x_4(E))$.

\medskip
{\bf Case $Z_{21}^-$.}
\begin{figure}[h]
		\center
		\includegraphics[width=10cm]{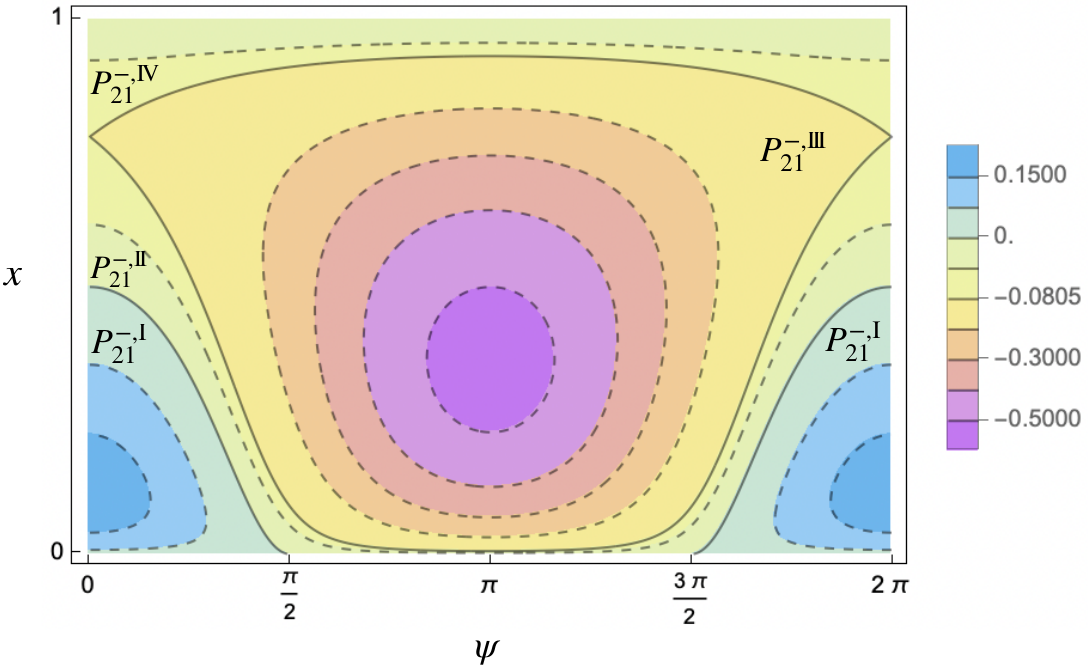}
		\caption{Phase portrait of $Z_{21}^-$: $(\psi,x)\in [0,2\pi] \times (0,1)$. The zone $Z_{21}^-$
			for $a_1=-1$ and $a_2=2$ filled by the level curves of $-x+x^2+b(x)\cos \psi=E$, for different values
			of the energy $E$. }
	\end{figure}
	The zero level separatrix  and the two separatrices 
emanating from the saddle point
$(0,x_2^{(0)})$ with energy
 $F(0,x_2^{(0)})=\Esad<0$ (recall \eqref{parigi})
 separate
 the phase space $\PP$
into 4 open connected components:
 \begin{equation}\label{P21-}
\PP_{21}^{-,\rmuno}:=\{F>0\}\,,\qquad 
\PP_{21}^{-,\rmtre}:=\{F<\Esad
\}\,,
\end{equation}
while $\PP_{21}^{-,\rmdue}$ and $\PP_{21}^{-,\rmquattro}$
are the two open connected components of
$\{\Esad<F<0\}$ with $\PP_{21}^{-,\rmdue}$
containing $(0,\pi)$ in its closure.
We immediately see that 
$$
\PP=\PP_{21}^{-,\rmuno}\cup\PP_{21}^{-,\rmdue}
\cup\PP_{21}^{-,\rmtre}\cup\PP_{21}^{-,\rmquattro}\cup\{F=0\}\cup\{F=\Esad\}\,.
$$
The level curves  in $\PP_{21}^{-,\rmdue}$ and $\PP_{21}^{-,\rmquattro}$ have the form in
case (ii) above, the ones in $\PP_{21}^{-,\rmuno}$  have the form in 
case (i), finally the ones in $\PP_{21}^{-,\rmtre}$ are as in (iii).
 In particular the level curves  in $\PP_{21}^{-,\rmuno}$
 pass through the points
 $(0,x_1(E))$ and $(0,x_2(E))$; the ones in 
 $\PP_{21}^{-,\rmdue}$ through 
 $(\pi,x_1(E))$ and $(0,x_2(E))$;
 the ones in 
 $\PP_{21}^{-,\rmtre}$ through 
 $(\pi,x_1(E))$ and $(\pi,x_2(E))$;
 the ones in 
 $\PP_{21}^{-,\rmquattro}$ through 
 $(0,x_3(E))$ and $(\pi,x_4(E))$.

\medskip
{\bf Case $Z_{12}^+$.}	
	\begin{figure}[h]
	\center
			\includegraphics[width=10cm,keepaspectratio]{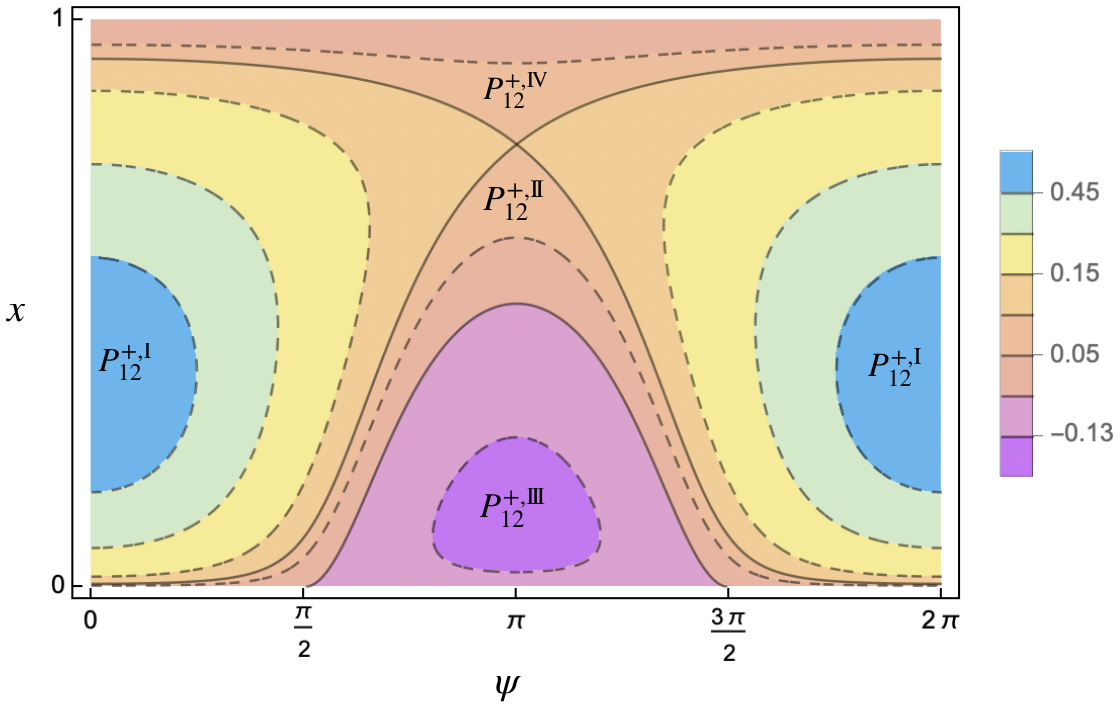}
			\caption{Phase portrait of $Z_{12}^+$: $(\psi,x)\in [0,2\pi] \times (0,1)$. The zone $Z_{12}^+$ for $a_1=1$ and $a_2=-2$ filled by the level curves of $x-x^2+b(x)\cos \psi=E$, for different values
			of the energy $E$.}
	\end{figure}
	The zero level separatrix  and the two separatrices 
emanating from the saddle point
$(\pi,x_2^{(\pi)})$ with energy
 $F(\pi,x_2^{(\pi)})=\Esad>0$ (recall \eqref{parigi})
 separate
 the phase space $\PP$
into 4 open connected components:
\begin{equation}\label{P12+}
\PP_{12}^{+,\rmuno}:=\{F>\Esad\}\,,\qquad 
\PP_{12}^{+,\rmtre}:=\{F<0
\}\,,
\end{equation}
while $\PP_{12}^{+,\rmdue}$ and $\PP_{12}^{+,\rmquattro}$
are the two open connected components of
$\{0<F<\Esad\}$ with $\PP_{12}^{+,\rmdue}$
containing $(0,0)$ in its closure.
We note that 
$$
\PP=\PP_{12}^{+,\rmuno}\cup\PP_{12}^{+,\rmdue}
\cup\PP_{12}^{+,\rmtre}\cup\PP_{12}^{+,\rmquattro}\cup\{F=0\}\cup\{F=\Esad\}\,.
$$
The level curves  in $\PP_{12}^{+,\rmdue}$ and $\PP_{12}^{+,\rmquattro}$ have the form
case (ii) above, the ones in $\PP_{12}^{+,\rmuno}$  have the form in 
case (i), finally the ones in $\PP_{12}^{+,\rmtre}$ are as in (iii).
 In particular the level curves  in $\PP_{12}^{+,\rmuno}$
 pass through the points
 $(0,x_1(E))$ and $(0,x_2(E))$; the ones in 
 $\PP_{12}^{+,\rmdue}$ through 
 $(0,x_1(E))$ and $(\pi,x_2(E))$;
 the ones in 
 $\PP_{12}^{+,\rmtre}$ through 
 $(\pi,x_1(E))$ and $(\pi,x_2(E))$;
 the ones in 
 $\PP_{12}^{+,\rmquattro}$ through 
 $(\pi,x_3(E))$ and $(0,x_4(E))$.
	
	\medskip
{\bf Case $Z_{12}^-$.}
	 \begin{figure}[h]
		\center
		\includegraphics[width=10cm,keepaspectratio]{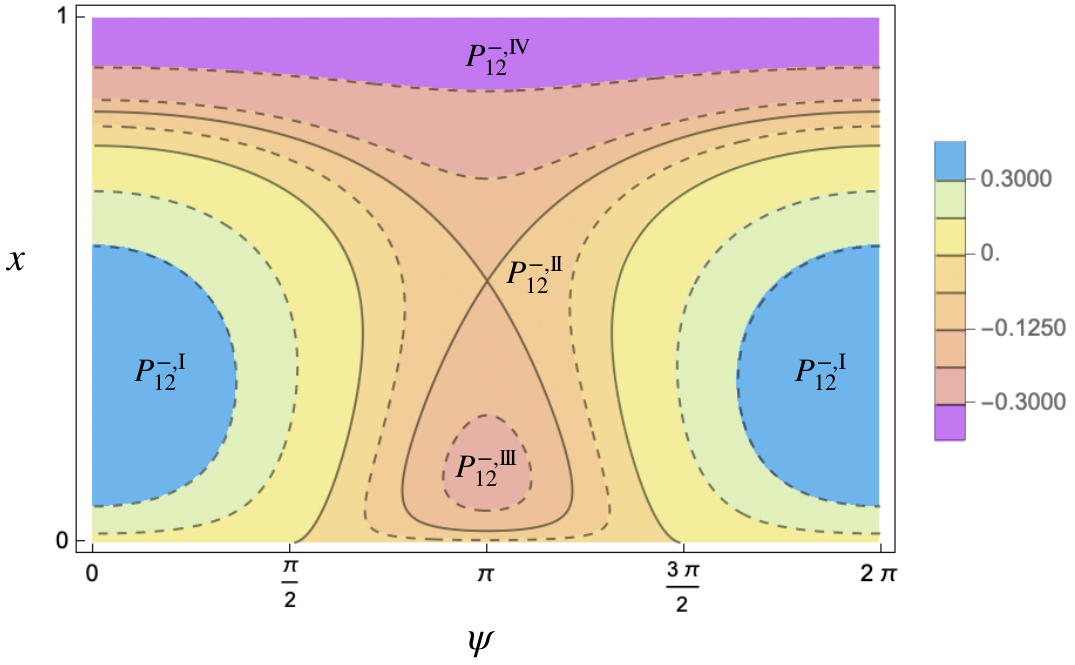}
		\caption{Phase portrait of $Z_{12}^-$: $(\psi,x)\in [0,2\pi] \times (0,1)$. The zone $Z_{12}^-$ for $a_1=1$ and $a_2=-3$ filled by the level curves of $x-\frac 32 x^2+b(x)\cos \psi=E$, for different values
			of the energy $E$.}
\label{Z12Meno_Mathematica}
	\end{figure}
The zero level separatrix  and the two separatrices 
emanating from the saddle point
$(\pi,x_2^{(\pi)})$ with energy
 $F(\pi,x_2^{(\pi)})=\Esad<0$ (recall \eqref{parigi})
 separate
 the phase space $\PP$
into 4 open connected components:
 \begin{eqnarray}\label{P12-}
&&\PP_{12}^{-,\rmuno}:=\{F>0\}\,,\qquad 
\PP_{12}^{-,\rmtre}:=\{F<\Esad\ \mbox{containing}\  
(x_1^{(\pi)},\pi)\}\,,
\nonumber
\\
&& 
\PP_{12}^{-,\rmdue}:=\{\Esad<F<0\}\,,\quad
\PP_{12}^{-,\rmquattro}:=\{F<\Esad\ \mbox{not containing}\  
(x_1^{(\pi)},\pi)\}\,,
\end{eqnarray}
with $\PP=\PP_{12}^{-,\rmuno}\cup\PP_{12}^{-,\rmdue}
\cup\PP_{12}^{-,\rmtre}\cup\PP_{12}^{-,\rmquattro}\cup\{F=0\}\cup\{F=\Esad\}$.
The level curves  in $\PP_{12}^{-,\rmdue}$ and $\PP_{12}^{-,\rmquattro}$ have the form in
case (ii) above, the ones in $\PP_{12}^{-,\rmuno}$  have the form in 
case (i), finally the ones in $\PP_{12}^{-,\rmtre}$ are as in (iii).
 In particular the level curves  in $\PP_{12}^{-,\rmuno}$
 pass through the points
 $(0,x_1(E))$ and $(0,x_2(E))$; the ones in 
 $\PP_{12}^{-,\rmdue}$ through 
 $(\pi,x_1(E))$ and $(0,x_2(E))$;
 the ones in 
 $\PP_{12}^{-,\rmtre}$ through 
 $(\pi,x_1(E))$ and $(\pi,x_2(E))$;
 the ones in 
 $\PP_{12}^{-,\rmquattro}$ through 
 $(\pi,x_3(E))$ and $(0,x_4(E))$.

\subsection{Degenerate cases}\label{degenere}

\begin{figure}[h!]
		\center
		\includegraphics[width=17.2cm,keepaspectratio]{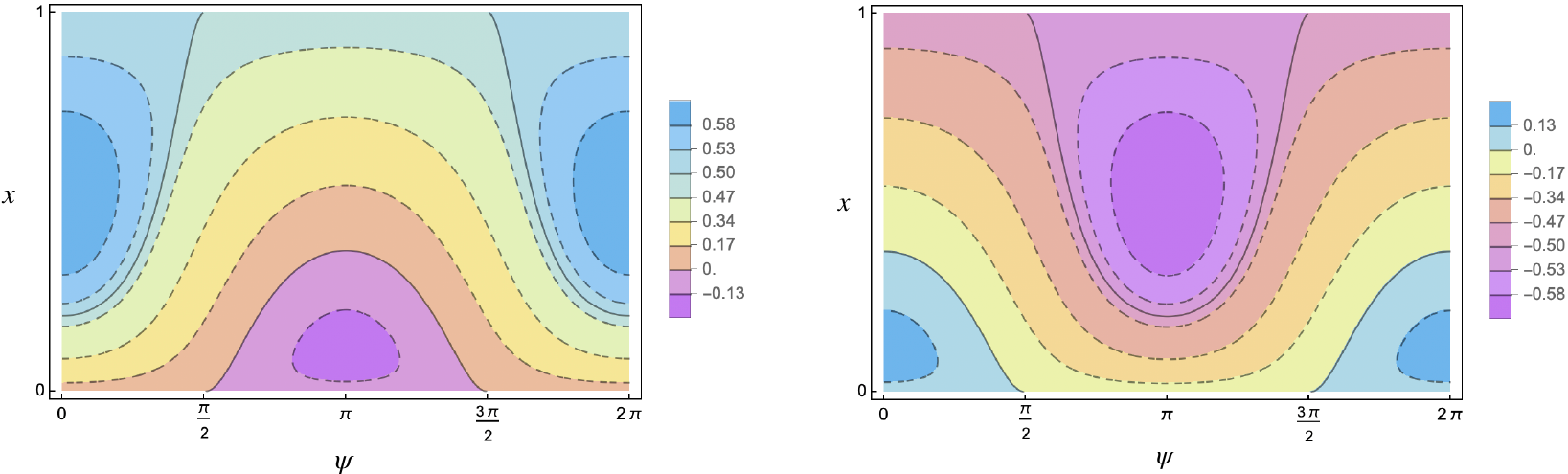}
					\caption{The degenerate case $a_1=-a_2$.	
			Phase portraits: $(\psi,x)\in [0,2\pi] \times (0,1)$. 
On the left, the case
with $a_1=1$, $a_2=-1$,
 filled by the level curves $x-1/2 x^2+b(x)\cos \psi=E$.
 On the right, the case
 with $a_1=-1$, $a_2=1$,
  filled by the level curves $-x+1/2 x^2
			+b(x)\cos \psi=E$.
			In both cases a new separatrix appears 
			approaching the line $x=1$ at $\psi=\pi/2, 3\pi/2$.
			}
			\label{deg12}
	\end{figure}

\begin{figure}[h!]
		\center
		\includegraphics[width=17.2cm,keepaspectratio]{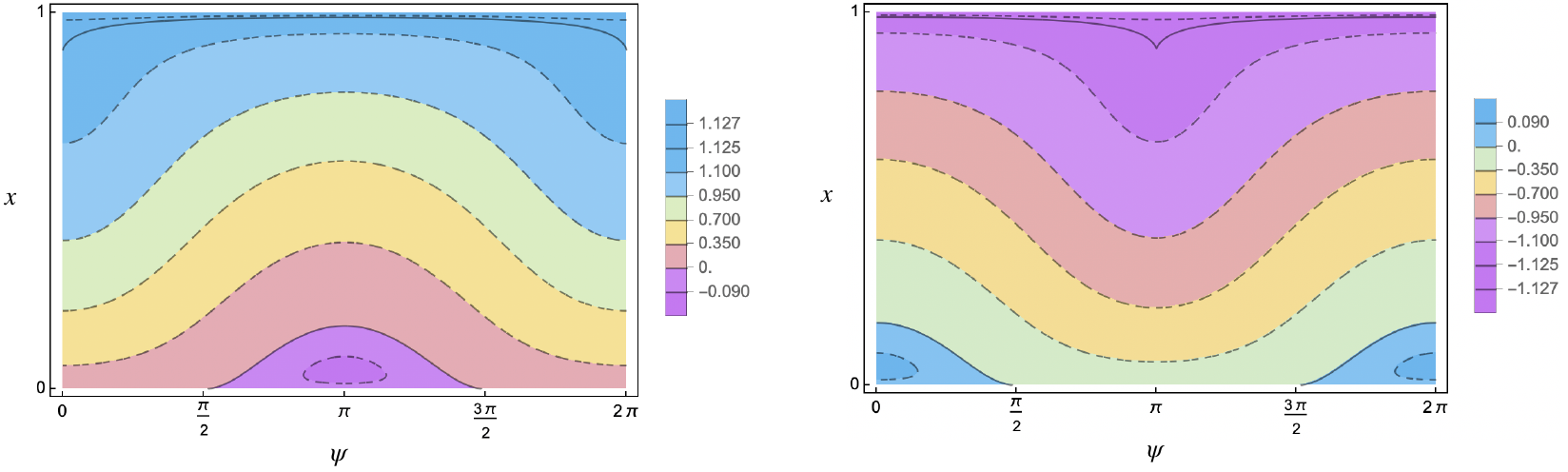}
			\caption{The degenerate cases $a_2=g(a_1)$
			on the left and $a_2=-g(-a_1)$.	
			Phase portraits: $(\psi,x)\in [0,2\pi] \times (0,1)$. 
On the left,  $a_1=2$, $a_2=\tilde g(2)=-\frac{47}{27}$,
 filled by the level curves $2x-\frac {47}{54} x^2+b(x)\cos \psi=E$.
On the right,  $a_1=-2$, $a_2=-\tilde g(2)=\frac{47}{27}$,
 filled by the level curves $-2x+\frac {47}{54} x^2+b(x)\cos \psi=E$.
Note that a non-smooth curve appears with a cusp.}
			\label{deg34}
	\end{figure}

\begin{figure}[h!]
		\center
	\includegraphics[width=17.2cm,keepaspectratio]{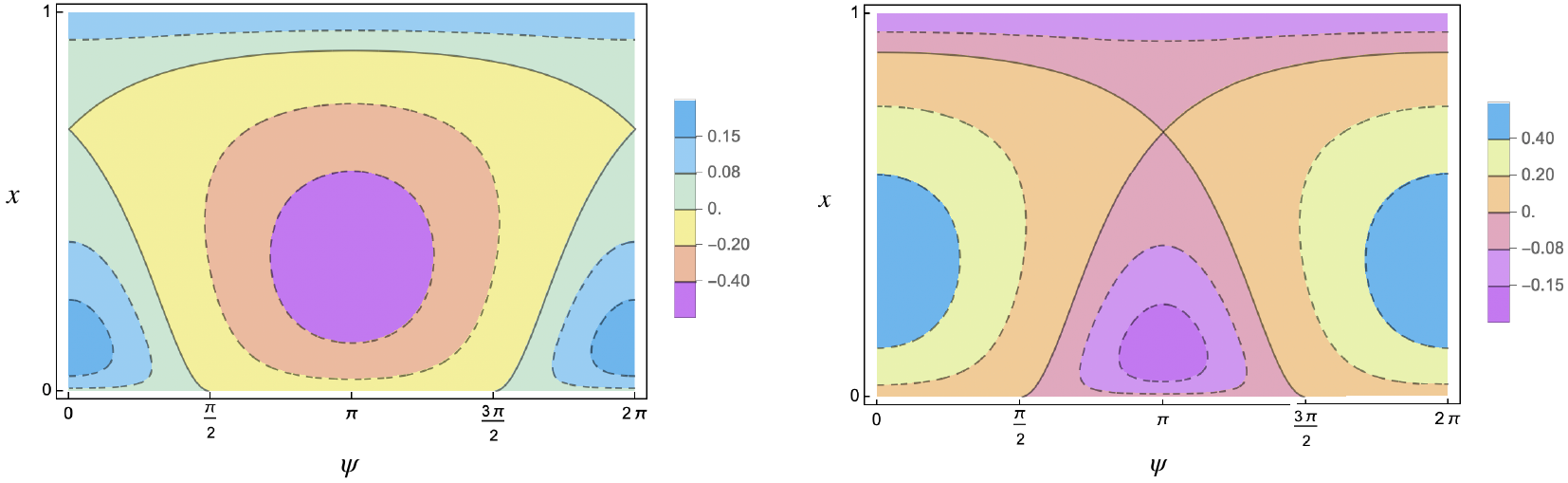}
			\caption{The degenerate case $a_2=\tilde g(a_1)$.	
			Phase portraits: $(\psi,x)\in [0,2\pi] \times (0,1)$. 
On the left, the zone $Z_{21}^0$
with $a_1=-1$, $a_2=\tilde g(-1)=\frac{62}{27}$,
 filled by the level curves $-x+\frac {31}{27} x^2+b(x)\cos \psi=E$.
 On the right, the zone $Z_{12}^0$
 with $a_1=1$, $a_2=\tilde g(1)=-\frac{62}{27}$,
  filled by the level curves $x-\frac {31}{27} x^2
			+b(x)\cos \psi=E$.
			In both cases the stable and unstable
			manifolds of the saddle point
			have zero energy.}
			\label{deg56}	
		\end{figure}

Recalling Remark \ref{maniavanti},
we briefly illustrate in Figures \ref{deg12}-\ref{deg56}  the degenerate cases:
 $a_2=-a_1$,
 when $x=1$ is a solution of
 \eqref{gingerina1}-\eqref{gingerina2},
 $a_2=g(a_1)$ and $a_2=-g(-a_1)$, 
 when two solutions coincide, finally 
 $a_2=\tilde g(a_1)$, when the separatrix and the
 stable and unstable manifolds of the
 saddle point coincide and have zero energy.


 \section{Explicit formulae of  the nonlinear frequencies}\label{sediablu}
 
 In this section we first write
 the integrating action 
 $\II_1$ as a function of the energy
 $E$ in terms of integrals 
 in the $x$ variables with endpoints
 given by  the roots of the quartic 
 polynomial 
 $\mathbf P(x)$ in \eqref{penisola},
 studied in the previous section.
 In addition to energy, these representation 
 formulae  depend on the
 values of the parameters $a_1$
 and $a_2$, according to the
 resulting different
 topologies of the phase space
 described above.
 
 The final integrated Hamiltonian
 is  the inverse
 $\cE:\II_1\to \cE(\II_1;\II_2)$
  of the function
 $E\to \II_1(E;\II_2)$ in
 \eqref{Action(E)}.
 Its derivatives with respect
 to $\II_1$ and $\II_2$ are the nonlinear frequencies
 and can be written in terms of the derivatives of $ \II_1(E;\II_2)$
 with respect to $E$ and $\II_2$,
 see  \eqref{501bis}.
 These derivatives are 
 expressed  in terms of elliptic
 integrals in Proposition
 \ref{WWW}.
 The integrals are explicitly
 evaluated
  by means of suitable Moebius transformations in Subsections
 \ref{Sec:elliptic} and \ref{somorta},
 in the case that $\mathbf P(x)$
has four  or two real roots, respectively.
In the last subsection we
consider  the exact 3:1 resonance 
case, where the above formulae
simplify a bit.

	\subsection{Construction of the integrating action variables}\label{scrauso}

Since $\hat H$ has two independent integrals 
of motions: the Hamiltonian itself and $\Ju$,
 by the Arnold-Liouville theorem the Hamiltonian 
$\hat H$ is integrable.
A part from $\II_2:=\Ju$ the construction of the 
other action $\II_1$ as function of the energy $E$ is as follows.
$\II_1(E)$ is simply  the area
 enclosed by the level curves of $\hat H=E$ 
 divided by $2\pi$. Such level curves coincide with the ones of
 $F$.

\noindent
Our aim  is to find
a symplectic map 
		$
		\Psi:(\II,\theta)=(\II_1,\II_2,\theta_1,\theta_2)\to(J_1,\Ju,\psi_1,\psi_2)
		$,
fixing 
\begin{equation}\label{Idue}
\II_2=\Ju\,,
\end{equation}
such that, in the new coordinates,
the Hamiltonian $\hat H$ is integrated, 
namely\footnote{Recalling \eqref{pluto}
note that $\cE$ is adimensional.}
	\begin{equation}\label{bicicletta}
				\hat H\circ \Psi=:\cE(\II)
	\end{equation}
depends only on the new actions 
		$
			\II=(\II_1,\II_2)
		$. 
		
	\noindent
	Note that the same transformation 
		$\Psi$ also integrates 
		$\mathbb H_{\Ju}=\mathbb H_{\II_2}$ in \eqref{pluto} and 
		$\hat{\mathbb H}_{\rm res}$ in \eqref{secularRES}.
		Indeed
		\begin{equation}\label{bicicletta2}
		\mathbb H_{\II_2}\circ \Psi=\chi \,
				\II_2^2\, \cE(\II)
				\qquad
				\mbox{and}
				\qquad
				\hat{\mathbb H}_{\rm res}\circ \Psi=
				\mathbb E(\II):=
				\omega_- \II_2+
				\chi \,
				\II_2^2\, \cE(\II)\,.
	\end{equation}
In the new coordinates, the actions are constants of motion
and the angles perform a linear motion 
		$
		\theta(t)=\theta(0)+\omega t
		$ 
with frequencies
		\begin{equation}\label{501}
				\omega=(\omega_1,\omega_2):=
				(\partial_{\II_1}\mathbb E,
				\partial_{\II_2}\mathbb E)
				=(\chi \,
				\II_2^2\, \partial_{\II_1}\cE\,,\  \ 
				\omega_-+2\chi\,\II_2  \cE+\chi\,
				\II_2^2\, \partial_{\II_2}\cE)\,.
		\end{equation}	
The classical construction of the 	Hamiltonian $\cE$, ``the adimensional energy'',
is as follows. First one constructs, for every fixed value of $\II_2=\Ju$,
 the action function  	
$\II_1:E\to \II_1(E;\Ju)$ defined as the area enclosed by the level
curve $\gamma_E:=\{\hat H_\Ju=E+a_0\}$ normalised by $2\pi$.
Then, since the function  	
$\II_1:E\to \II_1(E;\Ju)$ turns out to be monotone
(being $|\partial_E\II_1(E;\Ju)|>0$), one defines $\cE:\II_1\to \cE(\II_1,\Ju)$
as its inverse.
Namely, in view of \eqref{loacker},
\begin{equation}\label{levis}
\cE\big(\II_1(E;\Ju),\Ju\big)=\cE\big(\II_1(E;\II_2),\II_2\big)=E+a_0\,.
\end{equation}
So the level curves of $\hat H_\Ju$ play a crucial role here.	
Note that by \eqref{pluto} the level curves of 
			$\hat H_\Ju$ are the same as the ones of $\mathbb H_{\Ju}$,
			moreover by \eqref{loacker} they are simple related
			to the ones of $F$.
\\
More precisely the new action is defined as\footnote{Recall \eqref{loacker}.} 
	\begin{equation}\label{Action(E)}
		\II_1(E)=\II_1(E;\II_2):=\frac{\II_2}3 \Ac(E;\II_2)\,,
	\end{equation}
where, recalling the notation introduced in
 Remark \ref{badedas}, $\Ac$ is the area (normalised by $2\pi$)
enclosed by the $E$-level curve  in the cases (i) and (iii),
 and below the level curve
in the case (ii).
In particular we have four cases indexed by
 $\rmuno,\rmdue,\rmtre,\rmquattro$,
 according if one is in the zones
$\PP_{ij}^{\rmuno},\PP_{ij}^{\rmdue},
\PP_{ij}^{\rmtre},\PP_{ij}^{\pm,\rmuno},
\PP_{ij}^{\pm,\rmdue},\PP_{ij}^{\pm,\rmtre},
\PP_{ij}^{\pm,\rmquattro}$.
\\
Case $\rmuno$. The level curve makes a loop around the maximum
$(0,x_1^{(0)})$ then\footnote{Recall \eqref{psi(x)}.}
	\begin{equation}\label{AreaI}
		\Ac(E)=\Ac^{\rmuno}(E;\II_2):=\frac{\mathrm{Area}}{\pi}=\frac{1}{\pi}\int_{x_1(E)}^{x_2(E)}\psi(x;E)\,dx\,.
	\end{equation}
	This holds in the zones: 
	$\PP_{10}^{\rmuno},
	\PP_{ij}^{\pm,\rmuno}.$
	\\
Case $\rmdue$. 
The level curve wraps  on the cylinder
	\begin{equation}\label{AreaII}
		\Ac(E)=\Ac^\rmdue(E;\II_2):=\frac{\mathrm{Area}}{\pi}=
		\left\{
			\begin{array}{ll}
	x_1(E)+
		\frac{1}{\pi}\int_{x_1(E)}^{x_2(E)}\psi(x;E)\,dx\,,
	& {\rm if}\ \ \psi(x_1)=\pi	\\
	x_2(E)-
		\frac{1}{\pi}\int_{x_1(E)}^{x_2(E)}\psi(x;E)\,dx\,,
		& {\rm if}\ \ \psi(x_1)=0		\,.
			\end{array}
			\right.
	\end{equation}
 in the cases $\PP_{ij}^{\rmdue},
 \PP_{ij}^{\pm,\rmdue}$.
\\
Case $\rmtre$.
The level curve makes a loop around the minimum
$(\pi,x_1^{(\pi)})$
\begin{equation}\label{AreaIII}
		\Ac(E)=\Ac^\rmtre(E;\II_2):=\frac{\mathrm{Area}}{\pi}=\frac{1}{\pi}\int_{x_1(E)}^{x_2(E)}\big(\pi-\psi(x;E)\big)\,dx\,.
	\end{equation}
	This holds in the zones: $
	\PP_{01}^{\rmtre},\PP_{ij}^{\pm,\rmtre}$.
	\\
Case $\rmquattro$.  The level curve wraps  on the cylinder
\begin{equation}\label{AreaIV}
		\Ac(E)=\Ac^\rmquattro(E;\II_2):=\frac{\mathrm{Area}}{\pi}=
		\left\{
			\begin{array}{ll}
	x_3(E)+
		\frac{1}{\pi}\int_{x_3(E)}^{x_4(E)}\psi(x;E)\,dx\,,
	& {\rm if}\ \ \psi(x_3)=\pi	\\
	x_4(E)-
		\frac{1}{\pi}\int_{x_3(E)}^{x_4(E)}\psi(x;E)\,dx\,,
		& {\rm if}\ \ \psi(x_3)=0		\,.
			\end{array}
			\right.
	\end{equation}
in the cases $
 \PP_{ij}^{\pm,\rmquattro}$.

\begin{rem}[KAM Theory]\label{KAM}
The above integrating construction holds
for the truncated Hamiltonian $\hat{\mathbb H}_{\rm res}$ in \eqref{secularRES}
but it does not work for 
the complete Hamiltonian.
In fact 
 the complete system is genuinely two dimensional
 and, therefore, not integrable. In particular
  $\Ju$ is not more a constant of motion.
  One might wonder whether, for $\epsilon$ small enough, the invariant structures, both NNMs and stable and unstable manifolds, that exist for the truncated Hamiltonian survive, slightly deformed, for the full Hamiltonian. The answer is substantially positive
  thanks to KAM Theory.
  More precisely, the hyperbolic periodic orbit and its (local) stable and 
  unstable manifolds survive
  as can be demonstrated following, e.g., 
  \cite{Graff} and \cite{Val}.
  The conservation of two dimensional invariant tori
  is ensured  when the frequencies are strongly rationally independent. This implies that the majority of invariant tori
  still exist in the complete system, whereas
  a minority is destroyed.
  However we note that, in this resonant case,
  the application of KAM Theory is not straightforward.
  In fact the standard KAM theory only regards
  the persistence of the so called {\sl primary} tori,
  namely tori that are graphs over the angles.
  However, as we have already shown, in the 
  resonant case also the so called {\sl secondary} tori
  appear (the blue and the yellow tori in Figure \ref{dama2}).
  All our analysis can bee seen as a 
  necessary 
  preparatory 
  step in view of the application of 
  KAM techniques, since it integrates the 
  resonant BNF up to order four. This means that,
  in the final  action angle variables,  the invariant tori 
  {\sl are graphs} over the angles
  and KAM methods can be applied.
  For a KAM result in presence of resonances
  and the persistence of secondary tori see
  \cite{MNT}.
  \\
  Finally we note that, since the complete system is, in general, 
  not integrable, KAM tori do not completely fill the phase space
  but some gaps appear between them. In these gaps chaotic
  behaviour may occur. However one has to notice that, since we 
  are in two degrees of freedom, every orbit is perpetually stable
  in the  sense that the solutions exist for all times
  and the values of the action variables remain close to the initial ones
  forever. The argument is standard in KAM Theory: 
  the orbits evolve on the three dimensional energy surface
  we have two cases. 1) If on orbit starts on a KAM torus, 
  then it remains on it forever, since the torus is invariant for the 
  Hamiltonian flow. 2) If an orbit starts in a gap between 
  two KAM tori then, since the tori are invariant and bidimensional
  and the energy surface is three dimensional,
  the orbit cannot cross them and
   it remains trapped between them forever.
 \end{rem}

\subsection{Evaluation of the nonlinear frequencies as functions of the energy}
\label{sec:newfrequencies}

 In evaluating the new frequencies in \eqref{501},  it is convenient to use $(E;\II_2)$ as independent
variables, rather than $(\II_1,\II_2)$.
In particular, we have to evaluate
$\partial_{\II_1}\cE\big(\II_1(E;\II_2),\II_2 \big)$
and $\partial_{\II_2}\cE\big(\II_1(E;\II_2),\II_2 \big)$.
Deriving \eqref{levis} with respect to  $E$ we get
$$
\partial_{\II_1}\cE\big(\II_1(E;\II_2),\II_2 \big)\partial_E \II_1(E;\II_2)=1\,.
$$
Then
	\begin{equation}\label{peperosa}
\partial_{\II_1} \cE\big(\II_1(E;\II_2),\II_2\big)=
		\frac{1}{\partial_E \II_1(E;\II_2)}\,.
\end{equation}
Analogously, deriving \eqref{levis} with respect to  $\II_2$, we get
	$$
\partial_{\II_1}\cE\big(\II_1(E;\II_2),\II_2 \big)\partial_{\II_2} \II_1(E;\II_2)
+\partial_{\II_2}\cE\big(\II_1(E;\II_2),\II_2 \big)=0\,,
$$
and, therefore,
\begin{equation}\label{peperosa2}
		\partial_{\II_2} \cE\big(\II_1(E;\II_2),\II_2\big)=
		-\partial_{\II_1}\cE\big(\II_1(E;\II_2),\II_2 \big)
		\partial_{\II_2} \II_1(E;\II_2)
		\stackrel{\eqref{peperosa}}=
		-\frac{\partial_{\II_2} \II_1(E;\II_2)}{\partial_{E} \II_1(E;\II_2)}\,.
\end{equation}
Then, using \eqref{levis}, we rewrite \eqref{501} as
$$
\omega_1(E,\II_2) 
=\chi \,\II_2^2\, \frac{1}{\partial_E \II_1(E;\II_2)}\,,
\qquad
\omega_2(E,\II_2)
=
				\omega_-+2\chi\,\II_2 
				(E+a_0)
				-\chi\,
				\II_2^2\,
				\frac{\partial_{\II_2} \II_1(E;\II_2)}{\partial_{E} \II_1(E;\II_2)}\,,
$$
namely, recalling
\eqref{Action(E)},
\begin{equation}\label{501bis}
\omega_1(E,\II_2) 
=3\chi \,\II_2\, \frac{1}{\partial_E \Ac(E;\II_2)}\,,
\qquad
\omega_2(E,\II_2)
=
				\omega_-+2\chi\,\II_2  
				(E+a_0)
				-\chi\,
				\II_2^2\,
				\frac{\partial_{\II_2} 
				\Ac(E;\II_2)}{\partial_{E} \Ac(E;\II_2)}\,.
\end{equation}
As a final symplectic change of variables we consider the inverse
of the map in \eqref{salmone}, namely
the map $\tilde \Phi:(\tilde I,\tilde \varphi)\to(\II,\theta)$
\begin{equation}\label{salmoneBIS}
	\left\{
			\begin{array}{l}
			\II_1=\tilde I_2\\
			\II_2=\tilde	I_1+ 3\tilde I_2	
						\,,
			\end{array}
			\right.
			\qquad 
			\left\{
			\begin{array}{l}
			\theta_1=\tilde\varphi_2-3\tilde\varphi_1\\
				\theta_2=\tilde\varphi_1\,.
			\end{array}
			\right.
	\end{equation}
Applying the above map to the Hamiltonian $\mathbb E$ 
in \eqref{bicicletta2} we get 
$
\tilde{\mathbb E}:=\mathbb E\circ \tilde \Phi
$, namely
\begin{equation}\label{lisbona}
\tilde{\mathbb E}(\tilde I)=\mathbb E(\tilde I_2,
			\tilde	I_1+ 3\tilde I_2)\,.
\end{equation}	
In order to describe the frequencies of 
$\tilde{\mathbb E}(\tilde I)$ it is convenient to use
$(E,\II_2)$ as variables
instead of $(\tilde I_1,\tilde I_2)$.
The (invertible) relation between the two set of variable is the following
\begin{equation}\label{marasma}
E=\frac{\mathbb E(\tilde I_2,
			\tilde	I_1+ 3\tilde I_2)-\omega_-
			(\tilde I_1+ 3\tilde I_2)}{\chi\,(\tilde	I_1+ 3\tilde I_2)^2}
-a_0
\,,\qquad
\II_2=\tilde	I_1+ 3\tilde I_2\,,
\end{equation}
(recalling \eqref{bicicletta2}, \eqref{levis}).
We are now able to evaluate
the final nonlinear  frequencies, namely the partial derivatives
of $\tilde{\mathbb E}(\tilde I)$ in \eqref{lisbona}, namely
\begin{equation}\label{marasma2}
\omega_-^{\rm nlr}
:=
\partial_{\tilde I_1}\tilde{\mathbb E}
=\omega_2\,,\qquad
\omega_+^{\rm nlr}
:=
\partial_{\tilde I_2}\tilde{\mathbb E}
=\omega_1+3\omega_2
\,.
\end{equation}
Indeed, recalling \eqref{501} and \eqref{501bis}, we have
\begin{eqnarray}
\omega_-^{\rm nlr}(E,\II_2)
&:=&
\omega_-+
				\chi\,\II_2\left(  
				2(E+a_0)
				-
				\II_2\,
				\frac{\partial_{\II_2} 
				\Ac(E;\II_2)}{\partial_{E} \Ac(E;\II_2)}\right)
\,,
\nonumber
\\
\omega_+^{\rm nlr}(E,\II_2)
&:=&
3\omega_-+ 3\chi\,\II_2\left(  
				2(E+a_0)
				+\frac{1}{\partial_E \Ac(E;\II_2)}
				-
				\II_2\,
				\frac{\partial_{\II_2} 
				\Ac(E;\II_2)}{\partial_{E} \Ac(E;\II_2)}\right)\,.
\label{california}				
\end{eqnarray}

\noindent
It remains 
 to evaluate
$\partial_{E} \Ac(E;\II_2)$ and
$\partial_{\II_2} \Ac(E;\II_2)$.
\begin{pro}\label{WWW}
Set\footnote{$\mathbf P(x)$
was defined in \eqref{penisola}. Note that $\mathbf P(x_i)=0$
for $i=1,2,3,4$.}
\begin{equation}\label{W}
W(x,E,\II_2)
:=\frac{1}{\pi\sqrt{(1-x)^3x-\big(E-\frac12 a_2x^2
-a_1 x\big)^2}}
=\frac{1}{\pi\sqrt{-
\mathbf P(x)}}\,.
\end{equation}
 In the zones labelled by $\rmuno,\rmdue,\rmtre$
 \begin{eqnarray}\label{Der}
		\partial_E\Ac(E;\II_2)
		&=&
		\pm\int_{x_1(E,\II_2)}^{x_2(E,\II_2)}
		W(x,E,\II_2)\,dx\,,
		\nonumber
		\\
		\partial_{\II_2}\Ac(E;\II_2)
		&=&
		\pm
		\frac{\sigma}{3 \chi
		\II_2^2}\int_{x_1(E,\II_2)}^{x_2(E,\II_2)}
		x\,W(x,E,\II_2)\,dx\,,
	\end{eqnarray}
	where the $+$ sign holds in the zones labelled by 
	$\rmtre$ and 
	$\PP_{01}^{\rmdue}, \PP_{21}^{+,\rmdue},  \PP_{12}^{+,\rmdue}$,
	while the $-$ sign
	in the zones labelled by 
	$\rmuno$ and $\PP_{10}^{\rmdue}, \PP_{21}^{-,\rmdue}
	, \PP_{12}^{-,\rmdue}$.
Finally	
  \begin{eqnarray}\label{DerBIS}
		\partial_E\Ac^{\rmquattro}(E;\II_2)
		&=&
		\pm\int_{x_3(E,\II_2)}^{x_4(E,\II_2)}
		W(x,E,\II_2)\,dx\,,
		\nonumber
		\\
		\partial_{\II_2}\Ac^{\rmquattro}(E;\II_2)
		&=&
		\pm
		\frac{\sigma}{3 \chi
		\II_2^2}\int_{x_3(E,\II_2)}^{x_4(E,\II_2)}
		x\,W(x,E,\II_2)\,dx\,.
	\end{eqnarray}
	where the $+$ sign holds in the zones
	$\PP_{21}^{\pm,\rmquattro}$
	and the $-$ one in $\PP_{12}^{\pm,\rmquattro}$.
\end{pro}
\proof
First note that
from \eqref{psi(x)} and \eqref{sublime} we get
	\begin{eqnarray}\label{derivopsi(x)}
		\partial_E \psi(x,\II_2)
		&=&
		-\frac 1{b(x)}\cdot \frac{1}{\sqrt{1-\left(\displaystyle\frac{E-a(x)}{b(x)}\right)^2}}=-\frac{1}{\sqrt{b(x)^2-(E-a(x))^2}}
		=-W(x,E,\II_2)\,,
		\nonumber
		\\
		\partial_{\II_2} \psi(x,\II_2)
		&=&
		\frac{3\omega_- -\omega_+}{3 \chi
		\II_2^2}\,x\, W(x,E,\II_2)\,.
	\end{eqnarray}
Case $\rmuno$. Since 
$\psi(x_2(E,\II_2);E,\II_2)
=\psi(x_1(E,\II_2);E,\II_2)=0$, we 
have\footnote{For brevity
we omit to write the dependence on $\II_2$.}
	\begin{eqnarray*}
		\partial_E \Ac^\rmuno(E)&=&
		\frac 1\pi \left[\psi(x_2(E);E)
		\partial_E 
		x_2(E)-\psi(x_1(E);E)\partial_E x_1(E)\right]+\frac 1\pi\int_{x_1(E)}^{x_2(E)}\partial_E \psi(x;E)\,dx\\
		&=&\frac 1\pi\int_{x_1(E)}^{x_2(E)}\partial_E \psi(x;E)\,dx
		\stackrel{\eqref{derivopsi(x)}}=-\int_{x_1(E)}^{x_2(E)}
		W(x,E,\II_2)\,dx\,,\,,
	\end{eqnarray*}
and, analogously,
\begin{eqnarray*}
		\partial_{\II_2} \Ac^\rmuno(E)&=&
		\frac 1\pi \left[\psi(x_2(E);E)
		\partial_{\II_2} 
		x_2(E)-\psi(x_1(E);E)\partial_{\II_2}
		 x_1(E)\right]+\frac 1\pi\int_{x_1(E)}^{x_2(E)}\partial_{\II_2} \psi(x;E)\,dx\\
		&=&\frac 1\pi\int_{x_1(E)}^{x_2(E)}\partial_{\II_2} \psi(x;E)\,dx\,.
	\end{eqnarray*}	
	Then \eqref{Der} follows by 
\eqref{derivopsi(x)}.
\\
Case $\rmdue$. 
We have two sub-cases: 
$\psi(x_2(E,\II_2);E,\II_2)=0\,,
\psi(x_1(E,\II_2);E,\II_2)=\pi$
or
$\psi(x_2(E,\II_2);E,\II_2)=\pi\,,
\psi(x_1(E,\II_2);E,\II_2)=0$.
In the first sub-case by the first formula in \eqref{AreaII}
we 
have\footnote{For brevity
we omit to write the dependence on $E$ and $\II_2$.}
	\begin{eqnarray*}
		\partial_E \Ac^\rmdue(E)&=&
		\partial_E x_1+
		\frac 1\pi \left[\psi(x_2)
		\partial_E 
		x_2-\psi(x_1)\partial_E x_1\right]+\frac 1\pi\int_{x_1}^{x_2}\partial_E \psi\,dx\\
		&=&\frac 1\pi\int_{x_1}^{x_2}\partial_E \psi\,dx
		\stackrel{\eqref{derivopsi(x)}}=-\int_{x_1}^{x_2}
		W(x)\,dx\,,\,,
	\end{eqnarray*}
and, analogously,
\begin{eqnarray*}
		\partial_{\II_2} \Ac^\rmdue(E)&=&
		\partial_{\II_2} x_1+
		\frac 1\pi \left[\psi(x_2)
		\partial_{\II_2} 
		x_2-\psi(x_1)\partial_{\II_2}
		 x_1\right]+
		 \frac 1\pi\int_{x_1}^{x_2}\partial_{\II_2} \psi(x)\,dx\\
		&=&\frac 1\pi\int_{x_1}^{x_2}\partial_{\II_2} \psi(x)\,dx\,.
	\end{eqnarray*}	
In the second sub-case by the second formula in \eqref{AreaII}
we 
have
	\begin{eqnarray*}
		\partial_E \Ac^\rmdue(E)&=&
		\partial_E x_2-
		\frac 1\pi \left[\psi(x_2)
		\partial_E 
		x_2-\psi(x_1)\partial_E x_1\right]-\frac 1\pi\int_{x_1}^{x_2}\partial_E \psi\,dx\\
		&=&-\frac 1\pi\int_{x_1}^{x_2}\partial_E \psi\,dx
		\stackrel{\eqref{derivopsi(x)}}=\int_{x_1}^{x_2}
		W(x)\,dx\,,\,,
	\end{eqnarray*}
and, analogously,
\begin{eqnarray*}
		\partial_{\II_2} \Ac^\rmdue(E)&=&
		\partial_{\II_2} x_1-
		\frac 1\pi \left[\psi(x_2)
		\partial_{\II_2} 
		x_2-\psi(x_1)\partial_{\II_2}
		 x_1\right]-
		 \frac 1\pi\int_{x_1}^{x_2}\partial_{\II_2} \psi(x)\,dx\\
		&=&-\frac 1\pi\int_{x_1}^{x_2}\partial_{\II_2} \psi(x)\,dx\,.
	\end{eqnarray*}		
	We conclude by \eqref{derivopsi(x)}.
	\\
Case $\rmtre$. Since 
$\psi(x_2(E,\II_2);E,\II_2)
=\psi(x_1(E,\II_2);E,\II_2)=\pi$, by \eqref{AreaIII} we 
have
	\begin{eqnarray*}
		\partial_E \Ac^\rmtre(E)&=&
		\frac 1\pi \left[\big(\pi-\psi(x_2)\big)
		\partial_E 
		x_2-\big(\pi-\psi(x_1)\big)\partial_E x_1\right]
		-\frac 1\pi\int_{x_1}^{x_2}\partial_E \psi(x)\,dx\\
		&=&-\frac 1\pi\int_{x_1}^{x_2}\partial_E \psi(x)\,dx
		\stackrel{\eqref{derivopsi(x)}}=\int_{x_1}^{x_2}
		W(x)\,dx\,,\,,
	\end{eqnarray*}
and, analogously,
\begin{eqnarray*}
		\partial_{\II_2} \Ac^\rmtre(E)&=&
		\frac 1\pi \left[\big(\pi-\psi(x_2)\big)
		\partial_{\II_2} 
		x_2-\big(\pi-\psi(x_1)\big)\partial_{\II_2} x_1\right]
		-\frac 1\pi\int_{x_1}^{x_2}\partial_{\II_2} \psi(x)\,dx\\
		&=&-\frac 1\pi\int_{x_1}^{x_2}\partial_{\II_2} \psi(x)\,dx\,.
	\end{eqnarray*}	
	Again we conclude by \eqref{derivopsi(x)}.
	\\
Case $\rmquattro$ is analogous to case $\rmdue$
sending $1\to 3$ and $2\to 4$.	
\eproof

\begin{rem}\label{exactresonance}
In the case of exact 3:1 resonance, namely when $\omega_+=3\omega_-$
the functions $F$ and $a$ in \eqref{sublime}, 
$\psi$ in \eqref{psi(x)}, $\mathbf P$
 in \eqref{penisola} with its roots $x_i$,   do not depend
on $\II_2$.
As a consequence the functions $W$ and  $\Ac$ in Proposition
\ref{WWW} do not depend
on $\II_2$.
In particular $\partial_{\II_2} 
				\Ac(E;\II_2)=0$,
				$\Ac(E;\II_2)=\Ac(E)$
and
formula \eqref{california} simplifies 
\begin{eqnarray}
\omega_-^{\rm nlr}(E,\II_2)
&:=&
\omega_-+
				\chi\,\II_2\left(  
				2(E+a_0)
				\right)
\,,
\nonumber
\\
\omega_+^{\rm nlr}(E,\II_2)
&:=&
\omega_+ + 3\chi\,\II_2\left(  
				2(E+a_0)
				+\frac{1}{\partial_E \Ac(E)}
				\right)\,.
\label{californiaBIS}				
\end{eqnarray}	
 
\end{rem}

Let us now practically evaluate the 
elliptic integrals\footnote{For a wide treatment
of elliptic integrals see, e.g., \cite{Elliptic}.} $\int W(x) dx$ and $\int x W(x) dx$ in \eqref{Der}.
Assume that 
 the polynomial
$\mathbf P$
 in \eqref{penisola} has 4 distinct 
 roots: $x_1,x_2,x_3,x_4$, 
  namely\footnote{$(1+a_2^2/4)$
  is the coefficient of the fourth order term
  of $\mathbf P(x)$.}
 \begin{equation}\label{francia}
\mathbf P(x)=\left(1+\frac{a_2^2}{4}\right)(x-x_1)
(x-x_2)(x-x_3)(x-x_4)\,.
\end{equation}
We have two cases: 
\\
i) the four roots are real, 
$x_1<x_2<x_3<x_4$;
\\
ii) we have two real roots, $x_1<x_2$
and two complex conjugated roots
$x_3=\bar x_4$.

\subsection{Elliptic integrals: the case of four real roots}
\label{Sec:elliptic}

Let us define the {\sl cross ratio}\footnote{Note that
$\lambda\neq 0,1,\infty$, since $x_j$, $j=1,2,3,4,$ are distinct.}:
\begin{equation}\label{cross}
 \lambda
      :=
      \frac{(x_2 - x_1) (x_4 - x_3)}{(x_3 - x_1)(x_4 - x_2) }\,.
\end{equation}
Note that $0<\lambda<1$.
Define the {\sl elliptic modulus}:
\begin{equation}\label{ellipticmodulus}
\kt
:=
\frac{1-\sqrt\lambda}{1+\sqrt\lambda}\,.
\end{equation}
Note that $0<\kt<1$.
 We now construct a change of variable
$x=\Tt(z)$
given by a M\"obius transformation
\begin{equation}\label{Mobius}
\Tt(z):=\frac{\At z+\Bt}{\Ct z+\Dt}\,,
\end{equation}
such that\footnote{As is well known
the cross ratio is invariant under 
M\"obius transformations.
Then, by \eqref{2.3Takebe},
we get
$
\lambda=\frac{(\kt-1)^2}{(\kt+1)^2}
$, which is consistent with 
\eqref{ellipticmodulus}.
See Lemma 2.3 and  
Exercise 2.4 of \cite{Elliptic}.}
\begin{equation}\label{2.3Takebe}
\Tt(-1/\kt)=x_4\,,\quad
\Tt(-1)=x_3\,,\quad
\Tt(1)=x_2\,,\quad
\Tt(1/\kt)=x_1\,.
\end{equation}
It is simple to show (see formula (2.7) of \cite{Elliptic})
that the transformation $x=\Tt(z)$ can be 
construct as the solution of equation
\begin{equation}\label{frodo}
\frac{(x-x_1)(x_3-x_4)}{(x-x_4)(x_3-x_1)}
=
\frac{(z-1/\kt)(-1+1/\kt)}{(z+1/\kt)(-1-1/\kt)}\,.
\end{equation}
Then the (real) coefficients of $\Tt$
are given by
\begin{eqnarray}
\At
&:=&
-\kt  x_1 x_3 - \kt^2 x_1 x_3 + 2 \kt x_1 x_4 - \kt x_3 x_4 + \kt^2 x_3 x_4\,,
\nonumber
\\
\Bt
&:=&
-x_1 x_3 - \kt x_1 x_3 + 2 \kt x_1 x_4 + x_3 x_4 - \kt x_3 x_4\,,
\nonumber
\\
\Ct
&:=&
\kt x_1 - \kt^2 x_1 - 2 \kt x_3 + \kt x_4 + \kt^2 x_4 \,,
\nonumber
\\
\Dt
&:=&
-x_1 + \kt x_1 - 2 \kt x_3 + x_4 + \kt x_4
\,.
\label{gandalf}
\end{eqnarray}
Note that, since $\kt>0$ and
$x_1<x_2<x_3<x_4$ we 
have\footnote{Indeed
$ x_1 - \kt x_1 - 2  x_3 +  x_4 + \kt x_4=0$
implies, by \eqref{ellipticmodulus}, that
 $\sqrt\lambda x_1-(1+\sqrt\lambda)x_3+x_4=0$,
 namely $\sqrt\lambda=(x_4-x_3)/(x_3-x_1)$.
 Squaring, by \eqref{cross},  we get 
 the right hand side of \eqref{svezia}.
 }
\begin{equation}\label{svezia}
\Ct= 0\qquad\iff\qquad
(x_2-x_1)(x_3-x_1)=(x_4-x_2)(x_4-x_3)\,.
\end{equation}
Note also that $\Tt$ is invertible
(on the Riemann sphere $\mathbb C\cup\{\infty\}$)  
and
$
\Tt(\mathbb R)=\mathbb R$.
Note that, since  $x_1<x_2<x_3<x_4$
and
$0<\kt<1$,
then
\begin{equation}\label{spagna}
\At \Dt-\Bt \Ct=
2 \kt (1 - \kt^2) (x_1 - x_3) (x_1 - x_4) (x_3 - x_4)<0\,.
\end{equation}
We have
\begin{equation}\label{spagna2}
\frac{d\Tt}{dz}(z)=
\frac{\At \Dt-\Bt \Ct}{(\Ct z+\Dt)^2}<0\,.
\end{equation}
Since
$$
\Tt(z)-\Tt(\zeta)=\big( \At-\Ct\,\Tt(\zeta) \big)
\frac{z-\zeta}{\Ct z+\Dt}\,,
$$
recalling \eqref{francia} and 
\eqref{2.3Takebe}, the substitution 
$x=\Tt(z)$ gives
 \begin{eqnarray}\label{francia2}
\mathbf P(\Tt(z))
&=&
(1+\frac{a_2^2}{4})[\Tt(z)-\Tt(1/\kt)]
[\Tt(z)-\Tt(1)][\Tt(z)-\Tt(-1)][\Tt(z)-\Tt(-1/\kt)]
\nonumber
\\
&=&
\ct
\frac{p_\kt(z)}{(\Ct z+\Dt)^4}
\end{eqnarray}
where
\begin{equation}\label{pk}
p_\kt(z):=(1-z^2)(1-\kt^2 z^2)\,,
\qquad
\ct:=(1+a_2^2/4)
\kt^{-2}\prod_{1\leq j\leq 4}
\big( \At-\Ct\,x_j \big)
\,.
\end{equation}
Note that
\begin{eqnarray*}
\prod_{1\leq j\leq 4}
\big( \At-\Ct\,x_j \big)
&=&
16\sqrt\lambda 
\big(1 + \sqrt\lambda\big)^{-7}
\big(-1 + \sqrt\lambda\big)^4
 (x_1 - x_3)^2 (x_1 - x_4)^2 (x_3 - x_4)^2  \cdot
   \\
&& \cdot  \Big((x_1 - x_2) (x_3 - x_4) + (x_1 - x_3) (x_2 - x_4) 
\sqrt\lambda\Big)\, >\, 0\,,
\end{eqnarray*}
which implies that $\ct>0$.

\medskip

By \eqref{W}, \eqref{spagna2}, \eqref{2.3Takebe}
and \eqref{francia2}
we get
\begin{equation}\label{queen}
\int_{x_1}^{x_2}W(x)dx=
\frac{\Bt\Ct-\At\Dt}{\pi\sqrt\ct}\int_1^{1/\kt} 
\frac{dz}{\sqrt{-p_\kt(z)}}
=
\frac{\Bt\Ct-\At\Dt}{\pi\sqrt\ct}\int_{-1/\kt}^{-1}
\frac{dz}{\sqrt{-p_\kt(z)}}
=\int_{x_3}^{x_4}W(x)dx
\,,
\end{equation}
where the second equality holds since 
$p_\kt(z)$ is even.
It remains to evaluate
$\int_1^{1/\kt} 
\frac{dz}{\sqrt{-p_\kt(z)}}$,
which is an elliptic integral. We get
the {\sl complete elliptic integral of the first kind}\footnote{Note that 
${\tt EllipticK}:(-\infty,1)\to(0,+\infty)$
is an analytic strictly increasing function
with $\lim_{x\to-\infty}{\tt EllipticK}(x)=0$
and $\lim_{x\to 1^-}{\tt EllipticK}(x)=+\infty$.}
\begin{equation}\label{EllipticK}
\int_1^{1/\kt} 
\frac{dz}{\sqrt{-p_\kt(z)}}
=
\int_0^1 \frac{ds}{\sqrt{(1-s^2)(1-m_1 s^2)}}
=:{\tt EllipticK}(m_1)\,, 
\qquad m_1:=1-\kt^2\,,
\end{equation}
by the change of variable
 $z=\frac{1}{\sqrt{1-m_1 s^2}}$.
 Note that, since $0<\kt<1$ we have
 $0<m_1<1$.
 By \eqref{queen} and \eqref{EllipticK}
 we get
\begin{equation}\label{queen2}
\int_{x_1}^{x_2}W(x)dx=\int_{x_3}^{x_4}W(x)dx
=
\frac{\Bt\Ct-\At\Dt}{\pi\sqrt\ct}
{\tt EllipticK}(1-\kt^2)\,.
\end{equation}
\medskip

Similarly
\begin{equation}\label{queenx}
\int_{x_1}^{x_2}x\,W(x)dx=
\frac{\Bt\Ct-\At\Dt}{\pi\sqrt\ct}\int_1^{1/\kt} 
\frac{\At z+\Bt}{\Ct z+\Dt}
\frac{dz}{\sqrt{-p_\kt(z)}}
\,.
\end{equation}	

We have two cases: $\Ct\neq 0$ and $\Ct=0$.
In the first case	setting
\begin{equation}\label{sicilia}
\at:=\frac{\Dt}{\Ct}\,,\qquad
\bt:=\frac{\Bt\Ct-\At\Dt}{\Ct^2}\,,
\end{equation}
we have
\begin{equation}\label{sardegna}
\int_1^{1/\kt} 
\frac{\At z+\Bt}{\Ct z+\Dt}
\frac{dz}{\sqrt{-p_\kt(z)}}
=
\frac{\At}{\Ct}\int_1^{1/\kt} 
\frac{dz}{\sqrt{-p_\kt(z)}}
+ \bt \int_1^{1/\kt} 
\frac{1}{z+\at}
\frac{dz}{\sqrt{-p_\kt(z)}}
\,.
\end{equation}
Note that the real number 
$\at$ satisfies $|\at|>1/\kt$.
Otherwise, by contradiction, assume that
$|\at|\leq 1/\kt$. Since $\Tt(-\at)=\infty$
(in the Riemann sphere), by 
\eqref{2.3Takebe} we have that $|\at|=|-\at|<1/\kt$.
Since the real function $\Tt(z)$ has a vertical
asymptote at $z=-\at$, has  
$\At/\Ct$ as
 horizontal asymptote
 and is decreasing (recall \eqref{spagna2})
in the intervals $(-\infty,-\at)$ and $(-\at,+\infty)$,
we have that 
$\Tt(-1/\kt)<\At/\Ct<\Tt(1/\kt)$.
Then by \eqref{2.3Takebe}
we obtain $x_4<x_1$, which is a contradiction.
We conclude that $|\at|>1/\kt$.

The first integral on the right hand	
	side of \eqref{sardegna} has been evaluated in \eqref{EllipticK}.
	Regarding the second one we have
	\begin{equation}\label{integrale_ellittico_4radici_1}
		\int_1^{1/\kt} 
		\frac{1}{z+\at}
		\frac{dz}{\sqrt{-p_\kt(z)}}
				=
		\int_1^{1/\kt} 
		\frac{z}{z^2-\at^2}
		\frac{dz}{\sqrt{-p_\kt(z)}}
		-
		\int_1^{1/\kt} 
		\frac{\at}{z^2-\at^2}
		\frac{dz}{\sqrt{-p_\kt(z)}}\,.
	\end{equation}	
We have
	\begin{equation}\label{integrale_ellittico_4radici_2}
		\int_1^{1/\kt} 
		\frac{z}{z^2-\at^2}
		\frac{dz}{\sqrt{-p_\kt(z)}}\stackrel{z^2=t}=
		\frac12
		\int_1^{1/\kt^2}
		\frac{1}{t-\at^2}
		\frac{dt}{\sqrt{(t-1)(1-\kt^2 t)}}
		=-\frac{\pi}{2\sqrt{(\at^2-1)(\at^2\kt^2-1)}}
	\end{equation}	
and, by the change of variable
 $z=\frac{1}{\sqrt{1-m_1 s^2}}$,
we obtain
	\begin{eqnarray}
		\int_1^{1/\kt} 
		\frac{\at}{z^2-\at^2}
		\frac{dz}{\sqrt{-p_\kt(z)}}
		&=&\at
		\int_0^1 \frac{1-m_1 s^2}{(1-\at^2 +m_1\at^2 s^2)\sqrt{(1-s^2)(1-m_1 s^2)}}\, ds \label{integrale_ellittico_4radici_3}
			\\
			&=&
			-\frac{1}{\at}
			\int_0^1 \frac{1}{\sqrt{(1-s^2)(1-m_1 s^2)}}\, ds \label{integrale_ellittico_4radici_4}
				\\
			&&+\frac{1}{m_1 \at^3}
			\int_0^1 \frac{1}{\frac{1}{n_1}-s^2}\frac{1}{\sqrt{(1-s^2)(1-m_1 s^2)}}\, ds \,,\label{integrale_ellittico_4radici_5}
	\end{eqnarray}	
where
$$
0<n_1:=\frac{m_1\at^2}{\at^2-1}
\stackrel{\eqref{EllipticK}}=
\frac{\at^2-\kt^2\at^2}{\at^2-1}
<1\,,
$$
since  $|\at|>1/\kt>1$.
Recalling  that
\begin{equation}\label{EllipticPi} 
\mathtt{EllipticPi} (n_1,m_1):=
\int_0^1 \frac{1}{1-n_1 s^2}\frac{1}{\sqrt{(1-s^2)(1-m_1 s^2)}}\, ds 
\end{equation}
is the 
{\sl complete elliptic integral of the third
  kind},
by
\eqref{EllipticK} and
 \eqref{queenx}-\eqref{integrale_ellittico_4radici_5}
 and noting that $\frac{\At}{\Ct}+\frac{\bt}{\at}
 =\frac{\Bt}{\Dt},$
we get
\begin{eqnarray}\label{birra}
&&\int_{x_1}^{x_2}x\,W(x)dx=
\\
&&
\frac{\Bt\Ct-\At\Dt}{\pi\sqrt\ct}
\left(
\frac{\Bt}{\Dt} {\tt EllipticK}(m_1)
-\frac{\pi \bt}{2\sqrt{(\at^2-1)(\at^2\kt^2-1)}}
-
\frac{n_1\bt}{m_1\at^3} 
\mathtt{EllipticPi} (n_1,m_1)
\right)\,.
\nonumber
\end{eqnarray}

Let us now consider the case $\Ct=0$.
By \eqref{queenx} and
since
$$
\int_1^{1/\kt} 
\frac{z\,dz}{\sqrt{-p_\kt(z)}}
\stackrel{z^2=t}=
		\frac12
		\int_1^{1/\kt^2}
		\frac{dt}{\sqrt{(t-1)(1-\kt^2 t)}}
		=\frac{\pi}{2 \kt}\,,
$$
we have
\begin{eqnarray}\label{queenxBIS}
\int_{x_1}^{x_2}x\,W(x)dx
&=&
-\frac{\At\Bt}{\pi\sqrt\ct}{\tt EllipticK}(m_1)
-
\frac{\At^2}{\pi\sqrt\ct}\int_1^{1/\kt} 
\frac{z\,dz}{\sqrt{-p_\kt(z)}}
\nonumber
\\
&=&
-\frac{\At\Bt}{\pi\sqrt\ct}{\tt EllipticK}(m_1)
-
\frac{\At^2}{2\kt\sqrt\ct}
\,,
\end{eqnarray}	
which is exactly the limit for $\Ct\to 0$ of 
\eqref{birra}.

\medskip

Let us now evaluate
\begin{equation}\label{queenx3x4}
\int_{x_3}^{x_4}x\,W(x)dx=
\frac{\Bt\Ct-\At\Dt}{\pi\sqrt\ct}\int_{-1/\kt}^{-1} 
\frac{\At z+\Bt}{\Ct z+\Dt}
\frac{dz}{\sqrt{-p_\kt(z)}}
\,.
\end{equation}
Changing variable $z\to -z$ we get
$$
\int_{-1/\kt}^{-1} 
\frac{\At z+\Bt}{\Ct z+\Dt}
\frac{dz}{\sqrt{-p_\kt(z)}}
=
\int_1^{1/\kt}
\frac{\At z-\Bt}{\Ct z-\Dt}
\frac{dz}{\sqrt{-p_\kt(z)}}
$$
When $\Ct\neq 0$, recalling \eqref{sicilia}, 
we have
\begin{equation}\label{sardegna2}
\int_1^{1/\kt} 
\frac{\At z+\Bt}{\Ct z+\Dt}
\frac{dz}{\sqrt{-p_\kt(z)}}
=
\frac{\At}{\Ct}\int_1^{1/\kt} 
\frac{dz}{\sqrt{-p_\kt(z)}}
- \bt \int_1^{1/\kt} 
\frac{1}{z-\at}
\frac{dz}{\sqrt{-p_\kt(z)}}
\,.
\end{equation}
Reasoning as in derivation of \eqref{birra}
we get
\begin{eqnarray}\label{birra2}
&&\int_{x_3}^{x_4}x\,W(x)dx=
\\
&&
\frac{\Bt\Ct-\At\Dt}{\pi\sqrt\ct}
\left(
\frac{\Bt}{\Dt} {\tt EllipticK}(m_1)
+\frac{\pi \bt}{2\sqrt{(\at^2-1)(\at^2\kt^2-1)}}
-
\frac{n_1\bt}{m_1\at^3} 
\mathtt{EllipticPi} (n_1,m_1)
\right)\,.
\nonumber
\end{eqnarray}
The case $\Ct=0$ can be obtained taking the 
limit for $\Ct\to 0$ of 
\eqref{birra2}.

\subsection{Elliptic integrals: the case of two real roots}\label{somorta}

In this case define the  cross ratio
and the elliptic modulus
as:\footnote{Setting $w:=(x_1 - x_4)(x_2 - x_3)$
we have that $ \lambda_*=\bar w/w$
since $x_1,x_2\in\mathbb R$ and $\bar x_3=x_4$.
Then $\sqrt\lambda_*:=|w|/w$  satisfies
$(\sqrt\lambda_*)^2=
|w|^2/w^2=\lambda_*$.}
\begin{equation}\label{crossBIS}
 \lambda_*
      :=
      \frac{(x_1 - x_3) (x_2 - x_4)}{(x_1 - x_4)(x_2 - x_3) }\,,\quad
 \kt_*:=\frac{1-\sqrt \lambda_*}{1+\sqrt \lambda_*}
 \,,
 \quad 
 \sqrt \lambda_*:=
 \frac{|x_1 - x_4| |x_2 - x_3|}{(x_1 - x_4)(x_2 - x_3) }
     \,.
\end{equation}
Since $|\sqrt \lambda_*|=1$ there exists a real
$\theta$ such that $\sqrt\lambda_*=e^{\mathrm i \theta}$,
so that
 $\kt_*=-\mathrm i \tan(\theta/2)$,
namely $\kt_* $ is purely imaginary
and $  \kt_*^2<0$
(see page 40
of \cite{Elliptic} for details).	
	 We now construct 
 a M\"obius transformation
 \begin{equation}\label{MobiusBIS}
 \Tt_*(z):=\frac{\At_* z+\Bt_*}{\Ct_* z+\Dt_*}\,,
\end{equation}
such that
\begin{equation}\label{2.3TakebeBIS}
 \Tt_*(-1/ \kt_*)=x_4\,,\quad
 \Tt_*(-1)=x_2\,,\quad
 \Tt_*(1)=x_1\,,\quad
 \Tt_*(1/ \kt_*)=x_3\,.
\end{equation}
It is simple to show (see formula (2.7) of \cite{Elliptic})
that the transformation $x= \Tt_*(z)$ can be 
construct as the solution of equation
\begin{equation}\label{frodoBIS}
\frac{(x-x_3)(x_2-x_4)}{(x-x_4)(x_2-x_3)}
=
\frac{(z-1/ \kt_*)(-1+1/ \kt_*)}{(z+1/ \kt_*)(-1-1/ \kt_*)}\,.
\end{equation}
Note that $\Tt_*$ is invertible
(on the Riemann sphere $\mathbb C\cup\{\infty\}$)  
and
$
 \Tt_*(\mathbb R)=\mathbb R$.
Indeed the last claim is equivalent to show that
if $x\in\mathbb R$ in \eqref{frodoBIS}
then also $z\in\mathbb R$.
This can be proven 
taking the complex
conjugate of \eqref{frodoBIS}
 and inverting both sides\footnote{More precisely
 denoting by $\ell$ and $r$, respectively,
 the left and right hand side of 
 \eqref{frodoBIS}, we have that, if
 $x\in\mathbb R$ then $\ell=1/\bar\ell$
 (recall $\bar x_3=x_4$),
 which implies $r=1/\bar r$
 (recall $\bar\kt_*=-\kt_*$), namely
 $\frac{z-a}{z+a}=
 \frac{\bar z-a}{\bar z+a}$
 denoting for brevity $a:=1/\kt_*$.
 Then $z=\bar z$, namely $z\in\mathbb R$. }.
The  coefficients of $\Tt_*$, which
are given by
\begin{eqnarray}
\At_*
&:=&
- x_2 (x_3+x_4) + \kt_* x_2(x_4- x_3)   
 + 2  x_3 x_4
\,,
\nonumber
\\
\Bt_*
&:=&
 -  x_2 (x_3+x_4) 
+ x_2 (x_4-x_3)/\kt_*
+ 2  x_3 x_4
\,,
\nonumber
\\
\Ct_*
&:=&
-2  x_2 +  x_3  +  x_4 + \kt_* (x_4-x_3)
 \,,
\nonumber
\\
\Dt_*
&:=&
-2  x_2 +  x_3+  x_4 + (x_4-x_3)/\kt_* 
\,,
\label{gandalfBIS}
\end{eqnarray}
are real since, $x_2\in\mathbb R,$
$\bar x_3=x_4$ and $\kt_*$
is purely imaginary. 
We have that 
\begin{equation}\label{spagna2BIS}
\frac{d\Tt_*}{dz}(z)=
\frac{\At_* \Dt_*-\Bt_* \Ct_*}{(\Ct_* z+\Dt_*)^2}<0
\qquad
\mbox{for}\qquad
z\in\mathbb R\,,
\end{equation}
since 
$
 \Tt_*(-1)=x_2> \Tt_*(1)=x_1
 $ by \eqref{2.3TakebeBIS}.
 It follows that
 \begin{equation}\label{spagnaBIS}
\At_* \Dt_*-\Bt_* \Ct_*=
2 (\kt_*^2-1) 
(x_2 - x_3) (x_2 - x_4) (x_3 - x_4)
/\kt_*<0\,.
\end{equation}
Arguing as in \eqref{francia2},
 the substitution 
$x= \Tt_*(z)$ gives
 \begin{equation}\label{francia2BIS}
\mathbf P( \Tt_*(z))
=
-\ct_*
\frac{p_{\kt_*}(z)}{(\Ct_* z+\Dt_*)^4}
\end{equation}
where $p_{\kt_*}(z):=(1-z^2)(1-\kt_*^2 z^2)$
and
\begin{equation}\label{vallinsu}
\ct_*:=-(1+a_2^2/4)
\kt_*^{-2}\prod_{1\leq j\leq 4}
\big( \At_*-\Ct_*\,x_j \big)\,.
\end{equation}
Note that $\ct_*>0$; 
indeed
$\kt_*^{-2}<0$,  $(\At_*-\Ct_* x_3)
(\At_*-\Ct_* x_4)=|\At_*-\Ct_* x_3|^2>0 $
(since\footnote{Note that $\At_*-\Ct_* x_j\neq 0$
for $j=1,2,3,4$, since $\At_*-\Ct_* x_j=0$
implies $\Tt_*^{-1}(x_j)=\infty$ that contradicts
\eqref{2.3TakebeBIS}.
} $\At_*,\Ct_*\in\mathbb R$ and
$\bar x_3=x_4$), finally,
denoting for brevity $w:=(x_1-x_4)(x_2-x_3)$,
we have 
\begin{eqnarray*}
&&
(\At_*-\Ct_* x_1)
(\At_*-\Ct_* x_2)
=4 (x_2-x_3)^2 (x_1-x_4)(x_2-x_4)
\frac{1+\bar w/|w|}{1+ w/|w|}
\\
&&
=4 |x_2-x_3|^2 w\frac{1+\bar w/|w|}{1+ w/|w|}
=4 |x_2-x_3|^2 |w|>0\,.
\end{eqnarray*}
By \eqref{W}, \eqref{spagna2BIS}, 
\eqref{2.3TakebeBIS}
and \eqref{francia2BIS}
we get
\begin{equation}\label{queenBIS}
\int_{x_1}^{x_2}W(x)dx
=
\frac{\Bt_*\Ct_*-\At_*\Dt_*}{\pi\sqrt{\ct_*}}
\int_{-1}^{1}
\frac{dz}{\sqrt{p_\kt(z)}}
=
2\frac{\Bt_*\Ct_*-\At_*\Dt_*}{\pi\sqrt{\ct_*}}
\int_{0}^{1}
\frac{dz}{\sqrt{p_{\kt_*}(z)}}
\,,
\end{equation}
since $p_{\kt_*}(z)$ is even; in particular
\begin{equation}\label{EllipticKBIS}
\int_0^{1} 
\frac{dz}{\sqrt{p_{\kt_*}(z)}}
=
\int_0^{1} 
\frac{dz}{\sqrt{(1-z^2)(1-\kt_*^2 z^2)}}
=:{\tt EllipticK}(\kt_*^2)\,.
\end{equation}
 By \eqref{queenBIS} and \eqref{EllipticKBIS}
 we get
\begin{equation}\label{queen2BIS}
\int_{x_1}^{x_2}W(x)dx
=
2\frac{\Bt_*\Ct_*-\At_*\Dt_*}{\pi\sqrt{\ct_*}}
{\tt EllipticK}(\kt_*^2)\,.
\end{equation}
Arguing as in \eqref{queenBIS}	
and recalling the definition of $p_{\kt_*}(z)$ in \eqref{pk}
we obtain
\begin{equation}\label{queenBIS2}
\int_{x_1}^{x_2} x W(x)dx
=
\frac{\Bt_*\Ct_*-\At_*\Dt_*}{\pi\sqrt{\ct_*}}
\int_{-1}^{1}
\frac{\At_* z+\Bt_*}{\Ct_* z+\Dt_*}
\frac{dz}{\sqrt{p_{\kt_*}(z)}}
=
\frac{\Bt_*\Ct_*-\At_*\Dt_*}{\pi\sqrt{\ct_*}}
\int_{-1}^{1}
\Tt_*(z)
\frac{dz}{\sqrt{p_{\kt_*}(z)}}
\,,
\end{equation}
	Since the last integration interval is 
	symmetric 
and $p_{\kt_*}(z)$ is an even function
we can substitute $\Tt_*(z)$ with its even part,
namely 
$$
\frac12\big(\Tt_*(z)+\Tt_*(-z)\big)=
\frac{\At_*\Ct_* z^2-\Bt_*\Dt_*}{\Ct_*^2 z^2 -\Dt_*^2}
=\frac{\At_*}{\Ct_*}+
\frac{\Bt_*\Ct_*-\At_*\Dt_*}{\Ct_*\Dt_*}
\frac{1}{1-(\Ct_*/\Dt_*)^2z^2}\,,
$$
obtaining (since the integrands are even)
$$
\int_{-1}^{1}
\Tt_*(z)
\frac{dz}{\sqrt{p_{\kt_*}(z)}}
=
2\frac{\At_*}{\Ct_*}
\int_{0}^{1}
\frac{dz}{\sqrt{p_{\kt_*}(z)}}
+
2\frac{\Bt_*\Ct_*-\At_*\Dt_*}{\Ct_*\Dt_*}
\int_{0}^{1}
\frac{1}{1-(\Ct_*/\Dt_*)^2z^2}
\frac{dz}{\sqrt{p_{\kt_*}(z)}}\,.
$$
Then, by \eqref{EllipticK} and \eqref{EllipticPi}
\begin{equation}\label{goldberg}
\int_{-1}^{1}
\Tt_*(z)
\frac{dz}{\sqrt{p_{\kt_*}(z)}}
=
2\frac{\At_*}{\Ct_*}
{\tt EllipticK}(\kt_*^2)
+
2\frac{\Bt_*\Ct_*-\At_*\Dt_*}{\Ct_*\Dt_*}
\mathtt{EllipticPi} (\Ct_*^2 \Dt_*^{-2},\kt_*^2)\,.
\end{equation}
Recalling \eqref{Der}, \eqref{queen2BIS}, \eqref{queenBIS2}, \eqref{goldberg}, in the case of two real roots, the last term in 
\eqref{california} writes
\begin{equation}\label{germania}
\chi\, \II_2^2\,
				\frac{\partial_{\II_2} 
				\Ac(E;\II_2)}{\partial_{E} \Ac(E;\II_2)}
=
\frac{\omega_+-3\omega_-}{3}	
\left(
\frac{\At_*}{\Ct_*}
+
\frac{\Bt_*\Ct_*-\At_*\Dt_*}{\Ct_*\Dt_*}
\frac{\mathtt{EllipticPi} (\Ct_*^2 \Dt_*^{-2},\kt_*^2)}{{\tt EllipticK}(\kt_*^2)}
\right)\,.
\end{equation}

\subsection{Explicit expression of the 
nonlinear frequencies for the exact 3:1 resonance}\label{sec:Explicit}

In this subsection we consider only
 the case of exact 3:1 resonance, namely when $\omega_+=3\omega_-$.
 Let the energy $E$ be such that	
 the polynomial
$\mathbf P$
 in \eqref{penisola}
  has 4 distinct 
 roots: $x_1(E),x_2(E),$ $x_3(E),x_4(E).$	
Recalling the definitions of $\kt$ in
\eqref{ellipticmodulus}, of $\At,\Bt,\Ct,\Dt$
in
\eqref{gandalf}, of $\ct$ in \eqref{pk},
of $\kt_*$ in
\eqref{crossBIS}, of $\At_*,\Bt_*,\Ct_*,\Dt_*$
in
\eqref{gandalfBIS}, of $\ct_*$ in \eqref{vallinsu},
note that all these quantities depend on $E$.
Recalling Proposition \ref{WWW}, \eqref{Der}, \eqref{DerBIS},
\eqref{queen2} and \eqref{queen2BIS},
formula \eqref{californiaBIS} in Remark \ref{exactresonance} becomes 
\begin{eqnarray}
\omega_-^{\rm nlr}(E,\II_2)
&:=&
\omega_-+
				2\chi\,\II_2 
				(E+a_0)
\,,
\nonumber
\\
\omega_+^{\rm nlr}(E,\II_2)
&:=&
\omega_+ + 6\chi\,\II_2\left(  
				E+a_0+ V(E)
				\right)\,,
\label{californiaTER}				
\end{eqnarray}	
where the function $V(E)$ is defined as follows:
	\begin{equation}\label{loto1}
	V(E)
:=\pm\frac{\pi\sqrt{\ct_*}}
{2(\Bt_*\Ct_*-\At_*\Dt_*){\tt EllipticK}(\kt_*^2)}
	\end{equation}
with the $+$ sign in the zones
$\PP_{01},$
$\PP_{21}^{+,\rmdue}$, $\PP_{21}^{+,\rmtre}$,
$\PP_{21}^{-,\rmtre}$, $\PP_{12}^{+,\rmtre}$,
and with 	$-$ sign in the zones
$\PP_{10},$
$\PP_{21}^{-,\rmuno}$, $\PP_{12}^{+,\rmuno}$,
$\PP_{12}^{-,\rmuno}$, $\PP_{12}^{-,\rmdue}$,
	moreover
	\begin{equation}\label{loto2}
	V(E):=
\pm \frac{\pi\sqrt\ct}{(\Bt\Ct-\At\Dt){\tt EllipticK}(1-\kt^2)}
\end{equation}
with the $+$ sign in the zones
$\PP_{21}^{+,\rmquattro}$, $\PP_{21}^{-,\rmquattro}$,
$\PP_{12}^{+,\rmdue}$, $\PP_{12}^{-,\rmtre}$,
and with 	$-$ sign in the zones	
$\PP_{21}^{+,\rmuno}$, $\PP_{21}^{-,\rmdue}$,
$\PP_{12}^{+,\rmquattro}$, $\PP_{12}^{-,\rmquattro}$.
\\
Note that, recalling \eqref{sublime}, in \eqref{californiaTER}	 we have that
$$
\chi a_0={\mathtt G}_{(2,0),(2,0)}\,.
$$
Then we can rewrite \eqref{californiaTER}
as
\begin{eqnarray}
\omega_-^{\rm nlr}(E,\II_2)
&:=&
\omega_-+
				2\II_2 
				(\chi E+{\mathtt G}_{(2,0),(2,0)})
\,,
\nonumber
\\
\omega_+^{\rm nlr}(E,\II_2)
&:=&
\omega_+ + 6\II_2\left(  
				\chi E+{\mathtt G}_{(2,0),(2,0)}+ 
				\chi V(E)
				\right)\,.
\label{californiaQUATER}				
\end{eqnarray}
Finally we can see the nonlinear resonant frequencies
as functions of the initial amplitudes
$a_-$ and $a_+$.
By \eqref{ampiezze} and
\eqref{salmone} we get
\begin{equation}\label{vlad}
J_1(0)=\frac12 \omega_+ a_+^2\,,\qquad
J_2(0)=\frac12(\omega_- a_-^2 
+3 \omega_+ a_+^2)\,,\qquad
\psi_1(0)=\psi_2(0)=0\,.
\end{equation}
By
\eqref{Idue}
we have
\begin{equation}\label{fortuna}
\II_2=\II_2(0)=\frac12(\omega_- a_-^2 
+3 \omega_+ a_+^2)
\end{equation}
and by \eqref{loacker}
and
 \eqref{sublime}
we get
\begin{equation}\label{sublime2}
				E=F\Big(3J_1(0)/\Ju(0),\psi_1(0);\Ju(0)\Big)
=
a(x_\dag;\II_2)+b(x_\dag)\,,\quad
\mbox{with}\quad
x_\dag:=\frac{3\omega_+ a_+^2}{\omega_- a_-^2 
+3 \omega_+ a_+^2}\,.
\end{equation}

\section{Nonlinear bandgap for the honeycomb metamaterial}
\label{sec:bandgap}

In this section we 
present some outcomes
 of our analysis and discuss its application to the 
honeycomb metamaterial
described in the introduction.
In particular we investigate
 the effect of nonlinearity
on the bandgap size,
highlighting the differences between the resonant and non resonant cases. First, we briefly recall what we proved
in \cite{DL}.


For a given pair 
$(\tilde M,\tilde K)$,
the bandgap is defined as
 the interval between the maximum of the acoustic
 frequency and the minimum of the optical frequency
 as the wave numbers run over the 
 Brillouin triangle.
In the linear case,
since the gradients of $\omega_-$
and $\omega_+$ (with respect to $(\tilde{k}_1,\tilde{k}_2)$)
	never vanish
	 in the interior of $\triangle$,
	maxima and minima are attained on 
the boundary	$\partial \triangle$.
In particular, 
 for every pair $(\tilde M,\tilde K)$, 
 the maximum of the linear acoustic frequency 
is attained at ${\bf X}$, while the minimum of
the linear optical frequency 
is attained at ${\bf \Gamma}$.
We anticipate that, in evaluating the {\sl nonlinear} bandgap, the point 
${\bf X}$ plays a crucial role, more important than ${\bf \Gamma}$.
Indeed, typically, in the set of parameters we are considering,
namely the rectangle  $[0.05, 0.3]\times[1, 20]$ in the
 $(\tilde M,\tilde K)$-plane,  the displacement  of the maximum
of the acoustic frequency due to the nonlinearity is more relevant than that 
of the minimum of the optical frequency.

\subsubsection*{Resonant parameters}

\begin{figure}[h!]
			\centering			
\includegraphics[width=8cm,height=91cm,keepaspectratio]{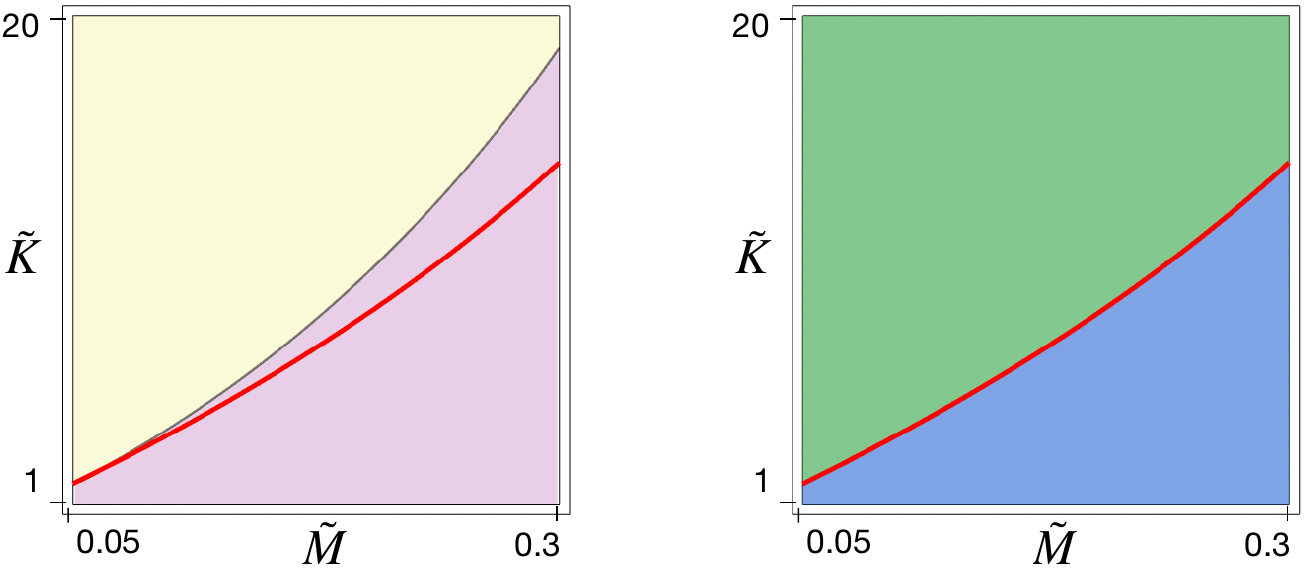}
			\caption{(Left) The reference rectangle
			$[0.05, 0.3]\times[1, 20]$ with the curve
			$\mathcal R$
			 plotted in red. The light yellow region 
			indicates 
			 the  non resonant set where the  representation formula \eqref{formulaSL1}
			is valid.
	In contrast, the light purple region
	shows the  set of resonant parameters,
	where  
	formula \eqref{californiaQUATER} is needed.
	(Right)	The reference rectangle
			$[0.05,0.3]\times[1,20]$. 
			For every pair $(\tilde M,\tilde K)$
			in the blue region, there exist two resonant curves intersecting $\triangle$ (as shown in Figure \ref{Brillouin_k1_k2}), while for pairs in the green region, there exists only one resonant curve intersecting $\triangle$.
						The red curve 
						$\mathcal R$ separates the two regions.	 }
			 \label{buoni}
		\end{figure}

As in
 \cite{DL}, 
within the reference  rectangle $[0.05, 0.3]\times[1, 20]$, we identify the curve
 $\mathcal R$ formed by the  pairs
  $(\tilde M,\tilde K)$ such that 
the linear acoustic and optical frequencies evaluated at $(\tilde k_1,\tilde k_2)={\bf X}$
are in 3:1 resonance, namely satisfy
$3 \omega_-=\omega_+$.
$\mathcal R$ is shown in Figure 
\ref{buoni}.
In \cite{DL}, we identify the set of
 nonresonant pairs $(\tilde M,\tilde K)$ 
within  the rectangle $[0.05, 0.3]\times[1, 20]$ (represented by the light yellow region in Figure
\ref{buoni} (left)), for which
 the maximum/minimum of the nonlinear acoustic/optical
 frequencies on the boundary
of the Brillouin triangle
 are attained at non resonant
wave numbers $(\tilde k_1,\tilde k_2)$, i.e.
at points where the quantity  $|3 \omega_--\omega_+|$
is not small.

Formula \eqref{formulaSL1} is valid in this  nonresonant set, 
allowing us to directly evaluate the bandgap in \cite{DL}.
In contrast,  in the complementary light purple zone in 
Figure
\ref{buoni}, formula \eqref{formulaSL1} is not applicable due to resonances and one has to use \eqref{californiaQUATER}
as we will show here.

\begin{figure}[h!]
			\centering			
			\includegraphics[width=7cm,height=7cm,keepaspectratio]{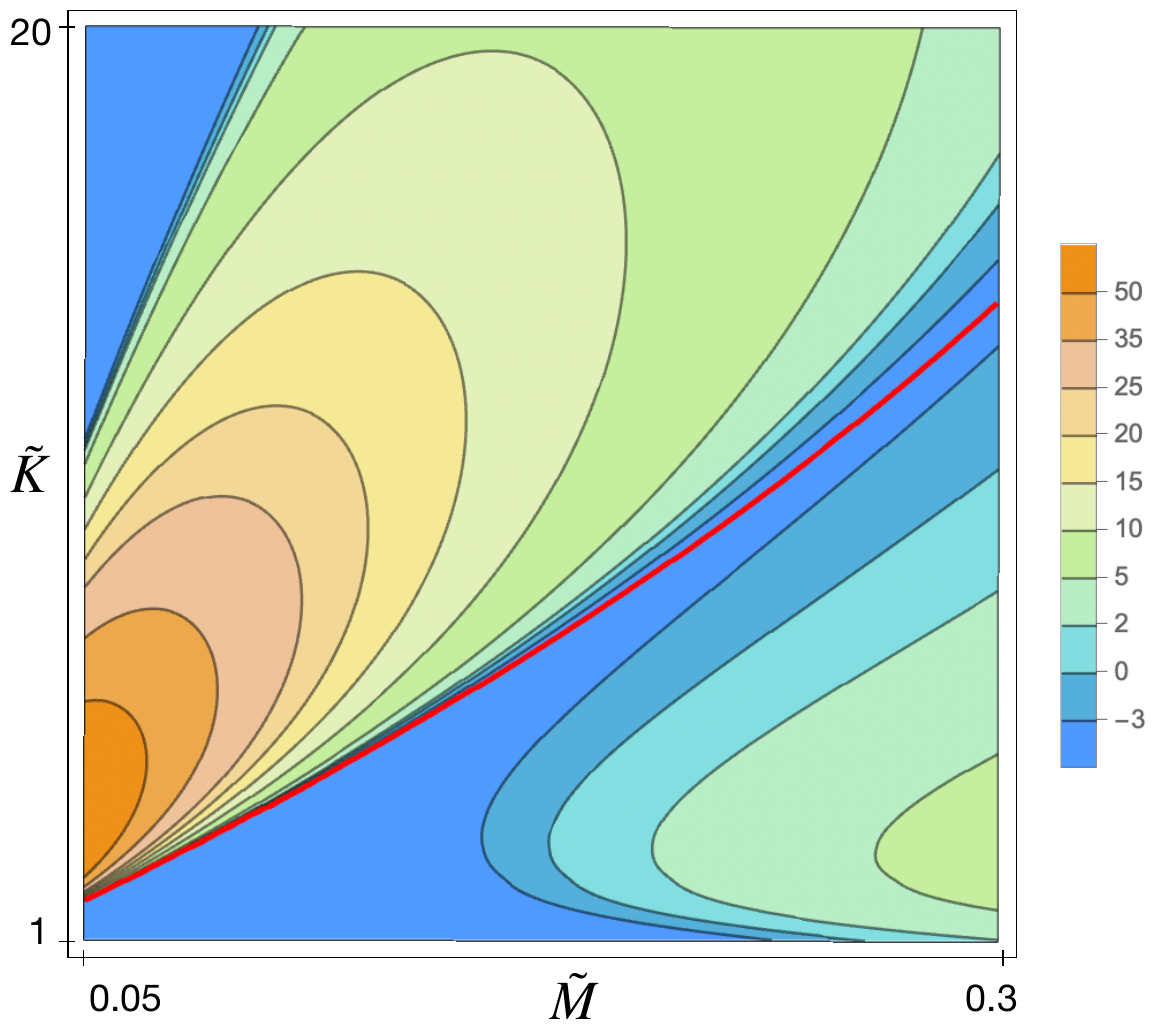} 
			\caption{
		Level curves of the maximum percentage 		difference between the nonlinear and linear bandgap are shown 
in the $(\tilde M,\tilde K)$-plane.
The  curve  $\mathcal R$
is plotted in red.
			Here,  $N_3=- 10^4$ (softening).
			In this softening case,
			the majority of parameter pairs above $\mathcal R$
			 result in
			an increase in the bandgap width,
			 while those below $\mathcal R$ either show a decrease or, at most, a very slight increase.
			 The region where the increase is most pronounced closely coincides with the set of nonresonant parameters highlighted in light yellow in	Figure \ref{buoni}.}
			 \label{ritz}
					\end{figure}

 The final result of our analysis is presented in 
  Figure \ref{ritz}, where the
 maximum percentage increment between the nonlinear and linear
 bandgap\footnote{Namely
 $100\times(W^{\rm nl}/W-1)$, where
 $W^{\rm nl}$ and $W$ denote the width of the nonlinear and linear bandgap, respectively.} is plotted as  the pair $(\tilde M,\tilde K)$
 varies over the rectangle $[0.05, 0.3]\times[1, 20]$
 in the softening case ($N_3=-10^4$).
We emphasize  that, while in \cite{DL}
we derived  Figure \ref{ritz} using \eqref{formulaSL1} only for 
 the pairs $(\tilde M,\tilde K)$ belonging
 to the light yellow set in Figure \ref{buoni},
 in this section, we show how to derive it
in the light purple set by \eqref{californiaQUATER}.

\begin{figure}[h!]
			\centering	\includegraphics[width=12cm,height=12cm,keepaspectratio]{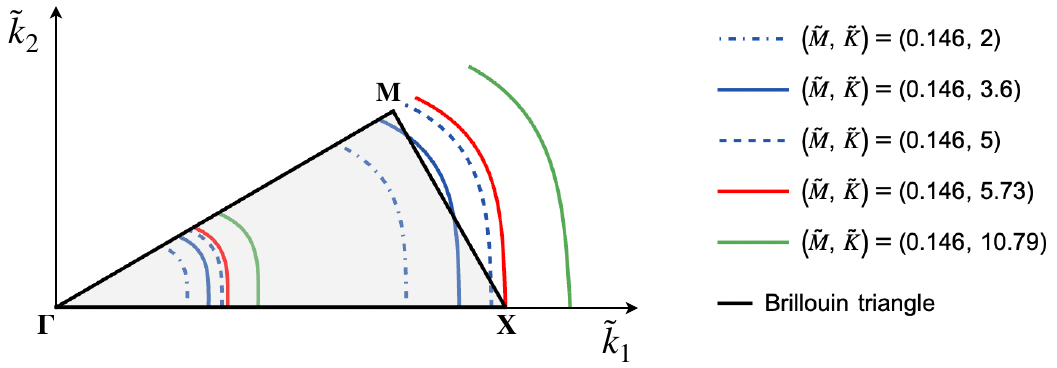}
			\caption{Resonant curves
in the $(\tilde{k}_1,\tilde{k}_2)$-plane	 for different fixed values of the pairs $(\tilde M,\tilde K)$ and their intersections with
the boundary of the Brillouin triangle $\triangle$.
The six blue curves  correspond to three
different
points $(\tilde M,\tilde K)$ in the blue region of
Figure \ref{buoni} (right).
 In particular  when 
$(\tilde M,\tilde K)=(0.146,2)$
the corresponding two blue curves, the
dot-dashed ones, 
 have 4 intersections. 
If 
$(\tilde M,\tilde K)=(0.146,3.6)$
the corresponding two blue curves, the
solid ones, 
 have 6 intersections.
 When 
$(\tilde M,\tilde K)=(0.146,5)$
the corresponding two blue curves, the
dashed ones, 
 have 4 intersections.
 The red curves, corresponding to
$(\tilde M,\tilde K)=(0.146,5.73)\in\mathcal R$, have 3 intersections.  Finally the green
curves, corresponding to
$(\tilde M,\tilde K)=(0.146,10.79)$, have only 2 intersections. 
			}\label{Brillouin_k1_k2}
		\end{figure}

Let us first recall how in \cite{DL} we identified
the two regions in Figure  \ref{buoni} (left).
Given a pair $(\tilde M,\tilde K)$, 
we define a set in the $(\tilde{k}_1,\tilde{k}_2)$-plane as {\sl resonant}
if every point in the set satisfies
the 3:1 resonance condition
 $3\omega_-
 (\tilde M,\tilde K,\tilde{k}_1,\tilde{k}_2)
 =\omega_+ (\tilde M,\tilde K,\tilde{k}_1,\tilde{k}_2)$.
For
a fixed pair
$(\tilde M,\tilde K)$ within the rectangle 
$[0.05, 0.3]\times[1, 20]$ (see Figure \ref{buoni}, (right))
 there are always one or two 
\textsl{resonant curves} in the
 $(\tilde{k}_1,\tilde{k}_2)$-plane,
that intersect
 the Brillouin triangle
$\triangle$ (see Figure \ref{Brillouin_k1_k2}).
The  curve $\mathcal R$
divides the rectangle $[0.05, 0.3]\times[1, 20]$ into two regions:  the one above and the one below $\mathcal R$, corresponding to
 the green region and the blue region in Figure
\ref{buoni} (right), respectively. 
For every fixed pair $(\tilde M,\tilde K)$ in the green region, there is only one 
resonant curve
in the plane of wave numbers $(\tilde k_1,\tilde k_2)$, that intersects the Brillouin triangle
(the green curve in Figure \ref{Brillouin_k1_k2}).
Conversely, for every fixed pair $(\tilde M,\tilde K)$ in the blue region, there  are two 
resonant curves
in the plane of wave numbers $(\tilde k_1,\tilde k_2)$, that intersect the Brillouin triangle
(the blue curves in Figure \ref{Brillouin_k1_k2}).
Finally, in the limit case when the pair $(\tilde M,\tilde K)$ belongs to the
curve $\mathcal R$,  there  are two 
resonant curves
in the $(\tilde k_1,\tilde k_2)$-plane, that intersect the Brillouin triangle, but one intersects
$\triangle$ only at ${\bf X}$
(see the red curves in Figure \ref{Brillouin_k1_k2}).

\subsubsection*{Admissible amplitudes}\label{sec:admissible}

Both in formula \eqref{formulaSL1}
and in formula \eqref{californiaQUATER},
(recall also \eqref{fortuna} and
\eqref{sublime2}), the nonlinear corrections
to the frequencies are essentially proportional
to the squares of the amplitudes $a_+$ and
$a_-$. Thus, the larger the amplitudes $a_\pm$, the greater the displacement of the nonlinear bandgap
relative to the linear one.
On the other hand, \eqref{formulaSL1} and
\eqref{californiaQUATER} are perturbative in nature,
as they are derived from the non resonant and resonant BNF, respectively.
Therefore,  $a_\pm$ must be sufficiently small for the
formulae to remain valid.
As shown in \cite{DL}, where they are analytically evaluated,
the ``admissible'' amplitudes are smaller
in the nonresonant case than in the resonant one.
Indeed, since the nonresonant BNF cancels more terms,
 it is ``stronger''
 than the resonant one. In particular 
 the admissible amplitudes in the nonresonant
 case approach zero  as the quantity $|3 \omega_- - \omega_+|$
vanishes.
For example, when taking $(\tilde k_1,\tilde k_2)=\mathbf X$,
the admissible amplitudes vanish for parameters values $(\tilde M, \tilde K)$
on the curve $\mathcal R$.
This is not the case of the admissible amplitudes 
in the resonant case, namely the ones appearing in  
 formulae \eqref{californiaQUATER}, \eqref{fortuna} and
\eqref{sublime2}).
Indeed they are bounded away from zero on the resonances.

\begin{figure}[h!]
\center
\includegraphics[width=16cm,height=16cm,keepaspectratio]{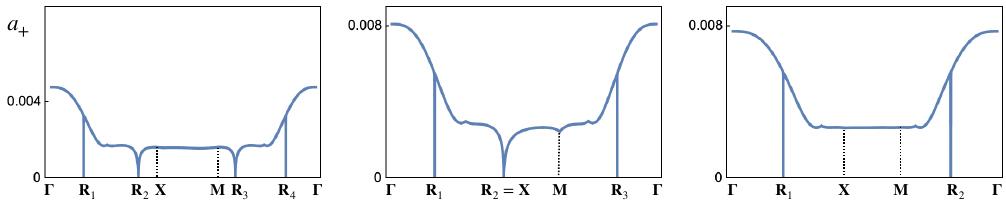}
\caption{Admissible 
initial amplitude (nonresonant case), $a_+$,  on the optical mode
as a function of the wave numbers 
 on  $\partial\triangle$ for $(\tilde M, \tilde K)=(0.146,2)$
 (on the left), $(\tilde M, \tilde K)=(0.146,5.73)$
 (in the middle)
 $(\tilde M, \tilde K)=(0.09,8)$
 (on the right).
 Here, $N_3=-10^4$.
 Note that the value
decreases to zero at
the four (on the left),    three (in the middle), and
two (on the right), resonant points (denoted by ${\bf R}_i$), respectively.
Referring to Figure 
\ref{Brillouin_k1_k2}, these four, three and two  points
correspond to the 
 intersections
of the dot-dashed blue, red, and green curves
 with $\partial\triangle$, respectively.
 The two points in the image  on the right are not
 shown in Figure 
\ref{Brillouin_k1_k2} but
they correspond to the same type of intersection that the green curves have with  $\partial\triangle$.
}\label{aPiuPar2}
\end{figure}

Shifting perspective, we can fix $(\tilde M, \tilde K)$
and observe at the variation of $a_\pm$
in the nonresonant case, 
as the wave numbers
$(\tilde k_1,\tilde k_2)$
 vary along the boundary $\partial\triangle$  of the Brillouin triangle 
$\triangle$. Notably,
  $a_\pm$
decreases to zero at certain resonant  points, denoted 
${\bf R}_i$. 
{\sl These points correspond to the intersections of
the boundary of the Brillouin triangle with 
the resonant curves plotted in  
 Figure \ref{Brillouin_k1_k2}.}
 Formula \eqref{formulaSL1} 
 loses validity in the vicinity of any point
 ${\bf R}_i$.
The values of the admissible 
initial amplitude $a_+$ (in the nonresonant case)
as $(\tilde k_1,\tilde k_2)$ traverses 
 $\partial\triangle$ are shown 
 in Figure \ref{aPiuPar2} for three different pairs of $(\tilde M, \tilde K)$.

\begin{figure}[h!]
\center
\includegraphics[width=16cm,height=16cm,keepaspectratio]{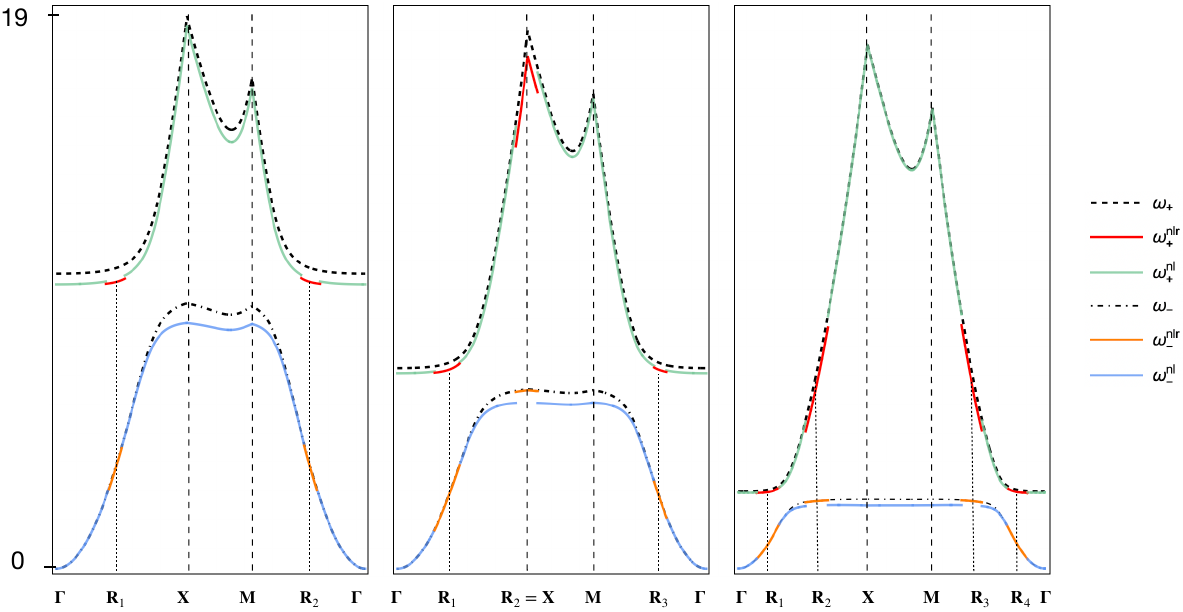}
			\caption{
			Linear, $\omega_\pm$,  and resonant,
			$\omega_\pm^{\rm nlr}$, as well as
			nonresonant, $\omega_\pm^{\rm nl}$, nonlinear dispersion curves versus wave numbers on $\partial\triangle$  for three different pairs of parameters
			$(\tilde M,\tilde K)$ in the softening case, $N_3=-10^4$.
			The points ${\bf R}_i$ are the resonant points, as seen in Figure
			\ref{aPiuPar2}.
			In  small  neighborhoods of these points, the expression
			$\omega_\pm^{\rm nl}$  is replaced by the resonant
			representation $\omega_\pm^{\rm nlr}$.
			Note that in all three cases, the minimum of the
			nonlinear optical frequency  essentially 
			coincides with $\omega_+^{\rm nl}$ evaluated at
			${\bf\Gamma}$.
			\\
			(On the left) Case (i):
			 $(\tilde M,\tilde K)=(0.09,8)$ belonging 
			 to the light yellow region in Figure \ref{buoni};  admissible
			 initial amplitudes  $a_-=0.0036$,
			 $a_+=0.0025$.
			 As in the linear case, the maximum of the acoustic
			 frequency is attained at ${\bf X}$, which is nonresonant.
			 The resulting percentage bandgap increment is around
			 $30\%$.
			\\
			(In the middle) Case (ii)
			$(\tilde M,\tilde K)=(0.146,5.73)$ belonging 
			 to the light purple region, more precisely to the red curve, in Figure \ref{buoni};  admissible
			 initial amplitudes  $a_-=0.002$,
			 $a_+=0.0012$.
			 As in the linear case the maximum of the acoustic
			 frequency is attained at ${\bf X}$, which, however, is now resonant. Since, at ${\bf X}$, $\omega_-^{\rm nlr}$ is very close to
			 $\omega_-$,
			 the  resulting nonlinear bandgap  undergoes a
			 slight decrement
			compared to the linear case.	
						 \\
(On the right) Case (iii)
$(\tilde M,\tilde K)=(0.2,1.1)$ belonging 
			 to the light purple region in Figure \ref{buoni};  
			 admissible
			 initial amplitudes  are $a_-=0.002$,
			 $a_+=0.001$.
			 Since $\omega_-$ is almost flat around its maximum, 
	 the maximum of the nonlinear acoustic
			 frequency is attained far from ${\bf X}$, more precisely
			 near the resonant point ${\bf R}_2$.
			  Since, close to ${\bf X}$, $\omega_-^{\rm nlr}$ is very close to
			 $\omega_-$,
			 the  resulting nonlinear bandgap  undergoes a
			  decrement
			 compared to the linear case.
			 }
			\label{casibrutti}
\end{figure}
In conclusion, due to the presence of the
3:1 resonance, formula \eqref{formulaSL1} becomes invalid
in the vicinity of the points ${\bf R}_i$,
when the parameters are resonant or nearly resonant. Specifically, this occurs when they
give rise to an exact, or nearly exact, 3:1 resonance between 
acoustic and optical  frequencies. 
In this resonant case the correct 
expression for the nonlinear frequencies  is
$\omega_\pm^{\rm nlr}$, 
as given by \eqref{californiaQUATER}.

\subsubsection*{Nonlinear bandgap}\label{sec:provola}

 Let us consider the softening case; the hardening case
can be treated analogously, 
leading to  a general decrement of the bandgap.
 We note that, since we are considering pairs $(\tilde M,\tilde K)$ 
 belonging to  the rectangle $[0.05, 0.3]\times[1, 20]$,
 the point ${\bf \Gamma}$, where the minimum of the linear acoustic
 frequency is attained, is always far from being resonant.  
Therefore, in the following discussion, we will focus on 
the maximum of acoustic frequency  because
it undergoes the most significant displacements
and may be resonant.
It turns out that, for the calculation of the nonlinear bandgap,
there are essentially three cases: 
\\
i) the maximum of the acoustic frequency
  and the minimum of the optical frequency
 are attained away from resonant
 points,
 \\
ii) $\mathbf X$ is resonant or nearly resonant,
\\
iii) $\omega_-$ has an almost flat maximum, 
 so that, even if
${\bf X}$ is away from resonance, 
the nonlinear acoustic frequency 
  may attain its maximum at 
some resonant (or nearly resonant) point away from ${\bf X}$.
\\
Note that
case i) corresponds to the 
light yellow region in Figure \ref{buoni},
while cases ii) and iii) correspond to the light purple
one.
These three cases are shown in Figure \ref{casibrutti}.

To summarize, one applies formula \eqref{formulaSL1} in case (i), as we did in \cite{DL}, and formula \eqref{californiaQUATER} in cases (ii) and (iii), as we do here. Using the expressions for the admissible amplitudes evaluated in \cite{DL}, we are able to compute the bandgap, thereby obtaining Figure \ref{ritz} in its entirety.


\section{Conclusions}

In this study, we investigated a broad range of structural engineering models by analyzing a general system of two coupled harmonic oscillators with cubic nonlinearity. Our examination revealed that, in the absence of damping, the system exhibits Hamiltonian dynamics, with an elliptic equilibrium at the origin characterized by two distinct linear frequencies.
 In particular, we focused on the resonant or nearly resonant case,
specifically when the two frequencies are close to a 3:1 resonance.
 \\
Our investigation involved employing
Hamiltonian  Perturbation Theory to transform the system into 
(resonant) Birkhoff Normal Form up to order 4. This transformation provided a new set of symplectic action-angle variables, 
on which the Hamiltonian, up to six-order terms, depends only
on the actions and the slow angle.  Notably, our analysis highlighted the dependency of the construction on the system's physical parameters, necessitating a meticulous case analysis of the phase portrait in the 3:1 resonant case. We found that
the system can exhibit up to six topologically different behaviors,
depending on the values of the physical parameters. 
 In each of these configurations,
 we described the nonlinear normal modes
 (elliptic/hyperbolic periodic orbits, invariant tori)
 and their stable and unstable
 manifolds of the truncated Hamiltonian (neglecting order six or higher terms).   This is a fundamental step for proving the
 persistence of the majority of these structures for the complete
 Hamiltonian by KAM Theory.

By using elliptic integrals, we derived explicit analytic formulas for the nonlinear frequencies.
 While this analytic expression was  already known
 away from resonances,
 it is, as far as we know,  new in this context
  for the  resonant or nearly resonant case.

As an application of our findings, we explored wave propagation in metamaterial honeycombs equipped with periodically distributed nonlinear resonators. 
Our investigation allowed us to examine the bandgap phenomenon in the presence of resonance. We found that while nonlinear effects far from resonances can significantly alter the bandgap, in the resonant case, the nonlinear frequencies, especially the acoustic one, closely align with the linear ones, resulting in a less pronounced variation in the bandgap.

\section{Appendix}

\subsection{Proof of Proposition \ref{pandivia}}

We first count the solutions of equation \eqref{gingerina2},
namely the intersections between the line 
$\ell(x):=a_2x+a_1$ and the function $b'(x)$
in \eqref{petrolio}.
We note that, since $b'$ is strictly convex,
  if $\ell(1)>0$ there is only one intersection.
  Note that condition $\ell(1)=a_2+a_1>0$ is equivalent
  to $a_2>-a_1$.
  Since $g(a_1)>-a_1$, condition
  $a_2>-a_1$ implies that we are in the zones
  $Z_{01}$ or $Z_{21}$ in which we have, indeed, 
  one intersection that we call
  $x_1^{(\pi)}$.
  \\
  Moreover, in this case $a_2>-a_1$, the function
  $x\to F(\pi,x)=a(x)-b(x)$ with $x\in(0,1)$ has only one critical point, which is exactly $x_1^{(\pi)}$.
  This critical point is a minimum
  since $\lim_{x\to 0^+}\partial_x F(\pi,x)
  =\lim_{x\to 0^+} a'(x)-b'(x)
  =-\infty$
  and $\lim_{x\to 1^-}\partial_x F(\pi,x)
  =\lim_{x\to 1^-} a'(x)-b'(x)
  =a_2+a_1>0$.
  \\
  Assume now that $a_2<-a_1$.
  Note that for every fixed $a_1\in\mathbb R$
   there exists a unique $a_2=h(a_1)$ such that
$a_2x+a_1$ is tangent to $b'(x)$ at some point
$0<x_0<1$.
In order to evaluate the function $h(a_1)$
above let us consider 
the tangent $r(x)$ in a point $x_0$ to 
$b'(x)$; namely:
	$$
		r(x)=b'(x_0)+b''(x_0)(x-x_0)\,.
	$$
Since we want that $r(x)=a_2 x+a_1$
we have to impose 
$r(0)=a_1$ and $b''(x_0)=a_2$.
Since
	\begin{equation}\label{mantova}
		b''(x)=\frac{8x^2-4x-1}{4x^{3/2}\sqrt{1-x}}\,,
	\end{equation}
imposing $r(0)=a_1$ we have
	\begin{eqnarray}
		a_1=r(0)&=&b'(x_0)-b''(x_0)x_0
		=\frac{(1-4x_0)\sqrt{1-x_0}}{2\sqrt{x_0}}-\frac{8x_0^2-4x_0-1}{4x_0^{3/2}\sqrt{1-x_0}}x_0
		\nonumber\\
		&=&\frac{2(1-4x_0)(1-x_0)-(8x_0^2-4x_0-1)}{4\sqrt{1-x_0}\sqrt{x_0}}
		\nonumber\\
	&=&\frac{2-8x_0-2x_0+8x_0^2-8x_0^2+4x_0+1}{4\sqrt{1-x_0}\sqrt{x_0}}
	\nonumber\\
		&=&\frac{-6x_0+3}{4\sqrt{1-x_0}\sqrt{x_0}}=\frac{3(1-2x_0)}{4\sqrt{1-x_0}\sqrt{x_0}}\,.
		\label{verona}
	\end{eqnarray}
Note that:
\begin{equation}\label{venezia}
a_1>0,<0,=0\qquad
\implies\qquad
x_0<\frac 12,>\frac12,=
	\frac12\,.
\end{equation}
	Squaring we get 
	$$
			a_1^2=\frac{9(1+4x_0^2-4x_0)}{16(1-x_0)x_0}
	$$
namely
	$$
		(36+16a_1^2) x_0^2-(36+16a_1^2) x_0 +9=0\,.
	$$
	The solutions of the above second order equation
	are
	$$
		x_0=\frac 12\pm\frac{a_1}{\sqrt{9+4a_1^2}}\,,
	$$	
	but by \eqref{venezia} we have to choose the
	minus sign.
	Since by \eqref{verona}
	we have
	$$
	\frac{1}{4\sqrt{1-x_0}\sqrt{x_0}}
	=\frac{a_1}{3(1-2x_0)}
	$$
	by \eqref{mantova}
	and denoting for brevity $s:=\sqrt{4 a_1^2+9}$, we 
	get\footnote{Note that 
	$2x_0-1=-2a_1/\sqrt{4a_1^2+9}=-2a_1/s$.}
	\begin{eqnarray}
a_2&=&b''(x_0)=\frac{8 x_0^2-4 x_0-1}{x_0}
\frac{a_1}{3(1-2x_0)}
=
\frac{2(2x_0-1)^2+2(2x_0-1)-1}{x_0}
\frac{a_1}{3(1-2x_0)}
\nonumber
\\
&=&
\frac{8a_1^2-4 a_1 s-s^2}{3(s-2a_1)}=
\frac{8a_1^2-4 a_1 s-s^2}{27}
(s+2a_1)
\nonumber
\\
&=&
\frac{1}{27}\textstyle
(4 a_1^2-4a_1\sqrt{4a_1^2+9}-9)
(2a_1+\sqrt{4a_1^2+9})=:h(a_1)\,.
\label{acca}
\end{eqnarray}
Note that $h(a_1)=-g(-a_1)<-a_1$.
Since we are in the case $a_2<-a_1$ and
 we have proved that  the line 
$h(a_1)x+a_1$ is tangent to $b'(x)$,
we have that for $a_2<h(a_1)$
there are not intersections (zone $Z_{10}$) while for 	
$h(a_1)<a_2<-a_1$
there are two intersections (zone $Z_{12}$),
that we call $0<x_1^{(\pi)}<x_2^{(\pi)}<1$.
\\
In this last case, the function
  $x\to F(\pi,x)=a(x)-b(x)$ with $x\in(0,1)$ has 
  two critical points, which are exactly $x_1^{(\pi)}$
  and
 $ x_2^{(\pi)}$.
  Since $\lim_{x\to 0^+}\partial_x F(\pi,x)
  =\lim_{x\to 0^+} a'(x)-b'(x)
  =-\infty$
  and $\lim_{x\to 1^-}\partial_x F(\pi,x)
  =\lim_{x\to 1^-} a'(x)-b'(x)
  =a_2+a_1<0$,
  $x_1^{(\pi)}$ must be a minimum 
  and
 $ x_2^{(\pi)}$	a maximum.
\\	
Finally the case
of equation \eqref{gingerina1}
and the critical points of the function
$F(0,x)$
can be studied in the same way
sending $a_2\to-a_2$ and 
$a_1\to-a_1$.
\qed

%
%


\begin{thebibliography}{99}

\bibitem[B20]{B20} Bukhari M., Barry O. \textit{Spectro-spatial analyses of a nonlinear metamaterial with multiple nonlinear local resonators}, Nonlinear Dynamics, 99, pp. 1539--1560, 2020.

\bibitem[F22]{F22} Fortunati A., Bacigalupo A., Lepidi M., Arena A., Lacarbonara W. \textit{Nonlinear wave propagation in locally dissipative metamaterials via Hamiltonian perturbation approach}, Nonlinear Dynamics 108, n.2, pp.765--787, 2022.

\bibitem[M23]{M23} Murer M., Guruva S. K., Formica G., Lacarbonara W. 
\textit{A multi-bandgap metamaterial with multi-frequency resonators}, Journal of Composite Materials  57(4), 783-804 (2023).

\bibitem[SW23mssp]{SW23mssp} Shen Y., Lacarbonara Y. \textit{Nonlinear dispersion properties of metamaterial beams hosting nonlinear resonators and stop band optimization}, Mechanical Systems and Signal Processing 187, 2023.

\bibitem[SW23jsv]{SW23jsv} Shen Y., Lacarbonara W. \textit{Nonlinearity-enhanced wave stop bands in honeycombs
embedding spider web-like resonators}, Journal of Sound and Vibration 562, 2023.
 
\bibitem[Guo22]{Guo22} Wenjie G., Zhou Yang, Qingsong Feng, Chengxin Dai, Jian Yang, Xiaoyan Lei \textit{A new method for band gap analysis of periodic structures using virtual spring model and energy functional variational principle}, Mechanical Systems and Signal Processing 168,  2022.

\bibitem[Liu21]{Liu21} Liu Lei, Sridhar A., Geers M.G.D., Kouznetsova V.G. \textit{Computational homogenization of locally resonant acoustic metamaterial panels towards enriched continuum beam/shell structures}, Computer Methods in Applied Mechanics and Engineering 387,  2021.

\bibitem[Cai22]{Cai22} Cai Changqi, Zhou Jiaxi, Wang Kai, Pan Hongbin, Tan Dongguo, Xu Daolin, Wen Guilin \textit{Flexural wave attenuation by metamaterial beam with compliant quasi-zero-stiffness resonators}, Mechanical Systems and Signal Processing 174,  2022.

\bibitem[B16]{B16} Bacigalupo A., Gambarotta L. \textit{Simplified modelling of chiral lattice materials with local resonators}, International Journal of Solids and Structures 83, 126--141, 2016.

\bibitem[Comi18]{Comi18} Comi C., Driemeier L. \textit{Wave propagation in cellular locally resonant metamaterials}, Latin American Journal of Solids and Structures 15,  2018.

\bibitem[M22]{M22} Miranda Jr. E.J.P., Rodrigues S.F., Aranas Jr. C., Dos Santos, J.M.C. \textit{Plane wave expansion and extended plane wave expansion formulations for Mindlin-Reissner elastic metamaterial thick plates}, Journal of Mathematical Analysis and Applications 2, 505, 2022.

\bibitem[Fan21]{Fan21} Fan Lei, He Ye, Chen Xiao-an, Zhao Xue \textit{A frequency response function-based optimization for metamaterial beams considering both location and mass distributions of local resonators}, Journal of Applied Physics 11, 130, 2021.

\bibitem[Wang21]{Wang21} Wang Qiang, Li Jinqiang, Zhang Yao, Xue Yu, Li Fengming \textit{A frequency response function-based optimization for metamaterial beams considering both location and mass distributions of local resonators}, Mechanical Systems and Signal Processing 151, 2021.

\bibitem[CP23]{CP23} Ch\`{a}vez-Pichardo M., Mart\`{i}nez-Cruz M.A., Trejo-Mart\`{i}nez A., Vega-Cruz A.B., Arenas-Resendiz T. \textit{On the Practicality of the Analytical Solutions for all Third- and Fourth-Degree Algebraic Equations with Real Coefficients}
Mathematics 11, 1147, 2023.

\bibitem[W]{W} Lacarbonara W. \textit{Nonlinear Structural Mechanics: Theory, Dynamical Phenomena and Modeling}, Springer, New-York, 2013.

\bibitem[Elliptic] {Elliptic} Takebe T. \textit{Elliptic Integrals and Elliptic Functions}, Moscow Lectures, Springer, 2022.

\bibitem[V07]{V07} Sanders J.A., Verhulst F., Murdock J. \textit{Averaging Methods in Nonlinear Dynamical Systems, Revised 2nd Edition}, Springer, New York, 2007.

\bibitem[L19jsv]{L19jsv} Fronk M. D., Leamy M. J. \textit{Direction-dependent invariant waveforms and stability in two-dimensional, weakly nonlinear lattices}, Journal of Sound and Vibration 447, pp. 137--154, 2019.

\bibitem[M15]{M15} Malek S., Gibson L. \textit{Effective elastic properties of periodic hexagonal honeycombs}, Mechanics of Materials 91, pp. 226--240, 2015.

\bibitem[S18]{S18} Sorohan S., Constantinescu D.M.,  Sandu M., Sandu A.G. \textit{On the homogenization of hexagonal honeycombs under axial and shear loading. Part I: Analytical formulation for free skin effect}, Mechanics of Materials 119, pp. 74--91, 2018.

\bibitem[G97]{G97} Gibson L.J., Ashby M.F. \textit{Cellular solids: structure and properties}, Cambridge Solid State Science Series, Cambridge University Press, 1997.

\bibitem[Graff]{Graff}
Graff S.M.
\textit{On the conservation of hyperbolic invariant tori for Hamiltonian systems},
J. Differential Equations 15, 1-69, 1974.

\bibitem[Val]{Val}
Valdinoci, E.
\textit{Families of whiskered tori for a-priori stable/unstable Hamiltonian systems and construction of unstable orbits},
Math. Phys. Electron. J. 6, Paper 2, 31 pp., 2000.


\bibitem[MNT]{MNT}  Medvedev A.G., Neishtadt A.I., Treschev D.V. \textit{Lagrangian tori near resonances of near--integrable Hamiltonian systems}, 
Nonlinearity 28 (7), pp. 2105--2130, 2015.


\bibitem[DL]{DL} Di Gregorio L., Lacarbonara W. \textit{On bandgaps sensitivity to 3:1 interactions between acoustic and optical waves}, Preprint 2024.

\bibitem[H16]{H16} Haller G., Ponsioen S. \textit{Nonlinear normal modes and spectral submanifolds: existence, uniqueness and use in model reduction}, Nonlinear Dynamics 86, pp. 1493--1534, 2016.

\bibitem[Cabre05]{Cabre05} Cabre X., Fontich E., de la Llave R. \textit{The parametrization method for invariant manifolds III: overview and applications}, J. Differential Equations 218, pp. 444--515, 2005.

\bibitem[Celletti13]{Celletti13}  Calleja R.C., Celletti A., de la Llave R. \textit{A KAM theory for conformally symplectic systems: Efficient algorithms and their validation}, J. Differential Equations 255, pp. 978--1049, 2013.

\bibitem[Llave05]{Llave05} Haro A., de la Llave R. \textit{A parameterization method for the computation of invariant tori and their whiskers in quasi-periodic maps: Rigorous results}, J. Differential Equations 228, pp. 230--279, 2005.

\bibitem[Fontich23]{Fontich23} Fontich E., Vierio A. \textit{Dynamics near the invariant manifolds after a Hamiltonian-Hopf bifurcation}, Communications in Nonlinear Science and Numerical Simulation 117, 2023, 106971.

\bibitem[Llave06]{Llave06} Haro A., de la Llave R. \textit{A parameterization method for the computation of invariant tori and their whiskers in quasi-periodic maps: numerical algorithms}, Discrete and continuous dynamical systems Series B, 6 (6), pp. 1261--1300, 2006.

\bibitem[HW96]{HW96} Haller G., Wiggins S. \textit{Geometry and chaos near resonant equilibria of 3-DOF Hamiltonian systems}, Physica D 90, pp. 319--365, 1996.

\bibitem[HW95]{HW95} Haller G., Wiggins S. \textit{N-pulse homoclinic orbits in perturbations of resonant Hamiltonian systems}, Arch. Rat. Mech. Anal. 130, pp. 25--101, 1995.

\bibitem[HW93]{HW93} Haller G., Wiggins S. \textit{Orbits homoclinic to resonances: the Hamiltonian case}, Physica D 66, pp. 298--346, 1993.

\end{thebibliography}
\end{document}